\documentclass[a4paper,10pt, oneside]{article}
\usepackage{bm,amsmath,amssymb,amsthm,mathrsfs,graphicx}
\usepackage[font=small,labelfont=md,textfont=it]{caption}
\usepackage{subfigure}
\usepackage[pdftex, pdfpagelabels, bookmarks, hyperindex, hyperfigures,
colorlinks=red, pagebackref]{hyperref}
\usepackage{diagbox}
\usepackage{amssymb}
\usepackage{amsmath,bm}
\usepackage{multirow}
\usepackage{makecell}
\usepackage{graphicx}

\usepackage{float}
\usepackage[titletoc, title]{appendix}
\usepackage{color}
\usepackage{lipsum}
\usepackage{amsfonts}
\usepackage{graphicx}
\usepackage{epstopdf}
\usepackage{algorithmic}
\usepackage{amsopn}
\usepackage{amsthm}
\usepackage{cleveref}
\usepackage{amssymb}

\numberwithin{equation}{section}
\numberwithin{figure}{section}

\newtheorem{lemma}{Lemma}[section]
\newtheorem{theorem}{Theorem}[section]
\newtheorem{remark}{Remark}[section]
\newtheorem{definition}{Definition}[section]

\newtheorem{exam}{Example}[section]

\begin{document}

\title{\Large
Robust globally divergence-free weak Galerkin methods for   stationary incompressible convective Brinkman-Forchheimer equations
\thanks
{
	This work was supported in part by National Natural Science Foundation of China  (12171340).
}
}
	\author{Xiaojuan Wang \thanks{ Email: wangxiaojuan@scu.edu.cn }, \
		Xiaoping Xie \thanks{Corresponding author. 
			Email: xpxie@scu.edu.cn}	\\
			School of Mathematics, Sichuan University, Chengdu 610064, China	
	}
\maketitle

\begin{abstract}
This paper  develops a class of  robust weak  Galerkin   methods for the stationary incompressible convective Brinkman-Forchheimer   equations. The methods adopt piecewise polynomials of degrees  $m\ (m\geq1)$ and $m-1$  respectively  for the approximations of   velocity  and pressure variables inside the elements and piecewise polynomials of degrees  $k \ ( k=m-1,m)$ and $m$  respectively   for    their numerical traces   on the interfaces of elements, and are   shown  to   yield globally divergence-free   velocity approximation. Existence and uniqueness results  for the discrete schemes, as well as optimal a priori error estimates, are established.  A convergent linearized iterative algorithm is also presented. Numerical experiments are provided to verify the performance of the proposed    methods.
\end{abstract}


\text{\bf AMS 2010} \; 65M60, 65N30


%
\section{Introduction}
Let $\Omega\subset \mathbb{R}^n$ $(n=2,3)$ be a Lipschitz  polygonal/polyhedral domain.  We consider the following stationary incompressible convective Brinkman-Forchheimer model: 
\begin{eqnarray}\label{BF0}
\left\{
\begin{aligned}
-\nu \Delta \bm{u}+\nabla\cdot(\bm{u}\otimes\bm{u})+\alpha |\bm{u}|^{r-2}\bm{u}+\nabla p&=\bm{f},&\text{in} \ \Omega,\label{BF0}\\
\nabla\cdot\bm{u}&=0,&\text{in} \ \Omega,\\
\bm{u}&=\bm{0},&\text{on} \ \partial \Omega. \\
\end{aligned}
\right.
\end{eqnarray}
 Here  $\bm{u}=(u_1,u_2,\cdots,u_n)^T$ is the velocity vector, $p$     the pressure,  $\bm{f}$ a given forcing function,  $\nu$    the  Brinkman coefficient,  $\alpha>0$   the    Forchheimer coefficient, and $r\geq 2$. The operator $\otimes$ is defined by $\bm{u}\otimes\bm{v}=(u_iv_j)_{n\times n}$   for $\bm{v}=(v_1,v_2,\cdots,v_n)^T$.

The Brinkman-Forchheimer model, which  can be viewed as the Navier-Stokes equations with a nonlinear damping term, is  used to modelling fast flows in highly porous media \cite{D1982Nonlinear,K1981Boundary}.
 In recent years there have developed many numerical algorithms   for  Brinkman-Forchheimer equations,  such as  conforming mixed finite element  methods \cite{CaucaoSergio2022AtBs,Caucao;2021,2022Unconditional,2019Mixed,YangHuaijun2023Saot},
  nonconforming mixed finite element  methods \cite{QianLiu2021Saoa},
  stabilized mixed methods \cite{LiZhenzhen2019Smfe,Louaked;2017}, multi-level mixed methods \cite{LiMinghao2019Tmfe,QiuHailong2019Msaf,ZhengBo2021Tdsa,ZhengBo2023Tsaf}, parallel finite element algorithms \cite{WassimEid2022Lapf,WassimEid2023Aptm}.
 We refer to  \cite{CaiXiaojing2008Wass,  CelebiA.O.2006Ocdo, Djoko;2014, JiangZaihong2012Abos,LiKwang-Ok2022Grft,LiYuanfei2014Cdft,LiuHui2021Wotg,PayneL.E.1999CaCD,YangRong2021Aso3,ZhangZujin2011Otuo,ZhongXin2019Anot,ZhouYong2012Rauf} for the study of the  properties of  weak/strong solutions to the Brinkman-Forchheimer equations.

It is well-known that  the  divergence constraint  $\nabla\cdot\bm{u}=0$ corresponds to the conservation of   mass for  incompressible fluid flows, and   that   numerical methods   with poor conservation usually  suffer from instabilities   \cite{bookin1, JLMNR2017, add4,bookin3,  add2}.
Besides, the  numerical schemes  with   exactly divergence-free velocity  approximation may automatically lead to pressure-robustness in the sense that the velocity approximation  error   is independent of    the pressure approximation \cite{JLMNR2017,MR3564690,MuLin2023ApwG}.
We refer to \cite{CFX2016,Chen-Xie2023,CKS2007,HLX2019,HX2019,MR3511719,MuLin2018Addf,XZ2010,ZhangLi2019Agdw,ZCX2017}
 for some divergence-free finite element methods for  the incompressible fluid flows.
	
	In this paper we consider a robust globally divergence-free weak Galerkin finite element discretization of the Brinkman-Forchheimer model \eqref{BF0}. The WG framework was first proposed  in \cite{WangJunping2013AwGf,WangJunping2013AWGM}  for  second-order elliptic problems.  It  allows the use of totally discontinuous functions on meshes with arbitrary shape of polygons/polyhedra due to the introduction of weakly defined gradient/divergence operators over functions with discontinuity, and has  the local elimination property, i.e. the unknowns defined in the interior of elements can be locally eliminated by using the numerical traces  defined on the interfaces of  elements. We refer to \cite{CFX2016, Chen-Xie2023, HLX2019,HX2019,HuXiaozhe2019AwGf,MuLin2023ApwG,MuLin2018Addf,PengHui2022WGmf,WangJunping2016AwGf,WangRuishu2016AwGf, XZ2010,ZZX2023,ZCX2017}  some developments and applications of  WG methods for fluid flow  problems.  Particularly,  a class of robust globally divergence-free weak Galerkin methods were developed  in \cite{CFX2016} for Stokes equations, and later were extended  to solve incompressible quasi-Newtonian Stokes equations \cite{ZCX2017}, natural convection equations \cite{HLX2019,HX2019} and incompressible Magnetohydrodynamics flow equations \cite{ZZX2023}.

 The goal of this  contribution is to extend the WG methods of  \cite{CFX2016} to the discretization of  the Brinkman-Forchheimer model.
The main features of our  WG discretization  for the   model \eqref{BF0} are as follows:
\begin{itemize}
\item  The   discretization scheme is arbitrary order, which adopts   piecewise polynomials of degrees $m$ $(m\geq1)$ and $m-1$ to approximate the velocity and pressure inside the elements, respectively, and piecewise polynomials of degrees $k$ $(k=m-1,m)$ and $m$ to approximate  the traces of velocity and pressure on the interfaces of elements, respectively.

\item The scheme yields     globally divergence-free velocity approximation, which  automatically leads to   pressure-robustness.	

\item The scheme is    ``parameter-friendly", i.e.   the stabilization parameter in the scheme
does not require  to be  ``sufficiently large".

\item The unknowns of the velocity and pressure in the interior of  elements can be locally eliminated so as to obtain a reduced discrete system of smaller size.

\item The well-posedness  and optimal error estimates of the scheme  are established.

\end{itemize}
		
The rest of this paper is organized as follows.
Section 2 gives notations, weak formulations, the WG  scheme and some preliminary results.
Section 3  establishes  the well-posedness of the discrete scheme.
Sections 4 is devoted to the a priori error analysis.
Section 5 derives $L^2$ error estimate for the velocity.
Section 6 shows the local elimination property and proposes an iteration algorithm for the nonlinear WG scheme. Section 7 provides several numerical experiments.  Finally, Section 8 gives some concluding remarks.

\section{Weak Galerkin finite element scheme}
\subsection{Notation and weak problem}

For any bounded domain $\Lambda \subset \mathbb{R}^l$ $( l=n,n-1 )$, nonnegative  integer $s$ and real number $1\leq q<\infty$, let $W^{s,q}(\Lambda)$ and $W_0^{s,q}(\Lambda)$ be the usual Sobolev spaces defined on $\Lambda$ with norm $||\cdot||_{s,q, \Lambda}$ and semi-norm $|\cdot|_{s,q,\Lambda}$. In particular,   $H^s(\Lambda):=W^{s,2}(\Lambda)$ and $H^s_0(\Lambda):=W^{m,2}_0(\Lambda)$, with $||\cdot||_{s,\Lambda}:=||\cdot||_{s,2, \Lambda}$ and $|\cdot|_{s,\Lambda}:=|\cdot|_{s,2, \Lambda}$.   We use $(\cdot,\cdot)_{s,\Lambda}$ to denote the inner product of $H^s(\Lambda)$, with $(\cdot,\cdot)_\Lambda :=(\cdot,\cdot)_{0,\Lambda}$.
	When $\Lambda = \Omega$, we set $||\cdot||_s := ||\cdot||_{s,\Omega},|\cdot|_s := |\cdot|_{s,\Omega}$,
	and $(\cdot,\cdot):= (\cdot,\cdot)_{\Omega}$.
	Especially, when
	$\Lambda \subset \mathbb{R}^{n-1}$ we use $\langle\cdot,\cdot\rangle_\Lambda$ to replace $(\cdot,\cdot)_{\Lambda}$. For a nonnegative integer $m$, let $P_m(\Lambda)$  be the set of all polynomials defined  on $\Lambda$ with degree   no more than $m$. We also need the following Sobolev spaces:
\begin{eqnarray*}
	L^2_0(\Omega):=\{q\in L^2(\Omega):(q,1)=0\},
\end{eqnarray*}
$$ \bm{H}({\rm div};\Lambda):=\left\{\bm v\in [L^2(\Lambda)]^n:\nabla\cdot \bm{v}\in L^2(\Lambda)\right\}.$$

Let $\mathcal{T}_h$ be a shape regular  partition of $\Omega$ into closed simplexes,
and let $\mathcal{E}_h$ be the set of all edges (faces) of all the elements in $\Omega$.
	For any $K\in \mathcal{T}_h$, $e\in \mathcal{E}_h$, we denote by $h_K$ the diameter of $K$
	and by $h_e$ the diameter of $e$, and set $h=\max_{K\in \mathcal{T}_h}h_K$.
	Let $\bm{n}_K$  and $\bm{n}_e$ denote the outward unit normal vectors along the boundary $\partial K$ and $e$, respectively. We may abbreviate $\bm{n}_K$ as $\bm{n}$ when  there is no ambiguity.
	We use  $\nabla_h$ and $\nabla_h\cdot$  to denote  respectively  the operators of piecewise-defined gradient and divergence with respect to the decomposition $\mathcal{T}_h$.
	
For convenience, throughout the paper we use $x\lesssim y$ ($x\gtrsim y$)
to denote $x\le Cy$ ($x\ge Cy$), where   $C$  is a positive constant  independent of the mesh size $h$.

We introduce the spaces
\begin{eqnarray}
\bm{V}:=[H^1_0(\Omega)]^n,\quad
Q:=L^2_0(\Omega),\quad
\bm{V}_0:=\{\bm{v}\in\bm{V}:\nabla
\cdot \bm{v}=0\},\nonumber
\end{eqnarray}
and define   the following bilinear and trilinear forms: for $\bm{u}, \bm{v}\in \bm{V}$ and $q\in Q$,
\begin{align*}
&a(\bm{u},\bm{v}):=\nu ( \nabla \bm{u},\nabla \bm{v}), \quad b(\bm{v},q):=-( q, \nabla\cdot\bm{v}),\\
&c(\bm{\kappa};\bm{u},\bm{v}):= \alpha(|\bm{\kappa}|^{r-2}\bm{u},\bm{v}),\\
&d(\bm{\kappa};\bm{u},\bm{v}):=\frac{1}{2}(\nabla\cdot(\bm{u}\otimes\bm{\kappa}),\bm{v})
-\frac{1}{2}(\nabla\cdot(\bm{v}\otimes\bm{\kappa}),\bm{u}).
\end{align*}

Then  the weak form of \eqref{BF0} is given as follows: seek $(\bm{u},p )\in \bm{V}\times Q$ such that
\begin{subequations}\label{weak}
\begin{align}
a(\bm{u},\bm{v})+b(\bm{v},p )+c(\bm{u};\bm{u},\bm{v} )+d(\bm{u};\bm{u},\bm{v} )&=(f,\bm{v}),&\forall\bm{v} \in \bm{V},\\
b(\bm{u},q )&=0, &\forall q \in Q.
\end{align}
\end{subequations}

\begin{remark} As shown in \cite{2019Mixed},  
 the weak problem \eqref{weak} admits at least one solution $(\bm{u},p )\in \bm{V}\times Q$  when  $\Omega$ is a bounded Lipschitz domain and   $\bm{f}\in [H^{-1}(\Omega)]^n$, and there holds
 \begin{align}\label{boundedCondition}
 \|\nabla\bm{u}\|_{0}\leq \frac{\|\bm{f}\|_*}{\nu}.
 \end{align}
 In addition, if the smallness  condition
\begin{align}\label{Condition}
 \frac{\mathcal{N}\|\bm{f}\|_*}{\nu^2} <1 
\end{align}
holds,   then the solution of \eqref{weak} is unique. Here
  \begin{align*}
\|\bm{f}\|_*:=\sup_{\bm{0}\neq\bm{v}\in\bm{V}_0}\frac{(\bm{f},\bm{v})}
{\|\nabla\bm{v}\|_0}, \quad 
\mathcal{N}:=\sup_{\bm 0\neq\bm{u},\bm{v},\bm{\kappa}\in \bm{V}_0}\frac{d(\bm{\kappa};\bm{u},\bm{v})}
{\|\nabla\bm{\kappa}\|_0\|\nabla\bm{u}\|_0\|\nabla\bm{v}\|_0}.
\end{align*}
\end{remark}

\subsection{WG scheme}
In order to give the WG scheme to the system \eqref{BF0} we introduce, for integer $\gamma\geq 0$, the
discrete gradient operator $\nabla_{w,\gamma}$ and  the discrete weak divergence operator $\nabla_{w,\gamma}\cdot$ as follows.

\begin{definition}
For all $K\in \mathcal{T}_{h}$ and $v\in \mathcal{V}(K):=\{v=\{v_i,v_b\}:v_i\in L^2(K),v_b\in H^{1/2}(\partial K)\}$, the discrete weak gradient   $\nabla_{w,\gamma,K}v\in [P_{\gamma}(K)]^{n}$ of v on $K$ is defined by
\begin{align}
(\nabla_{w,\gamma,K}v, \bm{\varsigma})_K=-(v_i,\nabla \cdot \bm{\varsigma})_K+\langle v_b,\bm{\varsigma}\cdot \bm{n}_{K} \rangle_{\partial K}, \quad\forall  \bm{\varsigma}\in [P_{\gamma}(K)]^{n}.
\end{align}
Then  the global discrete weak gradient operator $\nabla_{w,\gamma}$ is defined as
$$\nabla_{w,\gamma}|_{K}:=\nabla_{w,\gamma,K},\quad \forall K\in \mathcal{T}_h. $$
Moreover, for a vector  $\bm{v}=(v_1,v_2,...,v_n)^{T}$ with $v_j|_K\in \mathcal{V}(K)$ for $j=1,...,n$,
the discrete weak gradient $\nabla_{w,\gamma}\bm{v}$ is defined as
$$ \nabla_{w,\gamma}\bm{v}:= (\nabla_{w,\gamma}v_1,\nabla_{w,\gamma}v_2,..., \nabla_{w,\gamma}v_n )^{T}.$$
\end{definition}

\begin{definition}
For all  $K\in \mathcal{T}_{h}$ and $\bm{w}\in \bm{\mathcal{W}}(K):=\{\bm{w}=\{\bm{w}_i,\bm{w}_b\}:\bm{w}_i\in [L^2(K)]^{n},\bm{w}_b\cdot\bm{n}_{K}\in H^{-1/2}(\partial K)\} $, the discrete weak divergence   $\nabla_{w,\gamma,K}\cdot\bm{w}\in P_{\gamma}(K)$ of $\bm{w}$ on $K$ is defined by
\begin{align}
(\nabla_{w,\gamma,K}\cdot\bm{w} , \varsigma)_K=-(\bm{w}_i,\nabla  \varsigma)_K+\langle \bm{w}_b\cdot \bm{n},\varsigma \rangle_{\partial K}, \quad\forall \varsigma\in P_{\gamma}(K).
\end{align}
Then the global discrete weak divergence operator $\nabla_{w,\gamma}\cdot $ is defined as $$\nabla_{w,\gamma}\cdot|_{K}:=\nabla_{w,\gamma,K}\cdot, \quad\forall K\in \mathcal{T}_h. $$
Moreover, for a tensor  $\widetilde{\bm{w}}=(\bm{w}_1,...,\bm{w}_n)^{T}$ with $\bm{w}_j|_K\in  \bm{\mathcal{W}}(K)$ for $j=1,...,n$, the discrete weak divergence $\nabla_{w,\gamma}\cdot \widetilde{\bm{w}}$ is defined as $$\nabla_{w,\gamma}\cdot \widetilde{\bm{w}}:=(\nabla_{w,\gamma}\cdot \bm{w}_{1},...,\nabla_{w,\gamma}\cdot \bm{w}_n)^{T}. $$
\end{definition}

For any $K\in\mathcal{T}_h$, $ e\in\mathcal{E}_h$ and nonnegative integer $j$, 
let $\Pi_j^\ast: L^2(K)\rightarrow P_{j}(K)$ and $\Pi_j^{B}: L^2(e)\rightarrow P_{j}(e)$ be the usual $L^2$-projection operators. We shall adopt $\bm{\Pi}_j^\ast$ to denote $\Pi_j^\ast$ for the vector form.

For any integer $m\ge 1$, and integer $k= m-1,m $,
we introduce the following finite dimentional spaces:
\begin{align*}
\bm{V}_h:=&\{\bm{v}_h=\{\bm{v}_{hi},\bm{v}_{hb}\}: \bm{v}_{hi}|_K\in [P_m(K)]^n,\bm{v}_{hb}|_e\in [P_{k}(e)]^n ,
\forall K\in\mathcal{T}_h,\forall e\in\mathcal{E}_h\},\\
\bm{V}_h^{0}:=&\{\bm{v}_h=\{\bm{v}_{hi},\bm{v}_{hb}\}\in\bm{V}_h: \bm{v}_{hb}|_{\partial\Omega}=\bm{0} \},\\
Q_h:=&\{q_h=\{q_{hi},q_{hb}\}: q_{hi}|_K\in P_{m-1}(K),q_{hb}|_e\in P_{m}(e),
\forall K\in\mathcal{T}_h,\forall e\in\mathcal{E}_h\},\\
Q_h^{0}:=&\{q_h=\{q_{hi},q_{hb}\}\in Q_h, q_{hi}\in L_0^2(\Omega)\}.
\end{align*}

For any $\bm{u}_{h}=\{\bm{u}_{hi},\bm{u}_{hb}\},\bm{v}_{h}=\{\bm{v}_{hi},\bm{v}_{hb}\},
\bm{\kappa}_{h}=\{\bm{\kappa}_{hi},\bm{\kappa}_{hb}\}\in \bm{V}_{h}^{0}$, and $p_h=\{p_{hi}, p_{hb}\}\in Q_h^{0},$
we shall define bilinear and trilinear terms as follows:
\begin{align*}
a_h(\bm{u}_h,\bm{v}_h)&:=\nu(\nabla _{w,m-1}\bm{u}_h,\nabla _{w,m-1}\bm{v}_h)
+s_h(\bm{u}_h,\bm{v}_h), \\
s_h(\bm{u}_h,\bm{v}_h)&:=
\nu \langle  \eta( \bm{\Pi}_{k}^{B}\bm{u}_{hi}-\bm{u}_{hb}),\bm{\Pi}_{k}^{B} \bm{v}_{hi}-\bm{v}_{hb} \rangle_{\partial \mathcal{T}_h},\\
b_h(\bm{v}_h,q_h)&:=(\nabla _{w,m}q_h,\bm{v}_{hi}), \\
c_h(\bm{\kappa}_h ;\bm{u}_h,\bm{v}_h )&:=\alpha( |\bm{\kappa}_{hi}|^{r-2}\bm{u}_{hi}, \bm{v}_{hi} ), \\
d_h(\bm{\kappa}_h ;\bm{u}_h,\bm{v}_h )
&:=\frac{1}{2}( \nabla _{w,m}\cdot\{\bm{u}_{hi}\otimes\bm{\kappa}_{hi}, \bm{u}_{hb}\otimes\bm{\kappa}_{hb}\},\bm{v}_{hi}  )\\
&\qquad
    -\frac{1}{2}( \nabla _{w,m}\cdot\{\bm{v}_{hi}\otimes\bm{\kappa}_{hi}, \bm{v}_{hb}\otimes\bm{\kappa}_{hb}\},\bm{u}_{hi}  ),
\end{align*}
where $ \langle \cdot,\cdot\rangle_{\partial\mathcal{T}_h}:=\sum_{K\in \mathcal{T}_h}\langle \cdot,\cdot\rangle_{\partial K}$, and
the stabilization parameter $\eta|_{\partial K}=h_K^{-1} ,$ $\forall K\in\mathcal{T}_h$.

In what follows we assume that $$\bm{f}\in \bm{L}^2(\Omega).$$
Based on the above definitions, the WG scheme for \eqref{BF0} reads:
seek $\bm{u}_h=\{\bm{u}_{hi},\bm{u}_{hb} \}\in \bm{V}_h^{0}$, $p_h=\{p_{hi}, p_{hb}\}\in Q_h^{0} $ such that
\begin{subequations}\label{WG}
		\begin{align}
		a_h(\bm{u}_h,\bm{v}_h)+b_h(\bm{v}_{h},p_h)+c_h(\bm{u}_h ;\bm{u}_h,\bm{v}_h )
		+d_h(\bm{u}_h ;\bm{u}_h,\bm{v}_h )
         =&(\bm{f},\bm{v}_{hi}), \ \forall \bm{v}_h \in \bm{V}_h^{0},\label{WG1}\\
		b_h(\bm{u}_{h},q_h)=&0, \ \forall q_h\in  Q_h^{0}. \label{WG2}
		\end{align}
	\end{subequations}

The following theorem shows that the scheme \eqref{WG} yields   globally divergence-free velocity approximation.
\begin{theorem}\label{TH2.2}
Let $\bm{u}_h=\{\bm{u}_{hi},\bm{u}_{hb} \}\in \bm{V}_h^{0}$ be the velocity solution of the WG scheme \eqref{WG}. Then there hold
 \begin{align}\label{divergence-free}
	&\bm{u}_{hi}\in \bm{H}({\rm div};\Omega),\quad  \nabla\cdot\bm{u}_{hi}=0.
	\end{align}
\end{theorem}
\begin{proof}
 Define  a function $\varphi_{hb}\in L^{2}(\mathcal{E}_{h})$ as follows: for any $e\in \mathcal{E}_{h}$,
\begin{eqnarray*}
\varphi_{hb}|_e:=
\left\{
\begin{aligned}
-\big((\bm{u}_{hi}\cdot \bm{n}_e)|_{K_1}\big)|_e-\big((\bm{u}_{hi}\cdot \bm{n}_e)|_{K_2}\big)|_e,\quad &\text{if } e=K_1\cap K_2, K_1,K_2\in \mathcal{T}_h,\\
0,\quad& \ \forall e\subset \partial\Omega.
\end{aligned}
\right.
\end{eqnarray*}
Let
$\varphi_{0}:=\frac{1}{|\Omega|}\int_\Omega \nabla_{h}\cdot \bm{u}_{hi} d\bm{ x}. $
Taking  $q_{hi}=\nabla_h\cdot\bm{u}_{hi} -\varphi_{0},q_{hb} =  \varphi_{hb}-\varphi_{0}  $ in \eqref{WG2}, we obtain
\begin{align*}
0=&-(\bm{u}_{hi},\nabla_{w,m}q_h )\nonumber\\
=&(\nabla_{h}\cdot \bm{u}_{hi},q_{hi} )-\sum_{K\in \mathcal{T}_h}\langle \bm{u}_{hi}\cdot\bm{n},q_{hb} \rangle_{\partial K}\nonumber\\
=&(\nabla_{h}\cdot \bm{u}_{hi},\nabla_h\cdot\bm{u}_{hi} -\varphi_{0})
 -\sum_{K\in \mathcal{T}_h}\langle \bm{u}_{hi}\cdot\bm{n},\varphi_{hb}-\varphi_{0} \rangle_{\partial K}\nonumber\\
 =&(\nabla_{h}\cdot \bm{u}_{hi},\nabla_h\cdot\bm{u}_{hi})
   -\sum_{K\in \mathcal{T}_h}\langle \bm{u}_{hi}\cdot\bm{n},\varphi_{hb} \rangle_{\partial K}  \nonumber\\
  =&\|\nabla_{h}\cdot \bm{u}_{hi} \|^{2}_{0}+\sum_{e\in \mathcal{E}_{h}, e\nsubseteq\partial\Omega}\|(\bm{u}_{hi}\cdot \bm{n}_e)|_{K_1}+(\bm{u}_{hi}\cdot \bm{n}_e)|_{K_2} \|_{0,e}^{2},
\end{align*}
which indicates  the desired conclusion \eqref{divergence-free}. 
\end{proof}

\subsection{Preliminary results}

We first introduce   two semi-norms $|||\cdot|||_V$ and $|||\cdot|||_Q$  on the spaces $ \bm{V}_h$ and $ Q_h$, respectively, as follows:
\begin{align*}
|||\bm{v}_h|||_V^2&:= \|\nabla _{w,m-1}\bm{v}_h\|^2_0
+\|\eta^{\frac 1 2}(\bm{\Pi}_{k}^{B}\bm{v}_{hi}-\bm{v}_{hb})\|^2_{0,\partial \mathcal{T}_h}, \quad \forall \bm{v}_h\in \bm{V}_h,\\
|||q_h|||_Q^2&:=\|{q}_{hi}\|^2_0+\sum_{K\in \mathcal{T}_h}h_{K}^{2}\|\nabla _{w,m} q_h\|_{0,K}^{2}, \quad \forall q_h \in Q_h,
\end{align*}
where $\|\cdot\|_{0,\partial \mathcal{T}_h}:=(\sum_{K\in\mathcal{T}_h}\|\cdot\|^2_{0,\partial K})^{1/2}$,
and we recall that $\eta|_{\partial K}=h_{ K}^{-1}$.
It is easy to see that $|||\cdot|||_V$ and $|||\cdot|||_Q$ are norms on $\bm{V}_h^{0}$ and $Q_h^{0}$, respectively(cf. \cite{  CFX2016}).

The following lemma follows from the trace theorem, the inverse inequality and scaling arguments    (cf. \cite{HX2019,bookFEM}).
\begin{lemma}\label{Lemma 2.1}
For all $ K\in\mathcal{T}_{h}, \omega\in H^{1}(K) $, and $ 1\leq q\leq \infty $, there holds
\begin{equation*}
\|\omega\|_{0,q,\partial K}\lesssim h_K^{-\frac{1}{q}}\|\omega\|_{0,q,K}+h_K^{1-\frac{1}{q}}|\omega|_{1,q,K}.
\end{equation*}
In particular, for all $ \omega\in P_{m}(K) $,
\begin{equation*}
\|\omega\|_{0,q,\partial K}\lesssim h_K^{-\frac{1}{q}}\|\omega\|_{0,q,K}.
\end{equation*}
\end{lemma}

For the projections $\Pi_j^\ast$ and $\Pi_j^B$ with $j\geq 0$, the following  approximation and stability  results are standard.
\begin{lemma}[cf. \cite{bookFEM}]\label{Lemma 2.2}
For $\forall K\in \mathcal{T}_{h}, \forall e\in \mathcal{E}_{h}$ and $1\leq l\leq j+1$, there hold 
\begin{align*}
|| \omega-\Pi_j^\ast \omega||_{0,K}+h_K|\omega-\Pi_j^\ast \omega|_{1,K}&\lesssim h_K^{l}|\omega|_{l,K},\quad\forall \omega\in H^{l}(K),\\
|| \omega-\Pi_j^\ast \omega||_{0,\partial K}+|| \omega-\Pi_j^B \omega||_{0,\partial K}&\lesssim h_K^{l-\frac{1}{2}}|\omega|_{l,K},\quad\forall \omega\in H^{l}(K),\\
\|\Pi_j^\ast \omega\|_{0,K}&\leq \|\omega\|_{0,K},\quad \forall v\in L^2(K),\\
\|\Pi_j^B \omega\|_{0,e}&\leq \|\omega\|_{0,e}, \quad\forall v\in L^2(e).
\end{align*}
\end{lemma}

In view of the definitions of the discrete weak gradient operator,  the Green's formula,
  the projection operator,  the Cauchy-Schwarz inequality, the inverse inequality and the trace inequality,   the following lemma holds (cf. \cite{CFX2016}).
\begin{lemma}\label{Lemma 2.3}
For any $K\in \mathcal{T}_h$ and $\bm{\omega}_{h}=\{\bm{\omega}_{hi},\bm{\omega}_{hb} \}\in [P_{m}(K)]^{n}\times[P_{k}(\partial K) ]^{n}$ with $0\leq m-1\leq k\leq m$, there hold
\begin{subequations}
\begin{align}
\|\nabla \bm{\omega}_{hi}\|_{0,K}\lesssim \|\nabla_{w,m-1}\bm{\omega}_{h}\|_{0,K}+h_K^{-\frac{1}{2}}\| \bm{\Pi}_{k}^{B}\bm{\omega}_{hi}-\bm{\omega}_{hb}\|_{0,\partial K},\label{2.10a}\\
\|\nabla_{w,m-1} \bm{\omega}_{h}\|_{0,K}\lesssim \|\nabla \bm{\omega}_{hi}\|_{0,K}+h_K^{-\frac{1}{2}}\| \bm{\Pi}_{k}^{B}\bm{\omega}_{hi}-\bm{\omega}_{hb}\|_{0,\partial K}.\label{2.10b}
\end{align}
\end{subequations}
\end{lemma}

By the definition of the norm  $|||\cdot|||_V$, we
further have the following conclusion (cf. \cite{HX2019,ZZX2023}):
\begin{lemma} \label{Lemma 2.4}
For any $ \bm{v}_{h}\in \bm{V}_{h}^{0}$, there hold
\begin{align}
 \|\nabla_{h} \bm{v}_{hi}\|_{0}\lesssim ||| \bm{v}_{h}|||_V
\end{align}
and
\begin{eqnarray}\label{Cr-est}
		\|\bm{v}_{hi}\|_{0,r}\leq C_{\widetilde{r}} |||\bm{v}_{h}||| _V
	\end{eqnarray}
	for   $r$ satisfying
	\begin{eqnarray*}
		\left\{
		\begin{aligned}
			&2\le r< \infty,&\text{ if } \ n=2,\\
			&2\le r\le 6 ,&\text{ if } \ n=3,
		\end{aligned}
		\right.
	\end{eqnarray*}
	where $C_{\widetilde{r}}>0$ is a positive constant only  depending on $r$.
\end{lemma}

For  any  integer $ j\geq 0$, we introduce the local Raviart-Thomas (RT) element space
\begin{eqnarray}
	\bm{RT}_j(K)=[P_j(K)]^n+\bm{x}P_j(K), \ \forall K\in \mathcal{T}_{h}\nonumber
\end{eqnarray}
and the RT  projection operator $\bm{P}^{RT}_j: [H^1(K)]^n\rightarrow  \bm{RT}_j(K)$ (cf. \cite{RT-1}) defined by
	\begin{align}
			\langle\bm{P}^{RT}_j\bm{\omega}\cdot \bm{n}_e,\sigma \rangle_e&=\langle \bm{\omega}\cdot \bm{n}_e,\sigma\rangle_e, \quad\forall \sigma\in P_j(e),  e\in\mathcal{E}_{h} \text{ and } e\subset \partial K, \text{ for } j\geq 0,\label{RT1}\\
			(\bm{P}^{RT}_j\bm{\omega},\bm{\sigma})_K&=(\bm{\omega},\bm{\sigma})_K,  \quad\forall \bm{\sigma} \in [P_{j-1}(K)]^n, \text{ for } j\geq 1.\label{RT2}
		\end{align}
	
The following lemmas show some   properties of  $\bm{P}^{RT}_j$.
 \begin{lemma} (cf. \cite{RT-1}) \label{Lemma 2.5}
 For any $\bm{\omega}_{hi}\in \bm{RT}_j(K),$ the relation $\nabla\cdot\bm{\omega}_{hi}|_K=0 $ implies $\bm{\omega}_{hi}\in [P_j(K)]^{n} $.
 \end{lemma}

\begin{lemma} (cf. \cite{RT-1}) \label{Lemma 2.6}
 For any $K\in \mathcal{T}_h$,  the following properties   hold:
 \begin{align}
    (\nabla\cdot\bm{P}^{RT}_j\bm{\omega}, q_h)_{K}&=(\nabla\cdot\bm{\omega}, q_h)_{K},& \forall \bm{\omega}\in[H^1(K)]^n,  q_h\in P_{j}(K),\\
		\|\bm{\omega}-\bm{P}^{RT}_j\bm{\omega}\|_{0,K}&\lesssim h_K^{l}|\bm{\omega}|_{l,K},&\forall\bm{\omega}\in [H^l(K)]^n, \forall 1\le l\le j+1.\label{RT3}
	\end{align}
 \end{lemma}

\begin{lemma}(cf. \cite{HX2019})  \label{Lemma 2.7}
For any $K\in \mathcal{T}_h$, $\bm{\omega}\in[H^l(K)]^n$ and $1\leq l\leq j+1$,
the following estimates hold:
\begin{align*}
|\bm{\omega}-P_j^{RT}\bm{\omega}|_{1,K}
&\lesssim
h_K^{l-1}|\bm{\omega}|_{l,K}, \\
|\bm{\omega}-P_j^{RT}\bm{\omega}|_{0,\partial K}
&\lesssim
h_K^{l-\frac{1}{2}}|\bm{\omega}|_{l,K},\\
|\bm{\omega}-P_j^{RT}\bm{\omega}|_{0,3,K}&\lesssim h_K^{l-\frac{n}{6}}|\bm{\omega}|_{l,K},\\
|\bm{\omega}-P_j^{RT}\bm{\omega}|_{0,3,\partial K}&\lesssim h_K^{l-\frac{1}{3}-\frac{n}{6}}|\bm{\omega}|_{l,K}.
\end{align*}
\end{lemma}
We have the following commutativity properties for the   RT projection, the $L^2$ projections and the discrete weak operators:
\begin{lemma}(cf. \cite{CFX2016})\label{Lemma 2.8}
 For $m\geq 1$, there hold
\begin{align*}
\nabla _{w,m-1}\{ \bm{P}_m^{RT}\bm{\omega}, \bm{\Pi}_{k}^{B}\bm{\omega} \}&=\bm{\Pi}_{m-1}^{\ast}(\nabla\bm{\omega} ),\quad \forall \bm{\omega}\in [H^{1}(\Omega)]^{n}, k=m, m-1,\\
\nabla _{w,m}\{\Pi_{m-1}^{\ast}q,\Pi _m ^{B}q\}&=\bm{\Pi}_{m}^{\ast}(\nabla q ),\quad\forall q\in H^{1}(\Omega).
\end{align*}
 \end{lemma}
 Finally, we give several inequalities to be used later (cf. \cite{BorggaardJeff2012Tdot,CiarletP.G1978TFEM,LiMinghao2019Tmfe}).

\begin{lemma} \label{Lemma 2.9}
For any $  \lambda,\mu\in \mathbb{R}^{n}$ and $r\geq2$, there hold 
\begin{align*}
|| \lambda|^{r-2}-|\mu |^{r-2}|&\leq C_{r} (|\lambda|^{r-3}+|\mu|^{r-3})|\lambda-\mu|,\\
|| \lambda|^{r-2}\lambda-|\mu |^{r-2}\mu|&\leq C_{r} (|\lambda|+|\mu|)^{r-2}|\lambda-\mu|,\\
|| \lambda|^{r-2}-|\mu |^{r-2}-(r-2)|\mu |^{r-4}\mu\cdot(\lambda-\mu)|&\leq C_{r} (|\lambda|^{r-4}+|\mu|^{r-4})|\lambda-\mu|^{2},\\
(| \lambda|^{r-2}\lambda-|\mu |^{r-2}\mu, \lambda-\mu)&\gtrsim |\lambda-\mu|^{r},
\end{align*}
where  $|\cdot|$ denotes the Euclid norm and $C_{r}$ is a positive constant only depending on $r$.
\end{lemma}

\section{ Well-posedness of discrete scheme }

Lemmas \ref{Lemma 3.1} and \ref{theoremLBB} give some stability  conditions for the discrete scheme \eqref{WG}. 
\begin{lemma}\label{Lemma 3.1}
For any $ \bm{\kappa}_h=\{\bm{\kappa}_{hi},\bm{\kappa}_{hb}\},
\bm{u}_h=\{\bm{u}_{hi},\bm{u}_{hb}\},
\bm{v}_h=\{\bm{v}_{hi},\bm{v}_{hb}\}
\in\bm{V}_h^0$, there hold
\begin{align}
a_h(\bm{u}_h,\bm{v}_h)&\lesssim \nu|||\bm{u}_{h}|||_{V}\cdot|||\bm{v}_{h}|||_{V},\label{a1}\\
a_h(\bm{v}_h,\bm{v}_h)&= \nu|||\bm{v}_{h}|||^{2}_{V},\label{a2}\\
c_{h}(\bm{v}_h;\bm{v}_h,\bm{v}_h)&= \alpha\|\bm{v}_{hi}\|_{0,r}^{r},\label{c2}\\
c_{h}(\bm{\kappa}_h;\bm{u}_h,\bm{v}_h)&\leq \alpha C_{\widetilde{r}}^{r} |||\bm{\kappa}_{h}|||^{r-2}_{V}|||\bm{u}_{h}|||_{V}\cdot|||\bm{v}_{h}|||_{V},\label{c1}\\
d_h(\bm{\kappa}_h;\bm{v}_h,\bm{v}_h)&=0 , \label{d0}\\
d_h(\bm{\kappa}_h;\bm{u}_h,\bm{v}_h)
 &\lesssim |||\bm{\kappa}_h|||_{V} \cdot |||\bm{u}_h|||_{V}\cdot |||\bm{v}_h|||_{V},\label{dh1}
\end{align}
where $C_{\widetilde{r}}$ is the same as  in \eqref{Cr-est}.
\end{lemma}
\begin{proof}
According to the definition of $a_h(\cdot,\cdot)$, the  Cauchy-Schwarz inequality and  Lemma \ref{Lemma 2.4}, we easily get \eqref{a1} - \eqref{a2}.
The results \eqref{c2} and \eqref{d0} follow from the definitions of $c_h(\cdot;\cdot,\cdot)$ and $d_h(\cdot;\cdot,\cdot)$, respectively.
 From the definition of  $c_h(\cdot;\cdot,\cdot)$, the H\"{o}lder's inequality and Lemma \ref{Lemma 2.4}
 we obtain
\begin{align*}
c_{h}(\bm{\kappa}_h;\bm{u}_h,\bm{v}_h)
\leq \alpha ||\bm{\kappa}_{hi}||^{r-2}_{0,r}
 ||\bm{u}_{hi}||_{0,r}
 ||\bm{v}_{hi}||_{0,r} 
\leq C_{\widetilde{r}}^{r} \alpha|||\bm{\kappa}_{h}|||^{r-2}_{V}|||\bm{u}_{h}|||_{V}\cdot|||\bm{v}_{h}|||_{V},
\end{align*}
i.e. \eqref{c1} holds.  The inequality \eqref{dh1} has been proved in \cite[Lemma 3.10]{HX2019}.  %
 \end{proof}

We also have the following discrete inf-sup inequality.
 \begin{lemma}(\cite{CFX2016})\label{theoremLBB}
 There holds
\begin{eqnarray*}
\sup_{ \bm{v}_{h}\in \bm{V}_h^{0}}\frac{b_h(\bm{v}_{h},p_h)}{|||\bm{v}_{h}|||_{V} }\gtrsim |||p_h|||_{Q}, \quad \forall p_h\in Q_h^{0}.
\end{eqnarray*}
\end{lemma}

Denote
\begin{eqnarray*}
\bm{V}_{0h}:=\{\bm{\kappa}_h\in\bm{V}_h^{0}:b_h(\bm{\kappa}_{h},q_h)=0,\forall q_h\in Q_h^{0}\}.
\end{eqnarray*}
From the proof of Theorem \ref{TH2.2} we easily see that
\begin{eqnarray*}
\bm{V}_{0h}=\{\bm{\kappa}_h\in\bm{V}_h^{0}:\bm{\kappa}_{hi}\in \bm{H}({\rm div};\Omega), \nabla\cdot \bm{\kappa}_{hi}=0\}.
\end{eqnarray*}
To prove the existence of solutions to the scheme \eqref{WG}, we     introduce the following auxiliary system: find $\bm{u}_h\in \bm{V}_{0h}$ such that
\begin{align}\label{wg1}
		\mathcal{B}_h(\bm{u}_h;\bm{u}_h,\bm{v}_h )=&(\bm{f},\bm{v}_{hi}), \quad\forall \bm{v}_h\in \bm{V}_{0h},
\end{align}
where  the trilinear form $\mathcal{B}_h(\cdot;\cdot,\cdot ): \bm{V}_{0h}\times \bm{V}_{0h}\times \bm{V}_{0h}\rightarrow \mathbb R$ is defined by
 \begin{align*}
\mathcal{B}_h(\bm{\kappa}_h;\bm{u}_h,\bm{v}_h ):=&a_{h}(\bm{u}_h, \bm{v}_h )+c_h(\bm{\kappa}_h;\bm{u}_h,\bm{v}_h )+  d_h(\bm{\kappa}_h;\bm{u}_h,\bm{v}_h )
\end{align*}
for any $ \bm{\kappa}_h,
\bm{u}_h,
\bm{v}_h
\in\bm{V}_{0h}$.

We have the following equivalence result:
\begin{lemma} \label{equivalent2}
The discrete problems \eqref{WG} and \eqref{wg1} are equivalent in the sense that both $(i)$ and $(ii)$ hold:\begin{itemize}
\item[(i)] If $(\bm{u}_h,p_h)\in \bm{V}_h^{0}\times Q_h^{0}$ solves \eqref{WG}, then $\bm{u}_h\in \bm{V}_{0h}$ solves   \eqref{wg1};
\item[(ii)]  If $\bm{u}_h\in \bm{V}_{0h}$ solves   \eqref{wg1}, then  $(\bm{u}_h,p_h)$ solves \eqref{WG}, where $p_h\in Q_h^{0} $ is determined by
\begin{align}\label{ph-problem}
		b_h(\bm{v}_{h},p_h)=(\bm{f},\bm{v}_{hi})-a_{h}(\bm{u}_h, \bm{v}_h )
		-c_h(\bm{u}_h ;\bm{u}_h ,\bm{v}_h)-d_h(\bm{u}_h ;\bm{u}_h ,\bm{v}_h), \forall \bm{v}_{h}\in \bm{V}_h^{0}.
	\end{align}
\end{itemize}
\end{lemma}

Define
$$\mathcal{N}_h:=\sup_{\bm 0\neq\bm{\kappa}_h,\bm{u}_h,\bm{v}_h\in\bm{V}_{0 h}}
\frac{d_h(\bm{\kappa}_h;\bm{u}_h,\bm{v}_h)}{|||\bm{\kappa}_h|||_{V}\cdot |||\bm{u}_h|||_{V}\cdot|||\bm{v}_h|||_{V} }, \quad \|\bm{f}\|_{*,h}:=  \sup\limits_{\bm 0\neq \bm{v}_h\in\bm{V}_{0h}}\frac{(\bm f,\bm v_{hi})}{|||\bm{v}_h|||_{V}} .$$
 From \eqref{dh1}   and Lemma \ref{Lemma 2.4} we easily know that
 $$\mathcal{N}_h\lesssim 1, \quad \|\bm{f}\|_{*,h}\lesssim \|\bm{f}\|_0.$$

Based on Lemma \ref{equivalent2}, we can obtain the following existence and boundedness results for the WG method:
\begin{theorem} \label{TH42}
The WG scheme \eqref{WG} admits at least a solution pair
	$(\bm{u}_h,p_h)\in \bm{V}_h^{0}\times Q_h^{0}$
and  there hold
\begin{align}
	|||\bm{u}_h|||_{V} &\le\frac{1}{\nu}\|\bm{f}\|_{*,h},\label{sta}\\
|||p_h|||_{Q} &\lesssim  \|\bm{f}\|_{*,h}+\|\bm{f}\|_{*,h}^2+\|\bm{f}\|_{*,h}^{r-1}. \label{p-sta}
\end{align}
\end{theorem}
\begin{proof}
We first show the problem \eqref{wg1}  admits at least one solution $\bm{u}_{h}\in\bm{V}_{0h}$. According to
\cite[Theorem 1.2]{Girault.V;Raviart.P1986}, it suffices to show that
the following two results hold:
\begin{itemize}
\item[(I)] $\mathcal{B}_h(\bm{v}_h;\bm{v}_h,\bm{v}_h )\geq\nu|||\bm{v}_h|||^2, \quad \forall  \bm{v}_h\in \bm{V}_{0h}$;
	
\item[(II)] $\bm{V}_{0h}$ is separable,  and the relation $\lim\limits_{l\to\infty}\bm{u}_h^{(l)}=\bm{u}_h$ (weakly in $\bm{V}_{0h}$) implies
\begin{eqnarray*}
\lim_{l\to\infty}\mathcal{B}_h(\bm{u}_h^{(l)};\bm{u}_h^{(l)},\bm{v}_h )
	=\mathcal{B}_h(\bm{u}_h;\bm{u}_h,\bm{v}_h ),\quad\forall\bm{v}_h\in \bm{V}_{0h}.
\end{eqnarray*}
\end{itemize}

In fact,   (I) follows from Lemma \ref{Lemma 3.1}    directly. We only need to show (II).
	Since $\bm{V}_{0h}$ is a finite dimensional space, we know that $\bm{V}_{0h}$ is separable and that  the weak convergence $\lim\limits_{l\to\infty}\bm{u}_h^{(l)}=\bm{u}_h$ on $\bm{V}_{0h}$ is equivalent to the strong convergence
	\begin{eqnarray}
	\lim_{l\to\infty}||| \bm{u}_h^{(l)}-\bm{u}_{h}|||_{V} =0.\label{45}
	\end{eqnarray}
On the other hand, by Lemmas \ref{Lemma 2.4}, \ref{Lemma 2.9}, \ref{Lemma 3.1} and the definition of $\mathcal{N}_{h}$, we have
	\begin{align}
	&|\mathcal{B}_h(\bm{u}_{h}^{(l)};\bm{u}_{h}^{(l)},\bm{v}_{h} )
	-\mathcal{B}_h(\bm{u}_{h};\bm{u}_{h},\bm{v}_{h} )|\nonumber\\
	=&|a_{h}(\bm{u}_{h}^{(l)}-\bm{u}_h, \bm{v}_h )
+\big(c_h(\bm{u}_{h}^{(l)} ;\bm{u}_{h}^{(l)} ,\bm{v}_{h} )-c_h(\bm{u}_{h} ;\bm{u}_{h} ,\bm{v}_{h})\big)
   +d_h(\bm{u}_{h}^{(l)}-\bm{u}_{h};\bm{u}_{h}^{(l)}-\bm{u}_{h},\bm{v}_{h})\nonumber\\
	&+ d_h(\bm{u}_{h};\bm{u}_{h}^{(l)}-\bm{u}_{h},\bm{v}_{h})
    +d_h(\bm{u}_{h}^{(l)}-\bm{u}_{h};\bm{u}_{h},\bm{v}_{h})|\nonumber\\
	\le& \nu|||\bm{u}_{h}^{(l)}-\bm{u}_{h}|||_{V}\cdot|||\bm{v}_{h}|||_{V}
    +C_{\widetilde{r}}^{r}C_{r}\alpha|||\bm{u}^{(l)}_{h}-\bm{u}_{h}|||_{V} (|||\bm{u}_{h}^{(l)}|||_{V}+ |||\bm{u}_{h}|||_{V} )^{r-2}|||\bm{v}_{h} |||_{V}\nonumber\\
    &+\mathcal{N}_h|||\bm{u}_{h}^{(l)}-\bm{u}_{h}|||^2_{V}
	|||\bm{v}_{h}|||_{V}
	+2\mathcal{N}_h|||\bm{u}_{h}^{(l)}-\bm{u}_{h}|||_{V}\cdot|||\bm{u}_{h}|||_{V}\cdot |||\bm{v}_{h}|||_{V},\nonumber
	\end{align}
which, together with \eqref{45}, yields
 $$\lim_{l\to\infty}\mathcal{B}_h(\bm{u}_{h}^{(l)};\bm{u}_{h}^{(l)},\bm{v}_{h} )
	=\mathcal{B}_h(\bm{u}_{h};\bm{u}_{h},\bm{v}_{h} ),\
	\forall\bm{v}_{h}\in \bm{V}_{0h},$$
	i.e.  (II)   holds.	Hence, \eqref{wg1}  has at least one solution $\bm{u}_{h}\in\bm{V}_{0h}$.
	
	For a given $\bm{u}_{h}\in\bm{V}_{0h}\subset \bm{V}_h^{0}$, by Lemma \ref{theoremLBB} there is a unique $p_h\in Q_h^{0}$ satisfying \eqref{ph-problem}.
	Thus, in light of  Lemma \ref{equivalent2} we know that $(\bm{u}_h,p_h)$
	is a solution of \eqref{WG}.
	
	Taking $\bm{v}_h=\bm{u}_h,q_h=p_h$ in \eqref{WG} and using
 Lemma \ref{Lemma 3.1}, we immediately get
\begin{eqnarray*}
\nu|||\bm{u}_h||| ^2_{V}+\alpha\|\bm{u}_{hi}\|_{L^{r}}^{r}
=(\bm{f},\bm{u}_{hi})\le \|\bm{f}\|_{*,h}|||\bm{u}_h|||_{V},
\end{eqnarray*}
which implies \eqref{sta}.  Finally,  the estimate \eqref{p-sta} follows from Lemmas \ref{Lemma 2.4},  \ref{Lemma 3.1} and \ref{theoremLBB}, the equation \eqref{ph-problem} and the estimate \eqref{sta}. This finishes the proof.
\end{proof}

Furthermore, we have the following uniqueness result:

\begin{theorem} \label{Theorem 4.3}
Assume that the smallness condition
\begin{eqnarray}
	\frac{\mathcal{N}_h}{\nu^{2}}\|\bm{f}\|_{*,h}<1 \label{uni-condi}
	\end{eqnarray}
holds, then the scheme \eqref{WG} admits a unique solution $(\bm{u}_h,p_h)\in \bm{V}_h^{0}\times Q_h^{0}$. 
\end{theorem}
\begin{proof}
	  Let   $(\bm{u}_{h1},p_{h1})$ and $(\bm{u}_{h2},p_{h2})$ be two solutions of \eqref{WG}, i.e.   for $j=1,2$ and $(\bm{v}_h,q_h)\in \bm{V}_h^{0}\times Q_h^{0}$ there hold
	\begin{align*}
	a_{h}(\bm{u}_{hj},\bm{v}_h)+c_h(\bm{u}_{hj};\bm{u}_{hj},\bm{v}_{h})
+d_h(\bm{u}_{hj};\bm{u}_{hj},\bm{v}_{h})
+b_h(\bm{v}_{hj},p_{hj})
&=(\bm{f},\bm{v}_{hi}),\\
b_h(\bm{u}_{hj},q_h)&=0,
\end{align*}
which give
\begin{subequations}
\begin{align}
&a_{h}(\bm{u}_{h1}-\bm{u}_{h2},\bm{v}_h)+c_h(\bm{u}_{h1};\bm{u}_{h1},\bm{v}_{h})-c_h(\bm{u}_{h2};\bm{u}_{h2},\bm{v}_{h})-b_h(\bm{v}_{h},p_{h1}-p_{h2})\nonumber\\
&\qquad\qquad\qquad \quad \ \ = d_h(\bm{u}_{h2};\bm{u}_{h2},\bm{v}_h)- d_h(\bm{u}_{h1};\bm{u}_{h1},\bm{v}_h),\label{Ua}\\
&b_h(\bm{u}_{h1}-\bm{u}_{h2},q_h)=0.\label{Ub}
\end{align}
\end{subequations}
	Taking $\bm{v}_h=\bm{u}_{h1}-\bm{u}_{h2}$ and $q_h=p_{h1}-p_{h2} $ in the above two equations and using the relation
	$$ d_h(\bm{u}_{h2};\bm{u}_{h1}-\bm{u}_{h2},\bm{u}_{h1}-\bm{u}_{h2})=0 $$
	due to \eqref{d0},  we  obtain
\begin{subequations}
\begin{align*}
&a_{h}(\bm{u}_{h1}-\bm{u}_{h2},\bm{u}_{h1}-\bm{u}_{h2})+c_h(\bm{u}_{h1};\bm{u}_{h1},\bm{u}_{h1}-\bm{u}_{h2})-c_h(\bm{u}_{h2};\bm{u}_{h2},\bm{u}_{h1}-\bm{u}_{h2}) \nonumber\\
 =& d_h(\bm{u}_{h2};\bm{u}_{h2},\bm{u}_{h1}-\bm{u}_{h2})- d_h(\bm{u}_{h1};\bm{u}_{h1},\bm{u}_{h1}-\bm{u}_{h2})\\
=& d_h(\bm{u}_{h2}-\bm{u}_{h1};\bm{u}_{h1},\bm{u}_{h1}-\bm{u}_{h2}).
\end{align*}
\end{subequations}	
This relation, together with \eqref{a2}, \eqref{dh1}, \eqref{sta},  and 
 the inequality
$$c_h(\bm{u}_{h1};\bm{u}_{h1},\bm{u}_{h1}-\bm{u}_{h2})-c_h(\bm{u}_{h2};\bm{u}_{h2},\bm{u}_{h1}-\bm{u}_{h2})\gtrsim \|\bm{u}_{hi1}-\bm{u}_{hi2}\|^{r}_{0,r}\geq 0$$
due to Lemma \ref{Lemma 2.9}, implies  
\begin{align*}
\nu|||\bm{u}_{h1}-\bm{u}_{h2}|||^2_{V}
\leq 
&d_h(\bm{u}_{h2}-\bm{u}_{h1};\bm{u}_{h1},\bm{u}_{h1}-\bm{u}_{h2})\nonumber\\
\le&\mathcal{N}_h||| \bm{u}_{h1}|||_{V} ||| \bm{u}_{h1}-\bm{u}_{h2}|||^2_{V} \nonumber\\
\le&\nu^{-1}\mathcal{N}_h\|\bm{f}\|_{*,h}|||\bm{u}_{h1}-\bm{u}_{h2}|||^2_{V},
\end{align*}
i.e.
\begin{eqnarray*}
(1-\frac{\mathcal{N}_h}{\nu^{2}}\|\bm{f}\|_{*,h})|||\bm{u}_{h1}-\bm{u}_{h2}|||^2_{V} \le 0.
\end{eqnarray*}
This  inequality plus the assumption \eqref{uni-condi} yields $\bm{u}_{h1}=\bm{u}_{h2}$. Then by
 \eqref{Ua}  we have $b_h(\bm{v}_{h},p_{h1}-p_{h2})=0$ which, together with  Lemma \ref{theoremLBB}, leads to   $p_{h1}=p_{h2}$. This finishes the proof.
\end{proof}

\section{A priori error estimates}
\label{section6}
This section is devoted to the error analysis   for the WG method \eqref{WG}.
To   this end, we first assume that the weak solution   $(\bm{u},p) $ of  $(\ref{BF0})$ satisfies the following regularity conditions:
\begin{equation}\label{regularity}
 \bm{u} \in [H^{m+1}(\Omega)]^n\cap \bm{V}, \quad p\in  H^{m}(\Omega)\cap L_0^2(\Omega).
  \end{equation}
Define
\begin{align}
\bm{\mathcal{I}}_{h}\bm{u}|_K:=\{\bm{P}_{m}^{RT}(\bm{u}|_K), \bm{\Pi} _{k}^{B}(\bm{u}|_K)\},\quad
\mathcal{P}_{h}p|_K:=\{ \Pi_{m-1}^{\ast}(p|_K),\Pi_{m}^{B}(p|_K) \}.
\end{align}
Here we recall that $m\geq 1$ and $k=m-1, m.$

\begin{lemma}\label{Lemma 4.1}
There hold
	\begin{eqnarray}
	\bm{P}^{RT}_m\bm{u}|_K\in [P_{m}(K)]^n, \quad \forall K\in\mathcal{T}_h \label{4.5}
	\end{eqnarray}
and, for  any $ (\bm{v}_h,q_h)\in \bm{V}_h^{0}\times Q_h^{0}$, 
\begin{subequations}
\begin{align}
&a_h(\bm{\mathcal{I}}_{h}\bm{u},\bm{v}_h)+b_h(\bm{v}_h,\mathcal{P}_{h}p)
+c_h(\bm{\mathcal{I}}_{h}\bm{u};\bm{\mathcal{I}}_{h}\bm{u},\bm{v}_h)
+ d_h(\bm{\mathcal{I}}_{h}\bm{u};\bm{\mathcal{I}}_{h}\bm{u},\bm{v}_h)\nonumber\\
&\qquad \qquad \quad =(\bm{f},\bm{v}_{hi}) +\xi_{I}(\bm{u};\bm{u},\bm{v}_h)+\xi_{II}(\bm{u},\bm{v}_h)+\xi_{III}(\bm{u};\bm{u},\bm{v}_h),\label{e2}\\
	&b_h(\bm{\mathcal{I}}_{h}\bm{u},q_h)=0,\label{e3}
\end{align}
\end{subequations}
	where 
	\begin{align*}
&\xi_{I}(\bm{u};\bm{u},\bm{v}_{h}):= -\frac{1}{2}(\bm{P}_m^{RT}\bm{u}\otimes \bm{P}_m^{RT}\bm{u}-\bm{u}\otimes\bm{u},\nabla_h\bm{v}_{hi}) 
+\frac{1}{2}\langle(\bm{\Pi}_{k}^{B}\bm{u}\otimes \bm{\Pi}_{k}^{B}\bm{u}
-\bm{u}\otimes\bm{u})\bm{n},\bm{v}_{hi}\rangle_{\partial\mathcal{T}_h}\nonumber\\
&\qquad\hskip2cm -\frac{1}{2}((\bm{u}\cdot\nabla)\bm{u}-(\bm{P}_m^{RT}\bm{u}\cdot\nabla_h)\bm{P}_m^{RT}\bm{u},\bm{v}_{hi})
 -\frac{1}{2}\langle(\bm{v}_{hb}\otimes \bm{\Pi}_{k}^{B}\bm{u})\bm{n},\bm{P}_m^{RT}\bm{u}\rangle_{\partial\mathcal{T}_h},\\
	&\xi_{II}(\bm{u},\bm{v}_h):= \nu\langle \big(\nabla\bm u  -\bm{\Pi}_{m-1}^\ast\nabla\bm u\big)\bm{n},\bm{v}_{hi}-\bm{v}_{hb} \rangle_{\partial\mathcal{T}_h} +\nu\langle \eta(\bm{P}^{RT}_m\bm{u}-{\bm{u}}),\bm{\Pi}_{k}^{B}\bm{v}_{hi}-\bm{v}_{hb} \rangle_{\partial\mathcal{T}_h},\nonumber\\
&\xi_{III}(\bm{u};\bm{u},\bm{v}_h)=
\alpha (|\bm{P}^{RT}_m\bm{u}|^{r-2}\bm{P}^{RT}_m\bm{u}-|\bm{u}|^{r-2}\bm{u},\bm{v}_{hi}).
	\end{align*}
	\end{lemma}
\begin{proof}
For $\forall K\in\mathcal{T}_h,\varrho_m\in P_{m}(K)$, by Lemma \ref{Lemma 2.6} we have
	\begin{eqnarray*}
		(\nabla\cdot\bm{P}^{RT}_{m}\bm{u},\varrho_m)_K=(\nabla\cdot\bm{u},\varrho_m)_K=0,
	\end{eqnarray*}
	which means that $\nabla\cdot\bm{P}^{RT}_m\bm{u}=0$, i.e. \eqref{4.5} holds.
	
	 From the definition  of discrete weak divergence, the  Green's formula and the definition of   the trilinear form $d_h(\cdot;\cdot,\cdot)$ we easily have
	 \begin{eqnarray*}
	d_h(\bm{\mathcal{I}}_{h}\bm{u} ;\bm{\mathcal{I}}_{h}\bm{u},\bm{v}_{h})
	=(\nabla\cdot(\bm{u}\otimes\bm{u}),\bm{v}_{hi})+\xi_{I}(\bm{u};\bm{u},\bm{v}_{h}). 
	\end{eqnarray*}
Thus, according to  the definition of discrete weak gradient, the  Green's formula, the   projection properties,    Lemma    \ref{Lemma 2.8}  and  the first equation of \eqref{BF0}, we get
\begin{align*}
&a_h(\bm{\mathcal{I}}_{h}\bm{u},\bm{v}_h)+b_h(\bm{v}_h,\mathcal{P}_{h}p)
+c_h(\bm{\mathcal{I}}_{h}\bm{u};\bm{\mathcal{I}}_{h}\bm{u},\bm{v}_h)
+ d_h(\bm{\mathcal{I}}_{h}\bm{u};\bm{\mathcal{I}}_{h}\bm{u},\bm{v}_h)\nonumber\\
=&\nu (\nabla_{w,m-1}\{ \bm {P}_{m}^{RT}\bm{u},\bm {\Pi}_{k}^{B}\bm{u} \},\nabla_{w,m-1}\bm{v}_{h}  )
+\nu\langle \eta \bm {\Pi}_{k}^{B}(\bm {P}_{m}^{RT}\bm{u}-\bm{u}), \bm{\Pi}_{k}^{B}\bm{v}_{hi}-\bm{v}_{hb}\rangle_{\partial\mathcal{T}_h}\nonumber\\
&+(\bm{v}_{hi}, \nabla_{w,m}\{  \Pi_{m-1}^{RT}p, \Pi_{m}^{B}p \} )
+\alpha(| \bm {P}_{m}^{RT}\bm{u}|^{r-2}\bm {P}_{m}^{RT}\bm{u},\bm{v}_{hi} )\nonumber\\
&+(\nabla\cdot(\bm{u}\otimes\bm{u}),\bm{v}_{hi})+\xi_{I}(\bm{u};\bm{u},\bm{v}_{h})\nonumber\\
=&\nu ( \bm {\Pi}_{m-1}^{\ast}(\nabla\bm{u}), \nabla_{w,m}\bm{v}_{h} )
+\nu\langle \eta (\bm {P}_{m}^{RT}\bm{u}- \bm{u}), \bm{\Pi}_{k}^{B}\bm{v}_{hi}-\bm{v}_{hb}\rangle_{\partial\mathcal{T}_h}\nonumber\\
&+(\bm{v}_{hi}, \bm {\Pi}_{m}^{\ast}(\nabla p))
+\alpha(| \bm {P}_{m}^{RT}\bm{u}|^{r-2}\bm {P}_{m}^{RT}\bm{u},\bm{v}_{hi} )\nonumber\\
&+(\nabla\cdot(\bm{u}\otimes\bm{u}),\bm{v}_{hi})+\xi_{I}(\bm{u};\bm{u},\bm{v}_{h})\nonumber\\
=&-\nu (\nabla_{h}\cdot \bm {\Pi}_{m-1}^{\ast}(\nabla\bm{u}), \bm{v}_{hi} )
+ \langle\bm {\Pi}_{m-1}^{\ast}(\nabla\bm{u})\bm{n}, \bm{v}_{hb} \rangle_{\partial\mathcal{T}_h}\nonumber\\
&+\nu\langle \eta (\bm {P}_{m}^{RT}\bm{u}- \bm{u}), \bm{\Pi}_{k}^{B}\bm{v}_{hi}-\bm{v}_{hb}\rangle_{\partial\mathcal{T}_h}\nonumber\\
&+(\bm{v}_{hi}, \bm {\Pi}_{m}^{\ast}(\nabla p))
+\alpha(| \bm {P}_{m}^{RT}\bm{u}|^{r-2}\bm {P}_{m}^{RT}\bm{u},\bm{v}_{hi} )\nonumber\\
&+(\nabla\cdot(\bm{u}\otimes\bm{u}),\bm{v}_{hi})+\xi_{I}(\bm{u};\bm{u},\bm{v}_{h})\nonumber\\
=&-\nu( \triangle\bm{u}, \bm{v}_{hi} )
+\nu\langle(\nabla \bm{u}-\bm {\Pi}_{m-1}^{\ast}\nabla\bm{u})\bm{n}, \bm{v}_{hi} -\bm{v}_{hb} \rangle_{\partial\mathcal{T}_h}\nonumber\\
&+\nu\langle \eta (\bm {P}_{m}^{RT}\bm{u}- \bm{u}), \bm{\Pi}_{k}^{B}\bm{v}_{hi}-\bm{v}_{hb}\rangle_{\partial\mathcal{T}_h}\nonumber\\
&+(\bm{v}_{hi}, \nabla p)
+\alpha(| \bm{u}|^{r-2}\bm{u},\bm{v}_{hi} )
+\alpha(| \bm {P}_{m}^{RT}\bm{u}|^{r-2}\bm {P}_{m}^{RT}\bm{u}-| \bm{u}|^{r-2}\bm{u},\bm{v}_{hi} )\nonumber\\
&+(\nabla\cdot(\bm{u}\otimes\bm{u}),\bm{v}_{hi})+\xi_{I}(\bm{u};\bm{u},\bm{v}_{h})\nonumber\\
=&(\bm{f}, \bm{v}_{hi})+\nu\langle (\nabla \bm{u}-\bm {\Pi}_{m-1}^{\ast}\nabla\bm{u})\bm{n}, \bm{v}_{hi} -\bm{v}_{hb} \rangle_{\partial\mathcal{T}_h}\nonumber\\
&+\nu\langle \eta (\bm {P}_{m}^{RT}\bm{u}- \bm{u}), \bm{\Pi}_{k}^{B}\bm{v}_{hi}-\bm{v}_{hb}\rangle_{\partial\mathcal{T}_h}\nonumber\\
&+\alpha(| \bm {P}_{m}^{RT}\bm{u}|^{r-2}\bm {P}_{m}^{RT}\bm{u}-| \bm{u}|^{r-2}\bm{u},\bm{v}_{hi} )
+\xi_{I}(\bm{u};\bm{u},\bm{v}_{h})\nonumber\\
=&(\bm{f}, \bm{v}_{hi})+\xi_{I}(\bm{u};\bm{u},\bm{v}_{h})
+\xi_{II}(\bm{u},\bm{v}_{h})+\xi_{III}(\bm{u};\bm{u},\bm{v}_{h}),
\end{align*}
which proves \eqref{e2}.

From  the definition of $\nabla_{w,m}$, the fact $\nabla\cdot\bm{P}^{RT}_{m}\bm{u}=0$ and  \eqref{RT1} it follows
\begin{align*}
b_h(\bm{\mathcal{I}}_{h}\bm{u},q_h  ) &=
-(\nabla\cdot\bm{P}^{RT}_{m}\bm{u},q_{hi})+\langle\bm{P}^{RT}_{m}\bm{u}\cdot\bm{n},q_{hb} \rangle_{\partial\mathcal{T}_h}
=\langle\bm{u}\cdot\bm{n},q_{hb} \rangle_{\partial\mathcal{T}_h}=0.
\end{align*}
i.e. the relation \eqref{e3} holds.
This completes the proof.
\end{proof}

By following a similar line as in the proofs of \cite[Lemma 4.3]{HX2019} and \cite[Lemma 5.2]{ZZX2023}, we can obtain the  estimates of  $\xi_{I}$, $\xi_{II}$, and $\xi_{III}$.

\begin{lemma} \label{Lemma 4.2}
For any $\bm{v}_h\in\bm{V}_h^{0}$,  there hold
\begin{subequations}
\begin{align}
|\xi_{I}(\bm{u},\bm{u};\bm{v}_h)|&\lesssim  h^{m}
\|\bm{u}\|_{2}\|\bm{u}\|_{m+1}|||\bm{v}_h|||_{V},\label{X1}\\
|\xi_{II}(\bm{u},\bm{v}_h)|&\lesssim   h^{m}\|\bm{u}\|_{ m+1}|||\bm{v}_h|||_{V},\label{X2}\\
|\xi_{III}(\bm{u},\bm{u}; \bm{v}_h)|&\lesssim h^{m}\|\bm{u}\|_{2}^{r-2}\|\bm{u}\|_{m+1}|||\bm{v}_h|||_{V},\label{X3}
\end{align}
\end{subequations}
for $k=m, m-1$ when $n=2$ and $k=m$ when $n=3$.
\end{lemma}
\begin{proof}
 We first estimate the four terms of $\xi_{I}(\bm{u},\bm{u};\bm{v}_h)$   one by one.
 By the the   H\"{o}lder's inequality, the   Sobolev inequality and   Lemmas  \ref{Lemma 2.4}, \ref{Lemma 2.6} and \ref{Lemma 2.7},  we   have 
\begin{align*}
&|-\frac{1}{2}(\bm{P}_m^{RT}\bm{u}\otimes \bm{P}_m^{RT}\bm{u}-\bm{u}\otimes\bm{u},\nabla_h\bm{v}_{hi})|\nonumber\\
\lesssim &|(\bm{P}_m^{RT}\bm{u}\otimes (\bm{P}_m^{RT}\bm{u}-\bm{u}),\nabla_h\bm{v}_{hi})|
+|((\bm{P}_m^{RT}\bm{u}-\bm{u}) \otimes \bm{u} , \nabla_h\bm{v}_{hi})|\nonumber\\
\lesssim &\sum_{K\in \mathcal{T}_{h}}|\bm{P}_m^{RT}\bm{u}-\bm{u}|_{0,3,K}|\bm{P}_m^{RT}\bm{u} |_{0,6,K}\|\nabla_h\bm{v}_{hi}\|_{0,K}\nonumber\\
&+|\bm{u} |_{0,\infty,\Omega}\sum_{K\in \mathcal{T}_{h}}|\bm{P}_m^{RT}\bm{u}-\bm{u}|_{0,K}\|\nabla_h\bm{v}_{hi}\|_{0,K}\nonumber\\
\lesssim &\sum_{K\in \mathcal{T}_{h}}|\bm{P}_m^{RT}\bm{u}-\bm{u}|_{0,3,K}(|\bm{P}_m^{RT}\bm{u} -\bm{u}|_{0,6,K}+|\bm{u} |_{0,6,K})\|\nabla_h\bm{v}_{hi}\|_{0,K}\nonumber\\
&+|\bm{u} |_{0,\infty,\Omega}\sum_{K\in \mathcal{T}_{h}}|\bm{P}_m^{RT}\bm{u}-\bm{u}|_{0,K}\|\nabla_h\bm{v}_{hi}\|_{0,K}\nonumber\\
\lesssim &(|\bm{P}_m^{RT}\bm{u} -\bm{u}|_{0,6,\Omega}+|\bm{u} |_{0,6,\Omega})\sum_{K\in \mathcal{T}_{h}}|\bm{P}_m^{RT}\bm{u}-\bm{u}|_{0,3,K}\|\nabla_h\bm{v}_{hi}\|_{0,K}\nonumber\\
&+|\bm{u} |_{0,\infty,\Omega}\sum_{K\in \mathcal{T}_{h}}|\bm{P}_m^{RT}\bm{u}-\bm{u}|_{0,K}\|\nabla_h\bm{v}_{hi}\|_{0,K}\nonumber\\
\lesssim &\|\bm{u}\|_{1}\sum_{K\in \mathcal{T}_{h}}|\bm{P}_m^{RT}\bm{u}-\bm{u}|_{0,3,K}\|\nabla_h\bm{v}_{hi}\|_{0,K}
+h^{m+1}|\bm{u} |_{0,\infty,\Omega}| \bm{u}|_{m+1}||| \bm{v}_{h} |||_{V}  \nonumber\\
\lesssim& h^{m+1-\frac n 6} \|\bm{u}\|_{1}	| \bm{u}|_{m+1}\|\nabla_h\bm{v}_{hi}\|_{0,K}
+h^{m+1}|\bm{u} |_{0,\infty,\Omega}| \bm{u}|_{m+1}||| \bm{v}_{h} |||_{V}  \nonumber\\
\lesssim& h^{m} \|\bm{u} \| _{2}  \| \bm{u}\|_{m+1}  |||\bm{v}_h|||_{V}.
\end{align*}
Similarly, there hold
\begin{align*}
  &|\frac{1}{2}\langle(\bm{\Pi}_{}^{B}\bm{u}\otimes \bm{\Pi}_{k}^{B}\bm{u}
-\bm{u}\otimes\bm{u})\ \bm{n},\bm{v}_{hi}\rangle_{\partial\mathcal{T}_h}|\nonumber\\
=&|\frac{1}{2}\langle(\bm{\Pi}_{k}^{B}\bm{u}\otimes \bm{\Pi}_{k}^{B}\bm{u}
-\bm{u}\otimes\bm{u}) \bm{n},\bm{v}_{hi}-\bm{v}_{hb}\rangle_{\partial\mathcal{T}_h}|\nonumber\\
\lesssim &|\langle(\bm{\Pi}_{k}^{B}\bm{u}-\bm{u})\otimes( \bm{\Pi}_{m}^{\ast}\bm{u}
-\bm{u}) \bm{n},\bm{v}_{hi}-\bm{v}_{hb}\rangle_{\partial\mathcal{T}_h}|\nonumber\\
&+|\langle(\bm{\Pi}_{k}^{B}\bm{u}-\bm{u})\otimes( \bm{\Pi}_{m}^{\ast}\bm{u}
 \bm{n},\bm{v}_{hi}-\bm{v}_{hb}\rangle_{\partial\mathcal{T}_h}|\nonumber\\
&+|\langle (\bm{\Pi}_{k}^{\ast}\bm{u}-\bm {\Pi}_{k}^{B}\bm{u})\otimes (\bm{\Pi}_{k}^{B}\bm{u}-\bm{u})  \bm{n},
 \bm{v}_{hi}-\bm{v}_{hb} \rangle_{\partial\mathcal{T}_h}|\nonumber\\
&+|\langle (( \bm{\Pi}_{k}^{B}\bm{u} -\bm{u})\otimes\bm{\Pi}_{k}^{\ast}\bm{u}) \bm{n} , \bm{v}_{hi}-\bm{v}_{hb}\rangle_{\partial\mathcal{T}_h}|\nonumber\\
\lesssim&\sum _{K\in \mathcal{T}_{h}}(| \bm{\Pi}_{k}^{B}\bm{u}-\bm{u}|_{0,\partial K}| \bm{\Pi}_{m}^{\ast}\bm{u}
-\bm{u} |_{0,\partial K}\left(| \bm{v}_{hi}-\bm{\Pi}_{k}^{B}\bm{v}_{hi}|_{0,\infty,\partial K}+| \bm{\Pi}_{k}^{B}\bm{v}_{hi}-\bm{v}_{hb}|_{0,\infty,\partial K}\right) \nonumber\\
&+ | \bm{\Pi}_{k}^{B}\bm{u}-\bm{u} |_{0,\partial K}| \bm{\Pi}_{m}^{\ast}\bm{u}|_{0,\infty,\partial K}\left(| \bm{v}_{hi}-\bm{\Pi}_{k}^{B}\bm{v}_{hi}|_{0,\partial K}+| \bm{\Pi}_{k}^{B}\bm{v}_{hi}-\bm{v}_{hb}|_{0,\partial K}\right)\nonumber\\
&+ |\bm{\Pi}_{k}^{\ast}\bm{u}- {\Pi}_{k}^{B}\bm{u}|_{0,\partial K}|\bm{\Pi}_{k}^{B}\bm{u}-\bm{u}|_{0,\partial K}\left(| \bm{v}_{hi}-\bm{\Pi}_{k}^{B}\bm{v}_{hi}|_{0,\infty,\partial K}+| \bm{\Pi}_{k}^{B}\bm{v}_{hi}-\bm{v}_{hb}|_{0,\infty,\partial K}\right) \nonumber\\
&+|\bm{\Pi}_{k}^{B}\bm{u}-\bm{u}|_{0,\partial K}| \bm{\Pi}_{m}^{\ast}\bm{u}|_{0,\infty,\partial K} \left(| \bm{v}_{hi}-\bm{\Pi}_{k}^{B}\bm{v}_{hi}|_{0,\partial K}+| \bm{\Pi}_{k}^{B}\bm{v}_{hi}-\bm{v}_{hb}|_{0,\partial K}\right)\nonumber\\
\lesssim&h^{m} \|\bm{u} \| _{2}  \| \bm{u}\|_{m+1}  |||\bm{v}_h|||_{V},
\end{align*}
\begin{align*}
 &|-\frac{1}{2}((\bm{u}\cdot\nabla)\bm{u}-(\bm{P}_m^{RT}\bm{u}\cdot\nabla_h)\bm{P}_m^{RT}\bm{u},\bm{v}_{hi})|\nonumber\\
\lesssim &|((\bm{u}- \bm{P}_m^{RT}\bm{u})\cdot\nabla\bm{u}, \bm{v}_{hi})|
+|\bm{P}_m^{RT}\bm{u}\cdot(\nabla\bm{u}- \nabla_h\bm{P}_m^{RT}\bm{u}),\bm{v}_{hi}) |\nonumber\\
\lesssim &\sum_{K\in \mathcal{T}_{h}}|\bm{u}- \bm{P}_m^{RT}\bm{u} |_{0,3,K}|\nabla\bm{u} |_{0,K}\|\bm{v}_{hi} \|_{0,6,K}\nonumber\\
&+\sum_{K\in \mathcal{T}_{h}}|\bm{P}_m^{RT}\bm{u}|_{0,6,K} |\nabla\bm{u}- \nabla_h\bm{P}_m^{RT}\bm{u} |_{0,K}\|\bm{v}_{hi} \|_{0,3,K}\nonumber\\
\lesssim &h^{m} \|\bm{u} \| _{2}  \| \bm{u}\|_{m+1}  |||\bm{v}_h|||_{V},
\end{align*}
and
 \begin{align}
  &|-\frac{1}{2}\langle(\bm{v}_{hb}\otimes \bm{\Pi}_{k}^{B}\bm{u})\bm{n},\bm{P}_m^{RT}\bm{u}\rangle_{\partial\mathcal{T}_h}|\nonumber\\
 =&|\frac{1}{2}\langle(\bm{v}_{hb}\otimes \bm{\Pi}_{k}^{B}\bm{u})\bm{n},\bm{P}_m^{RT}\bm{u}-\bm{\Pi}_m^{B}\bm{u} \rangle_{\partial\mathcal{T}_h}|\nonumber\\
 \lesssim & | \langle (\bm{v}_{hi}-\bm{v}_{hb} )\otimes (\bm{\Pi}_{k}^{B}\bm{u} -\bm{\Pi}_{k}^{\ast}\bm{u})\bm{n}, \bm{P}_m^{RT}\bm{u}-\bm{\Pi}_m^{B}\bm{u} \rangle_{\partial\mathcal{T}_h} |\nonumber\\
 &+| \langle \bm{v}_{hi}\otimes (\bm{\Pi}_{k}^{B}\bm{u} -\bm{\Pi}_{k}^{\ast}\bm{u})\bm{n}, \bm{P}_m^{RT}\bm{u}-\bm{\Pi}_m^{B}\bm{u} \rangle_{\partial\mathcal{T}_h} |\nonumber\\
 &+|\langle (\bm{v}_{hi}-\bm{v}_{hb} )\otimes \bm{\Pi}_{k}^{\ast}\bm{u}\bm{n},\bm{P}_m^{RT}\bm{u}-\bm{\Pi}_m^{B}\bm{u}\rangle_{\partial\mathcal{T}_h} |\nonumber\\
 &+|\langle \bm{v}_{hi}\otimes \bm{\Pi}_{k}^{\ast}\bm{u}\bm{n},\bm{P}_m^{RT}\bm{u} -\bm{\Pi}_{m}^{B}\bm{u} \rangle_{\partial\mathcal{T}_h} |\nonumber\\
 \lesssim & \sum_{K\in \mathcal{T}_{h}}( (|\bm{v}_{hi}-\bm{\Pi}_{k}^{B}\bm{v}_{hi}  |_{0,\infty, \partial K}
  +|\bm{\Pi}_{k}^{B}\bm{v}_{hi}-\bm{v}_{hb}  |_{0,\infty, \partial K})|\bm{\Pi}_{k}^{B}\bm{u} -\bm{\Pi}_{k}^{\ast}\bm{u}|_{0, \partial K}|\bm{P}_m^{RT}\bm{u}-\bm{\Pi}_m^{B}\bm{u}  |_{0, \partial K}\nonumber\\
 &+|\bm{v}_{hi} |_{0,\infty, \partial K}|\bm{\Pi}_{k}^{B}\bm{u} -\bm{\Pi}_{k}^{\ast}\bm{u}|_{0, \partial K}|\bm{P}_m^{RT}\bm{u}-\bm{\Pi}_m^{B}\bm{u}|_{0, \partial K}\nonumber\\
 &+(|\bm{v}_{hi}-\bm{\Pi}_{k}^{B}\bm{v}_{hi}|_{0, \partial K}+|\bm{\Pi}_{k}^{B}\bm{v}_{hi}-\bm{v}_{hb}|_{0, \partial K})|\bm{\Pi}_{k}^{\ast}\bm{u} |_{0,6, \partial K}|\bm{P}_m^{B}\bm{u} -\bm{\Pi}_{k}^{B}\bm{u} |_{0,3, \partial K}\nonumber\\
 &+|  \bm{v}_{hi}|_{0,3, \partial K}|\bm{\Pi}_{k}^{\ast}\bm{u} |_{0,6, \partial K}|\bm{P}_m^{RT}\bm{u} -\bm{\Pi}_{m}^{B}\bm{u}  |_{0, \partial K}\nonumber\\
 \lesssim &h^{m} \|\bm{u} \| _{2}  \| \bm{u}\|_{m+1}  |||\bm{v}_h|||_{V}.
\end{align}
Combining the above four estimates leads to the  desired result \eqref{X1}.

Similarly, we have
\begin{align*}
&|\xi_{II}(\bm{u},\bm{v}_h)|\nonumber\\
\leq&  |\nu\langle \big(\nabla\bm u  -\bm{\Pi}_{m-1}^\ast\nabla\bm u\big) \bm{n},\bm{v}_{hi}-\bm{v}_{hb} \rangle_{\partial\mathcal{T}_h}|  +|\nu\langle \eta(\bm{P}^{RT}_m\bm{u}-{\bm{u}}),\bm{\Pi}_{k}^{B}\bm{v}_{hi}-\bm{v}_{hb} \rangle_{\partial\mathcal{T}_h}|\\
\lesssim &\sum_{K\in \mathcal{T}_{h}}\|\nabla\bm u  -\bm{\Pi}_{m-1}^\ast\nabla\bm u\|_{0,\partial K}
(\|\bm{v}_{hi}-\bm{\Pi}_{k}^{B}\bm{v}_{hi}\|_{0,\partial K}+\|\bm{\Pi}_{k}^{B}\bm{v}_{hi}-\bm{v}_{hb}\|_{0,\partial K} )\nonumber\\
&+\sum_{K\in \mathcal{T}_{h}}\|\eta^{\frac{1}{2}}(\bm{P}^{RT}_m\bm{u}-{\bm{u}} )\|_{0,\partial K}
\|\eta^{\frac{1}{2}}(\bm{\Pi}_{k}^{B}\bm{v}_{hi}-\bm{v}_{hb}  )\|_{0,\partial K}\nonumber\\
\lesssim &h^{m}\|\bm{u}\|_{m+1}(\| \nabla_{h}\bm{v}_{hi}\|_{0} +\|\eta^{\frac{1}{2}}(\bm{\Pi}_{k}^{B}\bm{v}_{hi}-\bm{v}_{hb}  )\|_{0,\partial K})
+h^{m}\|\bm{u}\|_{m+1}|||\bm{v}_{h}|||_{V}\nonumber\\
\lesssim &h^{m}\|\bm{u}\|_{m+1}|||\bm{v}_{h}|||_{V}
\end{align*}
and 
\begin{align*}
    |\xi_{III}(\bm{u},\bm{u};\bm{v}_h)| 
    \lesssim& \alpha\|\bm{P}^{RT}_m\bm{u}-\bm{u}\|_{0,3}(\|\bm{P}^{RT}_m\bm{u} \|_{0,2(r-2)}^{r-2}+\|\bm{u}\|_{0,2(r-2)}^{r-2})\|\bm{v}_{hi} \|_{0,6}\nonumber\\
    \lesssim & h^{m}\|\bm{u}\|_{2}^{r-2}\|\bm{u}\|_{m+1}|||\bm{v}_h|||_{V},
	\end{align*}
where in the estimate of $|\xi_{III}|$ we have used Lemma \ref{Lemma 2.9}. This finishes the proof.
\end{proof}
Based on Lemmas \ref{Lemma 4.1} and  \ref{Lemma 4.2}, we can obtain the following conclusion.
\begin{theorem} \label{Theorem 4.1}
Let $(\bm{u}_{h},p_h)\in \bm{V}_h^{0}\times Q^0_h$  be the solution to the WG scheme $(\ref{WG})$.
 Under the regularity assumption \eqref{regularity} and the discrete smallness condition  \eqref{uni-condi} with
 \begin{eqnarray}\label{UNI}
  \vartheta:=
  1-\frac{\mathcal{N}_h}{\nu^{2}}\|\bm{f}\|_{*,h}>0,
\end{eqnarray}
 there hold the following error estimates:
\begin{subequations}
\begin{align}
	||| \bm{\mathcal{I}}_{h}\bm{u}-\bm{u}_h|||_{V} \lesssim&  \mathcal{M}_{1}( \bm{u})h^{m},\label{78}\\
	|||\mathcal{P}_hp-p_h|||_{Q}\lesssim& \mathcal{M}_{2}( \bm{u})h^{m}+\mathcal{M}_{3}( \bm{u})h^{2m}  ,\label{79}
	\end{align}
\end{subequations}
where  $\mathcal{M}_{1}( \bm{u}):= \vartheta^{-1} 
(1+\|\bm{u}\|_{2}+\|\bm{u}\|_{2}^{r-2})\|\bm{u}\|_{m+1}$, and  $\mathcal{M}_{2}( \bm{u})$ and  $\mathcal{M}_{3}( \bm{u})$ are two positive constants depending only on $\vartheta, \nu, \|\bm{f}\|_{*,h} $, $\|\bm{u}\|_{2}$ and $\|\bm{u}\|_{m+1}$.
\end{theorem}
\begin{proof}
    Subtracting \eqref{WG1} and \eqref{WG2} from \eqref{e2} and \eqref{e3}, respectively,
    we have
\begin{subequations}\label{ee}
\begin{align}
&a_h(\bm{\mathcal{I}}_{h}\bm{u}-\bm{u}_{h},\bm{v}_h)+b_h(\bm{v}_h,\mathcal{P}_{h}p-p_h)
+c_h(\bm{\mathcal{I}}_{h}\bm{u};\bm{\mathcal{I}}_{h}\bm{u},\bm{v}_h)\nonumber\\
&\qquad -c_h(\bm{u}_{h};\bm{u}_{h},\bm{v}_{h})
+ d_h(\bm{\mathcal{I}}_{h}\bm{u};\bm{\mathcal{I}}_{h}\bm{u},\bm{v}_h)
-d_h(\bm{u}_{h};\bm{u}_{h},\bm{v}_{h})&\nonumber\\
 &\quad=\xi_{I}(\bm{u};\bm{u},\bm{v}_h)+\xi_{II}(\bm{u},\bm{v}_h)+\xi_{III}(\bm{u};\bm{u},\bm{v}_h),\forall \bm{v}_{h}\in \bm{V}_{h}^{0} .\label{ee1}\\
	&b_h(\bm{\mathcal{I}}_{h}\bm{u}-\bm{u}_{h},q_h)=0,\forall q_h\in Q_{h}^{0} .\label{ee2}
\end{align}
\end{subequations}
 Taking $\bm{v}_h=\bm{\mathcal{I}}_{h}\bm{u}-\bm{u}_{h} $  
 in equation \eqref{ee1}  and utilizing Lemmas  \ref{Lemma 2.9} and \ref{Lemma 3.1}, we obtain
	\begin{align*}
&\nu|||\bm{\mathcal{I}}_{h}\bm{u}-\bm{u}_{h} |||^{2}_{V}+\|\bm{P}_{m}^{RT}\bm{u}-\bm{u}_{hi} \|_{0,r}^{^{r}}\nonumber\\
\lesssim
&\nu|||\bm{\mathcal{I}}_{h}\bm{u}-\bm{u}_{h} |||^{2}_{V}
+(|\bm{P}_{m}^{RT}\bm{u}|^{r-2}\bm{P}_{m}^{RT}\bm{u}-|\bm{u}_{hi}|^{r-2}\bm{u}_{hi},\bm{P}_{m}^{RT}\bm{u}-\bm{u}_{hi} )\nonumber\\
=
&\xi_{I}(\bm{u},\bm{u};\bm{\mathcal{I}}_{h}\bm{u}-\bm{u}_{h})+\xi_{II}(\bm{u}; \bm{\mathcal{I}}_{h}\bm{u}-\bm{u}_{h})
+\xi_{III}(\bm{u},\bm{u};\bm{\mathcal{I}}_{h}\bm{u}-\bm{u}_{h})\nonumber\\
    &-\left\{d_h(\bm{\mathcal{I}}_{h}\bm{u};\bm{\mathcal{I}}_{h}\bm{u},\bm{\mathcal{I}}_{h}\bm{u}-\bm{u}_{h})
-d_h(\bm{u}_{h};\bm{u}_{h},\bm{\mathcal{I}}_{h}\bm{u}-\bm{u}_{h} )\right\}\nonumber\\
=&\xi_{I}(\bm{u},\bm{u}; \bm{\mathcal{I}}_{h}\bm{u}-\bm{u}_{h})
+\xi_{II}(\bm{u}; \bm{\mathcal{I}}_{h}\bm{u}-\bm{u}_{h})
+\xi_{III}(\bm{u},\bm{u};\bm{\mathcal{I}}_{h}\bm{u}-\bm{u}_{h})\nonumber\\
&-d_h(\bm{\mathcal{I}}_{h}\bm{u}-\bm{u}_{h} ;\bm{u}_h,\bm{\mathcal{I}}_{h}\bm{u}-\bm{u}_{h} ),
	\end{align*}
where in the  last $``="$ we have used the relation $d_h( \bm{\mathcal{I}}_{h}\bm{u}; \bm{\mathcal{I}}_{h}\bm{u}-\bm{u}_{h} ,\bm{\mathcal{I}}_{h}\bm{u}-\bm{u}_{h} )=0$.
In view of Lemma \ref{Lemma 4.2},   the definition of $\mathcal{N}_{h}$ with $\bm{u}_{h} , \bm{\mathcal{I}}_{h}\bm{u}\in  \bm{V}_{0h}$,   and the fact  that    $\|\bm{P}_{m}^{RT}\bm{u}-\bm{u}_{hi} \|_{0,r}^{r}\geq 0$,  we further obtain
\begin{align*}
\nu|||\bm{\mathcal{I}}_{h}\bm{u}-\bm{u}_{h} |||^{2}_{V}
\leq& C h^{m}\left(  \|\bm{u}\|_{2}\|\bm{u}\|_{ m+1}
+   \|\bm{u}\|_{m+1} +  \|\bm{u}\|_{2}^{r-2}\|\bm{u}\|_{m+1}\right)|||\bm{\mathcal{I}}_{h}\bm{u}-\bm{u}_{h}|||_{V}
\nonumber\\
&\quad +\mathcal{N}_h|||\bm{u}_h|||_{V} |||\bm{\mathcal{I}}_{h}\bm{u}-\bm{u}_{h}  |||^2_{V} ,
\end{align*}
which, together with \eqref{sta}, yields
\begin{align*}
&\nu(1-\frac{\mathcal{N}_h}{\nu^{2}}\|\bm{f}\|_{*,h})|||\bm{\mathcal{I}}_{h}\bm{u}-\bm{u}_{h} |||_{V}\leq  C h^{m}(1+ \|\bm{u}\|_{2}+ \|\bm{u}\|_{2}^{r-2})\|\bm{u}\|_{m+1}.
\end{align*}
Thus,  the desired estimate \eqref{78} follows.

Next we estimate the error of pressure.
Using the first equation in \eqref{ee} we get
\begin{align*}
&b_h(\bm{v}_{h},  \mathcal{P}_{h}p-p_h )\nonumber\\
=&-a_h(\bm{\mathcal{I}}_{h}\bm{u}-\bm{u}_{h},\bm{v}_h)
-c_h(\bm{\mathcal{I}}_{h}\bm{u};\bm{\mathcal{I}}_{h}\bm{u},\bm{v}_h)\nonumber\\
&+c_h(\bm{u}_{h};\bm{u}_{h},\bm{v}_{h})
- d_h(\bm{\mathcal{I}}_{h}\bm{u};\bm{\mathcal{I}}_{h}\bm{u},\bm{v}_h)
+d_h(\bm{u}_{h};\bm{u}_{h},\bm{v}_{h})\nonumber\\
&+\xi_{I}(\bm{u};\bm{u},\bm{v}_h)+\xi_{II}(\bm{u},\bm{v}_h)+\xi_{III}(\bm{u};\bm{u},\bm{v}_h),\forall \bm{v}_{h}\in \bm{V}_{h}^{0}\nonumber\\
=&-a_h(\bm{\mathcal{I}}_{h}\bm{u}-\bm{u}_{h},\bm{v}_h)
-\big(c_h(\bm{\mathcal{I}}_{h}\bm{u};\bm{\mathcal{I}}_{h}\bm{u},\bm{v}_h)-c_h(\bm{u}_{h};\bm{u}_{h},\bm{v}_{h})\big)\nonumber\\
&- d_h(\bm{\mathcal{I}}_{h}\bm{u}-\bm{u}_{h};\bm{\mathcal{I}}_{h}\bm{u}-\bm{u}_{h},\bm{v}_h)
-d_h(\bm{\mathcal{I}}_{h}\bm{u}-\bm{u}_{h};\bm{u}_{h},\bm{v}_{h})
-d_h(\bm{u}_{h};\bm{\mathcal{I}}_{h}\bm{u}-\bm{u}_{h},\bm{v}_{h})\nonumber\\
&+\xi_{I}(\bm{u};\bm{u},\bm{v}_h)+\xi_{II}(\bm{u},\bm{v}_h)+\xi_{III}(\bm{u};\bm{u},\bm{v}_h),\forall \bm{v}_{h}\in \bm{V}_{h}^{0}.
\end{align*}
In light  of   Lemmas  \ref{Lemma 2.4}, \ref{Lemma 2.9}, \ref{Lemma 3.1}, \ref{theoremLBB}, and \ref{Lemma 4.2} and the estimates \eqref{sta} and \eqref{78}, we have
	\begin{align*}
	|||\mathcal{P}_{h}p-p_h|||_{Q}
    \lesssim&\sup_{\bm 0\neq \bm{v}_{h}\in \bm{V}_h}\frac{b_h(\bm{v}_{h}, \mathcal{P}_{h}p-p_h)}{|||\bm{v}_h|||_V}\nonumber\\
    \lesssim& \nu|||\bm{\mathcal{I}}_{h}\bm{u}-\bm{u}_{h}|||_{V}+\alpha C_{r}C_{\widetilde{r}}^{r} ( |||\bm{\mathcal{I}}_{h}\bm{u}|||_{V}+|||\bm{u}_{h}|||_{V} )^{r-2}|||\bm{\mathcal{I}}_{h}\bm{u}-\bm{u}_{h}|||_{V}\\
    &+|||\bm{\mathcal{I}}_{h}\bm{u}-\bm{u}_{h}|||_{V}^{2}+|||\bm{u}_{h}|||_{V}|||\bm{\mathcal{I}}_{h}\bm{u}-\bm{u}_{h}|||_{V}\\
    &+\xi_{I}(\bm{u};\bm{u},\bm{v}_h)+\xi_{II}(\bm{u},\bm{v}_h)+\xi_{III}(\bm{u};\bm{u},\bm{v}_h)\\
   \lesssim&\left(\nu +  \frac{\|\bm{f}\|_{*,h}}{\nu}+\frac{\|\bm{f}\|_{*,h}^{r-2}}{\nu^{r-2}}+\|\bm{u}\|_{2}^{r-2}\right) |||\bm{\mathcal{I}}_{h}\bm{u}-\bm{u}_{h}|||_{V}+|||\bm{\mathcal{I}}_{h}\bm{u}-\bm{u}_{h}|||_{V}^{2}\nonumber\\
    &+h^{m}(1+ \|\bm{u}\|_{2}+ \|\bm{u}\|_{2}^{r-2})\|\bm{u}\|_{m+1},
	\end{align*}
	which shows \eqref{79}.
\end{proof}

Finally, based upon    Theorem \ref{Theorem 4.1},  Lemmas \ref{Lemma 2.2}  and \ref{Lemma 2.6} - \ref{Lemma 2.8}, we can obtain  the following main conclusion.
\begin{theorem} \label{Theorem 4.2}
 Under the same conditions as in Theorem \ref{Theorem 4.1},  there hold
	\begin{align}
	\|\nabla\bm{u}-\nabla_h\bm{u}_{hi}\|_0
+\|\nabla\bm{u}-\nabla_{w,m-1}\bm{u}_{h}\|_0&\lesssim \mathcal{M}_{1}(\bm{u} ) h^{m}\label{HU},\\
	\|p-p_{hi}\|_0&\lesssim(\mathcal{M}_{2}(\bm{u} )+\|p\|_{m})h^{m}+\mathcal{M}_{3}( \bm{u})h^{2m}. \label{LP}
	\end{align}
\end{theorem}
\begin{remark}
	 The result \eqref{HU} shows  that the velocity error estimate is independent of the pressure approximation, which means that the proposed WG scheme is pressure-robust.
\end{remark}

\section{$L^2$ error estimation for velocity }
\label{section6.2}
We follow standard dual arguments to derive an   $L^2$ error estimate for the  velocity solution of the WG scheme.  To this end, we introduce the following dual problem:  seek $(\bm{\phi},\psi)$ such that
\begin{eqnarray}\label{l2}
\left\{
\begin{aligned}
-\nu\triangle\bm{\phi}-(\bm{u}\cdot\nabla)\bm{\phi}+(\nabla\bm{u })^T\bm{\phi}+\alpha\mid\bm{u } \mid^{r-2}\bm{\phi}&\\
+\alpha(r-2)\mid\bm{u }\mid^{r-4}(\bm{u }\cdot\bm{\phi})\bm{u } +\nabla \psi=&\bm{e}_{hi}, \text{ in }\ \Omega,\label{l21} \\
\nabla\cdot \bm{\phi}=&0,\ \text{ in } \ \Omega, \label{l22} \\
\bm{\phi}=&\bm{0},\ \text{ on } \ \partial \Omega, \label{l23}\\
\end{aligned}
\right.
\end{eqnarray}
where $\bm{u}$ and $\bm{u}_{h}=\{\bm{u}_{hi},\bm{u}_{hb}\}$ are respectively  the solutions of (\ref{BF0}) and (\ref{WG}), and $\bm{e}_{hi}:= \bm{P}^{RT}_m\bm{u}-\bm{u}_{hi}.$
We assume  the following  regularity condition holds:
\begin{eqnarray}\label{sta-assum}
\|\bm{\phi}\|_{2}+\|\psi\|_{1}\lesssim\|\bm{e}_{hi}\|_0  .
\end{eqnarray}

The corresponding weak form of  \eqref{l2} reads: seek $(\bm{\phi},\psi) \in \bm{V} \times Q$ such that
\begin{subequations}\label{weakl22}
\begin{align}
 \mathcal{A}_u(\bm{\phi},\bm{v})+b(\bm{v},\psi )&=(\bm{e}_{hi},\bm{v}), \quad &\forall \bm{v} \in \bm{V} , \\
b(\bm{\phi},q )&=0, \quad & \forall q \in Q,
\end{align}
\end{subequations}
where  the bilinear form $\mathcal{A}_u(\cdot,\cdot): \bm{V}\times\bm{V}\rightarrow \mathbb{R}$ is defined by
\begin{align}\label{A}
 \mathcal{A}_u(\bm{\phi},\bm{v}):=&a(\bm{\phi},\bm{v})+c(\bm{u};\bm{\phi},\bm{v} )+d(\bm{u};\bm{v},\bm{\phi} )+d(\bm{v};\bm{u},\bm{\phi} )\nonumber\\
 &+\alpha(r-2)(|\bm{u }|^{r-4}(\bm{u }\cdot\bm{\phi})\bm{u } ,\bm{v}),
\end{align}
and  the bilinear forms, $a(\cdot,\cdot)$ and $ b(\cdot,\cdot)$, and the trilinear forms, $ c(\cdot;\cdot,\cdot ) $ and $d(\cdot;\cdot,\cdot)$, are given in subsection 2.1.

\begin{remark}
According to the H\"{o}lder's inequality, the Sobolev inequality and the boundedness result \eqref{boundedCondition},
we can get the boundedness result
\begin{align}
\mathcal{A}_u(\bm{\phi},\bm{v})\lesssim \|\nabla \bm{\phi}\|_{0}\| \nabla\bm{v}\|_{0}, \quad \forall  \bm{\phi},\bm{v}\in \bm{V}. \label{A1}
\end{align}
At the same time, under the uniqueness condition \eqref{Condition} we can  obtain the coercivity result
  \begin{align}
 \mathcal{A}_u(\bm{v},\bm{v})&\gtrsim  \| \nabla \bm{v}\|_{0}^{2}, \quad \forall \bm{v}\in \bm{V}. \label{A2}
\end{align}
 It is standard that the inf-sup inequality 
	\begin{eqnarray*}
		\sup_{ \bm{v}\in \bm{V}}\frac{b(\bm{v},q)}{||\nabla\bm{v}||_{0} }\gtrsim\|q \|_{0}, \quad \forall q \in Q
	\end{eqnarray*}
	holds. As a result, the problem \eqref{weakl22} admits a unique solution.
 \end{remark}

By taking   similar routines as in the proofs of Lemmas \ref{Lemma 4.1} and \ref{Lemma 4.2},  respectively, we can obtain Lemmas \ref{Lemma 51} and \ref{Lemma 52}.
\begin{lemma}\label{Lemma 51}
There hold
\begin{subequations}\label{L-2}
\begin{align}
&a_h(\bm{\mathcal{I}}_{h}\bm{\phi},\bm{v}_h)
+b_h(\bm{v}_{h},\mathcal{P}_{h}\psi)
+c_h(\bm{\mathcal{I}}_{h}\bm{u};\bm{\mathcal{I}}_{h}\bm{\phi},\bm{v}_h)
- d_h(\bm{\mathcal{I}}_{h}\bm{u};\bm{\mathcal{I}}_{h}\bm{\phi},\bm{v}_h)
+d_h(\bm{v}_h;\bm{\mathcal{I}}_{h}\bm{u},\bm{\mathcal{I}}_{h}\bm{\phi})\nonumber\\
 &\quad =(\bm{e}_{hi},\bm{v}_{hi}) -E_{I}(\bm{u};\bm{\phi},\bm{v}_{h})
 +E_{II}(\bm{\phi},\bm{v}_h)+E_{III}(\bm{u};\bm{\phi},\bm{v}_h)
  +E_{IV}(\bm{v}_h;\bm{u},\bm{\phi})\nonumber\\
 & \qquad
       -\alpha(r-2)(\mid\bm{u }\mid^{r-4}(\bm{u }\cdot\bm{\phi})\bm{u },\bm{e}_{hi}),\quad \forall \bm{v}_h \in \bm{V}_{h}^{0}, \label{error-2}\\
	&b_h(\bm{\mathcal{I}}_{h}\bm{\phi},q_h)=0,\quad \forall  q_h \in Q_{h}^{0},\label{error-3}
\end{align}
\end{subequations}
where
\begin{align*}
&E_{I}(\bm{u};\bm{\phi},\bm{v}_{h}):= -\frac{1}{2}(\bm{P}_m^{RT}\bm{\phi}\otimes \bm{P}_m^{RT}\bm{u}-\bm{\phi}\otimes\bm{u},\nabla_h\bm{v}_{hi}) 
+\frac{1}{2}\langle(\bm{\Pi}_{k}^{B}\bm{\phi}\otimes\bm{\Pi}^{B}_{k}\bm{u}
-\bm{\phi}\otimes\bm{u})\bm{n},\bm{v}_{hi}\rangle_{\partial\mathcal{T}_h}\nonumber\\
&\qquad\hskip2cm -\frac{1}{2}((\bm{u}\cdot\nabla)\bm{\phi}-(\bm{P}_m^{RT}\bm{ u}\cdot\nabla_h)\bm{P}_m^{RT}\bm{\phi},\bm{v}_{hi})
 -\frac{1}{2}\langle(\bm{v}_{hb}\otimes \bm{\Pi}_{k}^{B}\bm{u})\bm{n},\bm{P}_m^{RT}\bm{\phi}\rangle_{\partial\mathcal{T}_h},\\
	&E_{II}(\bm{\phi},\bm{v}_h):= \nu\langle \big(\nabla\bm \phi  -\bm{\Pi}_{m-1}^\ast\nabla\bm \phi\big)\bm{n},\bm{v}_{hi}-\bm{v}_{hb} \rangle_{\partial\mathcal{T}_h} +\nu\langle \eta(\bm{P}^{RT}_m\bm{\phi}-{\bm{\phi}}),\bm{\Pi}_{k}^{B}\bm{v}_{hi}-\bm{v}_{hb} \rangle_{\partial\mathcal{T}_h},\nonumber\\
&E_{III}(\bm{u};\bm{\phi},\bm{v}_h)=
\alpha (|\bm{P}^{RT}_m\bm{u}|^{r-2}\bm{P}^{RT}_m\bm{\phi}-|\bm{u}|^{r-2}\bm{\phi},\bm{v}_{hi}),\nonumber\\
&	E_{IV}(\bm{v}_h;\bm{u},\bm{\phi}):= \big( (\nabla_h\bm{P}_m^{RT}\bm{u})^T\bm{P}_m^{RT}\bm{\phi}-(\nabla\bm{u})^T\bm{\phi},\bm{v}_{hi}\big)
-\frac{1}{2}\langle(\bm{\Pi}_k^{B}\bm{\phi}\otimes\bm{v}_{hb}) \bm{n},\bm{P}_m^{RT}\bm{u}\rangle_{\partial\mathcal{T}_h}\nonumber\\
	&\qquad\hskip2cm+\frac{1}{2}\langle (\bm{\Pi}_k^{B}\bm{u}\otimes\bm{v}_{hb})\bm{n},\bm{P}_m^{RT}\bm{\phi}\rangle_{\partial\mathcal{T}_h}
	-\frac{1}{2}\langle(\bm{P}_m^{RT} \bm{u}\otimes\bm{v}_{hi})\bm{n},\bm{P}_m^{RT}\bm{\phi}\rangle_{\partial\mathcal{T}_h}.
	\end{align*}
\end{lemma}

\begin{lemma}
\label{Lemma 52}
For any $\bm{v}_h\in\bm{V}_{h}^0$, there hold
\begin{subequations}
	\begin{align}
	|\xi_{I}(\bm{u};\bm{u},\bm{\mathcal{I}}_{h}\bm{\phi})|\lesssim& h^{m+1}
	\|\bm{u}\|_{2}\|\bm{u}\|_{m+1}\|\bm{\phi}\|_{2},\label{e58a}\\
|E_{I}(\bm{u};\bm{\phi},\bm{v}_h)|\lesssim& h
	\|\bm{u}\|_{2}\|\bm{\phi}\|_{2}|||\bm{v}_h|||_{V},\label{e58b}\\
	|\xi_{II}(\bm{u},\bm{\mathcal{I}}_{h}\bm{\phi})|+
	|E_{II}(\bm{\phi},\bm{e}_{h})|\lesssim& h^{m+1}\|\bm{u}\|_{m+1}\|\bm{\phi}\|_{2},\label{e58c}\\
|\xi_{III}(\bm{u};\bm{u},\bm{\mathcal{I}}_{h}\bm{\phi})|\lesssim& h^{m+1}
	\|\bm{u}\|_{2}^{r-2}\|\bm{u}\|_{m+1}\|\bm{\phi}\|_{2},\label{e58d}\\
|E_{III}(\bm{u};\bm{\phi},\bm{v}_{h})|\lesssim& h
	\|\bm{u}\|_{2}^{r-2}\|\bm{\phi}\|_{2}|||\bm{v}_h|||_{V},\label{e58e}\\
    |E_{IV}(\bm{v}_h;\bm{u},\bm{\phi})|\lesssim& h\|\bm{u}\|_{2}\|\bm{\phi}\|_{2}|||\bm{v}_h|||_{V}\label{e58f}.
	\end{align}
\end{subequations}
\end{lemma}

\begin{theorem}\label{Theorem 5.1}
Under the regularity condition \eqref{sta-assum} and the same conditions as in Theorem \ref{Theorem 4.1}, there holds
\begin{align}
&\|\bm{u}-\bm{u}_{h}\|_{0}\lesssim \mathcal{M}_{4}(\bm{u}) h^{m+1},\label{LU}
\end{align}
where $\mathcal{M}_{4}(\bm{u})$ is a positive constant depending on $\|\bm{u}\|_{2}$ and $\|\bm{u}\|_{m+1}$. 
\end{theorem}
\begin{proof}
Denote $\bm{e}_h:=\bm{\mathcal{I}}_{h}\bm{u}-\bm{u}_h$ and  $\varepsilon_h:=\mathcal{P}_{h}p-p_h$.
Taking $\bm{v}_h=\bm{e}_{h}$ and $q_h=\varepsilon_{h} $  in \eqref{L-2}, we derive
\begin{subequations}\label{L21}
\begin{align}
\|\bm{e}_{hi}\|_{0}^{2}=&a_h(\bm{\mathcal{I}}_{h}\bm{\phi},\bm{e}_h)
+b_h(\bm{e}_{h},\mathcal{P}_{h}\psi)
+c_h(\bm{\mathcal{I}}_{h}\bm{u};\bm{\mathcal{I}}_{h}\bm{\phi},\bm{e}_{h})\nonumber\\
&+\alpha(r-2)(\mid\bm{u }\mid^{r-4}(\bm{u }\cdot\bm{\phi})\bm{u },\bm{e}_{hi})
+d_h(\bm{\mathcal{I}}_{h}\bm{u};\bm{e}_h, \bm{\mathcal{I}}_{h}\bm{\phi})
+d_h(\bm{e}_h;\bm{\mathcal{I}}_{h}\bm{u}, \bm{\mathcal{I}}_{h}\bm{\phi})\nonumber\\&
+E_{I}(\bm{u};\bm{\phi},\bm{e}_h )-E_{II}(\bm{\phi},\bm{e}_h )
-E_{III}(\bm{u};\bm{\phi},e_h)
-E_{IV}(\bm{e}_h;\bm{u},\bm{\phi} ),\\
&\quad\quad\quad\quad\quad\quad\quad\quad\quad\quad\quad\quad b_h( \bm{\mathcal{I}}_{h}\bm{\phi},\varepsilon_h)=0.
\end{align}
\end{subequations}
Taking $\bm{v}_h=\bm{\mathcal{I}}_{h}\bm{\phi}$ and $q_{h}=\mathcal{P}_{h}\psi$ in \eqref{ee}, respectively, we have
\begin{align*}
a_h(\bm{e}_{h},\bm{\mathcal{I}}_{h}\bm{\phi})
+b_h(\bm{\mathcal{I}}_{h}\bm{\phi},\varepsilon_{h})
+c_h(\bm{\mathcal{I}}_{h}\bm{u};\bm{\mathcal{I}}_{h}\bm{u},\bm{\mathcal{I}}_{h}\bm{\phi})
-c_h(\bm{u}_h;\bm{u}_h,\bm{\mathcal{I}}_{h}\bm{\phi})&\nonumber\\
+d_h(\bm{\mathcal{I}}_{h}\bm{u};\bm{\mathcal{I}}_{h}\bm{u},\bm{\mathcal{I}}_{h}\bm{\phi})
-d_h(\bm{u}_h;\bm{u}_h,\bm{\mathcal{I}}_{h}\bm{\phi})&\nonumber\\
=\xi_{I}(\bm{u},\bm{u},\bm{\mathcal{I}}_{h}\bm{\phi})+\xi_{II}(\bm{u};\bm{\mathcal{I}}_{h}\bm{\phi})+\xi_{III}(\bm{u},\bm{u}; \bm{\mathcal{I}}_{h}\bm{\phi}),& \\
b_h(\bm{e}_{h},\mathcal{P}_{h}\psi)=0,&
\end{align*}
which plus \eqref{L21}  give
\begin{align}\label{ine511}
\|\bm{e}_{hi}\|_{0}^{2}=&\big\{c_h(\bm{\mathcal{I}}_{h}\bm{u};\bm{\mathcal{I}}_{h}\bm{\phi},\bm{e}_{h})
-c_h(\bm{\mathcal{I}}_{h}\bm{u};\bm{\mathcal{I}}_{h}\bm{u},\bm{\mathcal{I}}_{h}\bm{\phi})
+c_h(\bm{u}_h;\bm{u}_h,\bm{\mathcal{I}}_{h}\bm{\phi})
\nonumber\\
&+(r-2)(\mid\bm{u }\mid^{r-4}(\bm{u }\cdot\bm{\phi})\bm{u },\bm{e}_{hi})\big\}
+\big\{-d_h(\bm{\mathcal{I}}_{h}\bm{u};\bm{\mathcal{I}}_{h}\bm{u},\bm{\mathcal{I}}_{h}\bm{\phi})
		+d_h(\bm{u}_h;\bm{u}_h,\bm{\mathcal{I}}_{h}\bm{\phi})\nonumber\\
&\qquad +d_h(\bm{\mathcal{I}}_{h}\bm{u};\bm{e}_{h}, \bm{\mathcal{I}}_{h}\bm{\phi})
+d_h(\bm{e}_{h};\bm{\mathcal{I}}_{h}\bm{u}, \bm{\mathcal{I}}_{h}\bm{\phi})\big\}
\nonumber\\
&+\big\{\xi_{I}(\bm{u},\bm{u},\bm{\mathcal{I}}_{h}\bm{\phi})+E_{I}(\bm{u};\bm{\phi},\bm{e}_{h} )
+\xi_{II}(\bm{u};\bm{\mathcal{I}}_{h}\bm{\phi})-E_{II}(\bm{\phi},\bm{e}_{h} )\nonumber\\&
\qquad +\xi_{III}(\bm{u},\bm{u}; \bm{\mathcal{I}}_{h}\bm{\phi})
-E_{III}(\bm{u};\bm{\phi},e_h)
-E_{IV}(\bm{e}_{h};\bm{u},\bm{\phi} )\big\}\nonumber\\
=: &\sum_{j=1}^{3}\mathfrak{R}_{j}.
\end{align}
Then let us estimate $\mathfrak{R}_{j}(j=1,...,3)$ one by one  by  using Lemmas \ref{Lemma 2.1}, \ref{Lemma 2.2}, \ref{Lemma 2.4}, \ref{Lemma 2.6}, \ref{Lemma 2.7}, \ref{Lemma 3.1} and \ref{Lemma 52}, and we have
\begin{align*}
|\mathfrak{R}_{1}|=&|\alpha\{-(|\bm{P}^{RT}_m\bm{u}|^{r-2}\bm{P}^{RT}_m\bm{u}, \bm{P}^{RT}_m\bm{\phi} )+(|\bm{u}_{h}|^{r-2}\bm{u}_{h}, \bm{P}^{RT}_m\bm{\phi} ) +(\mid\bm{P}^{RT}_m\bm{u } \mid^{r-2}\bm{P}^{RT}_m\bm{\phi},\bm{e}_{hi})
\nonumber\\
&+(r-2)(\mid\bm{u }\mid^{r-4}(\bm{u }\cdot\bm{\phi})\bm{u },\bm{e}_{hi})\}|\nonumber\\
=&|\alpha\{(|\bm{P}^{RT}_m\bm{u}|^{r-2}\bm{e}_{hi},\bm{\phi}-\bm{P}^{RT}_m\bm{\phi})+\big((\bm{P}^{RT}_m\bm{u}|^{r-2}-|\bm{u}_{h}|^{r-2})\bm{e}_{hi},\bm{P}^{RT}_m\bm{\phi}\big)\nonumber\\
&+\big((|\bm{P}^{RT}_m\bm{u}|^{r-2}-|\bm{u}_{h}|^{r-2})\bm{P}^{RT}_m\bm{u},\bm{\phi}-\bm{P}^{RT}_m\bm{\phi}\big)\nonumber\\
&-\big((|\bm{P}^{RT}_m\bm{u}|^{r-2}-|\bm{u}_{h}|^{r-2}-(r-2)\mid\bm{u }\mid^{r-4}\bm{u }\cdot \bm{e}_{hi},\bm{u}\cdot\bm{\phi}\big)\}|\nonumber\\
\lesssim&\|\bm{u}\|^{r-2}_{0,3(r-2)}\|\bm{e}_{hi}\|_{0,6}\|\bm{\phi}-\bm{P}^{RT}_m\bm{\phi}\|_{0}
+(\|\bm{u}\|^{r-3}_{0,3(r-3)}+\|\bm{u}_{h}\|^{r-3}_{0,3(r-3)})\|\bm{e}_{hi}\|_{0,6}^{2}\|\bm{P}^{RT}_m\bm{\phi}\|_{0,3}\nonumber\\
&+(\|\bm{u}\|^{r-3}_{0,3(r-3)}+\|\bm{u}_{h}\|^{r-3}_{0,3(r-3)})\|\bm{e}_{hi}\|_{0,6}\|\bm{P}^{RT}_m\bm{u}\|_{0,6}\|\bm{\phi}-\bm{P}^{RT}_m\bm{\phi}\|_{0,3}\nonumber\\
&+(\|\bm{u}\|^{r-4}_{0,3(r-4)}+\|\bm{u}_{h}\|^{r-4}_{0,3(r-4)})\|\bm{e}_{hi}\|_{0,6}^{2}\|\bm{u}\|_{0,6}\|\bm{\phi}\|_{0,3}\nonumber\\
\lesssim& h^{2-\frac{n}{6}}|||\bm{e}_{h}|||_{V} \cdot\|\bm{\phi}  \|_{2} ,\\
|\mathfrak{R}_{2}|=& |-d_h(\bm{\mathcal{I}}_{h}\bm{u};\bm{\mathcal{I}}_{h}\bm{u},\bm{\mathcal{I}}_{h}\bm{\phi})
		+d_h(\bm{u}_h;\bm{u}_h,\bm{\mathcal{I}}_{h}\bm{\phi})
+d_h(\bm{\mathcal{I}}_{h}\bm{u};\bm{e}_{h}, \bm{\mathcal{I}}_{h}\bm{\phi})
+d_h(\bm{e}_{h};\bm{\mathcal{I}}_{h}\bm{u}, \bm{\mathcal{I}}_{h}\bm{\phi})|\nonumber\\
=&|d_h(\bm{e}_{h};\bm{e}_{h}, \bm{\mathcal{I}}_{h}\bm{\phi})| \nonumber\\
\lesssim& |||\bm{e}_{h}|||^{2}_{V}||| \bm{\mathcal{I}}_{h}\bm{\phi} |||_{V}\nonumber\\
 \lesssim&h|||\bm{e}_{h}|||^{2}_{V}\| \bm{\phi} \|_{2}, \\
|\mathfrak{R}_{3}|=&|\xi_{I}(\bm{u},\bm{u},\bm{\mathcal{I}}_{h}\bm{\phi})+E_{I}(\bm{u};\bm{\phi},\bm{e}_{h} )
+\xi_{II}(\bm{u};\bm{\mathcal{I}}_{h}\bm{\phi})-E_{II}(\bm{\phi},\bm{e}_{h} )\nonumber\\&
+\xi_{III}(\bm{u},\bm{u}; \bm{\mathcal{I}}_{h}\bm{\phi})
-E_{III}(\bm{u};\bm{\phi},e_h)
-E_{IV}(\bm{e}_{h};\bm{u},\bm{\phi} )|\nonumber\\
\lesssim&   h^{m+1}\|\bm{u}\|_{m+1}\| \bm{u}\|_{2}\|\bm{\phi}\|_{2}
+h\|\bm{u}\|_{2}\|\bm{\phi}\|_{2}|||\bm{e}_{h}|||_{V}+h^{m+1}\|\bm{u}\|_{m+1}\|\bm{\phi}\|_{2}\nonumber\\
& +h^{m+1}\|\bm{u}\|_{m+1}\| \bm{u}\|_{2}^{r-2}\|\bm{\phi}\|_{2}
+h\| \bm{u}\|_{2}^{r-2}\|\bm{\phi}\|_{2}|||\bm{e}_{h}|||_{V}
+h \| \bm{u}\|_{2}\|\bm{\phi} \|_{2}|||\bm{e}_{h}|||_{V}.
\end{align*}
These three estimates, together with  \eqref{ine511}, \eqref{sta-assum},  Theorem \ref{Theorem 4.2} and the triangle inequality, yield  the desired conclusion \eqref{LU}.
\end{proof}

\section{Local elimination property and iteration scheme}

\subsection{Local elimination property}

In the subsection, we shall demonstrate that  in the WG scheme (\ref{WG})  the velocity and pressure approximations, $(\bm{u}_{hi}, p_{hi})$,
		defined in the interior of elements can be locally eliminated by the using the numerical traces $(\bm{u}_{hb}, p_{hb})$  defined  on the element interfaces.
		After the local elimination the resulting system  only includes the degrees of freedom of $(\bm{u}_{hb}, p_{hb})$ as unknowns.

For any $K\in \mathcal{T}_{h}$, taking $\bm{v}_{hi}|_{\mathcal{T}_{h}/K}= \bm{0}$, $\bm{v}_{hb}= \bm{0}$, $q_{hi}|_{\mathcal{T}_{h}/K}=0$ and $q_{hb}=0$  in \eqref{WG}, we obtain the following local problem: Seek $(\bm{u}_{hi},p_{hi})\in [P_{m}(K)]^{n}\times P_{m-1}(K)$ such that 
\begin{subequations}\label{WGhi}
\begin{align}
&a_{h,K}(\bm{u}_{hi},\bm{v}_{hi})+b_{h,K}(\bm{v}_{hi},p_{hi})+c_{h,K}(\bm{u}_{hi} ;\bm{u}_{hi},\bm{v}_{hi} )+d_{h,K}(\bm{u}_{hi} ;\bm{u}_{hi},\bm{v}_{hi} )\nonumber\\
& \qquad =F_{h,K}(\bm{v}_{hi}),\quad \forall \bm{v}_{hi} \in [P_{m}(K)]^{n},\label{WG-1}\\
	&b_{h,K}(\bm{u}_{hi},q_{hi})=0, \quad \forall  q_{hi}\in   P_{m-1}(K),\label{WG-2}
\end{align}
\end{subequations}
where
\begin{subequations}
\begin{align*}
&a_{h,K}(\bm{u}_{hi},\bm{v}_{hi}):= \nu (\nabla_{w,m-1}\{\bm{u}_{hi},\bm{0}\},\nabla_{w,m-1}\{\bm{v}_{hi},\bm{0}\}  ) _{K}+ s_{h,K}(\bm{u}_{hi},\bm{v}_{hi}),\\
&s_{h,K}(\bm{u}_{hi},\bm{v}_{hi}):=\nu \langle \eta \bm{\Pi}_{k}^{B}\bm{u}_{hi},\bm{\Pi}_{k}^{B}\bm{v}_{hi}\rangle_{\partial K},\\
&b_{h,K}(\bm{v}_{hi},p_{hi}):=(\nabla_{w,m}\{p_{hi},0\}, \bm{v}_{hi} )_{K},\\
&c_{h,K}(\bm{u}_{hi} ;\bm{u}_{hi},\bm{v}_{hi} ):= (\alpha|\bm{u}_{hi}|^{r-2}\bm{u}_{hi},\bm{v}_{hi}  )_{K},\\
&d_{h,K}(\bm{u}_{hi} ;\bm{u}_{hi},\bm{v}_{hi} ):=\frac{1}{2}(\nabla_{w,m}\cdot\{\bm{u}_{hi}\otimes \bm{u}_{hi} ,\bm{0}\otimes\bm{0}\},\bm{v}_{hi} )_{K}\nonumber\\
&\qquad \qquad\qquad\qquad \qquad -\frac{1}{2}(\nabla_{w,m}\cdot\{\bm{v}_{hi}\otimes \bm{u}_{hi} ,\bm{0}\otimes\bm{0}\},\bm{u}_{hi})_{K},\\
&F_{h,K}(\bm{v}_{hi}):= (\bm{f},\bm{v}_{hi})_{K}
- \nu (\nabla_{w,m-1}\{\bm{0},\bm{u}_{hb}\},\nabla_{w,m-1}\{\bm{v}_{hi},\bm{0}\}  ) _{K}+\nu \langle \eta \bm{u}_{hb},\bm{\Pi}_{k}^{B}\bm{v}_{hi}\rangle_{\partial K}\nonumber\\
& \qquad\qquad\qquad\qquad
-\frac{1}{2}(\nabla_{w,m}\cdot\{\bm{0}\otimes\bm{0},\bm{u}_{hb}\otimes \bm{u}_{hb} \},\bm{v}_{hi} )_{K}- (\nabla_{w,m}\{0,p_{hb}\}, \bm{v}_{hi} )_{K}.
\end{align*}
\end{subequations}

By following the same routines as in the proofs of Theorems \ref{TH42} and \ref{Theorem 4.3},
 we can get  existence and uniqueness results of  \eqref{WGhi}.
\begin{theorem}\label{th71}
For all $K\in \mathcal{T}_{h}$ and given numerical traces $\bm{u}_{hb}|_{\partial K}$ and $p_{hb}|_{\partial K} $, the local problem \eqref{WGhi} admits at least one solution. In addition, under the  smallness condition
\begin{align}
\frac{\mathcal{N}_{h,K}\|F_{h,K}\|_{\ast,h}}{\nu^{2}}<1,
\end{align}
the problem \eqref{WGhi} admits a unique solution.
Here
\begin{align*}
&\mathcal{N}_{h,K}:=\sup_{\bm 0\neq\bm{\kappa}_{hi},\bm{u}_{hi},\bm{v}_{hi}\in\bm{V}_{0h,K}}
\frac{d_{h,K}(\bm{\kappa}_{hi};\bm{u}_{hi},\bm{v}_{hi})}{|||\bm{\kappa}_{hi}|||_{V,K}\cdot |||\bm{u}_{hi}|||_{V,K}\cdot |||\bm{v}_{hi}|||_{V,K} },\label{NhK}\\
&\|F_{h,K}\|_{*,h}:=  \sup_{\bm 0\neq \bm{v}_{hi}\in\bm{V}_{0h,K}}\frac{F_{h,K}(\bm v_{hi})}{|||\bm{v}_h|||_{V,K}},\\
&\bm{V}_{0h,K}:=\{\bm{\kappa}_{hi}\in [P_{m} (K)]^{n}: b_{h,K}(\bm{\kappa}_{hi}, q_{hi} )=0, \forall q_{hi}\in P_{m-1} (K)\},\\
&||| \bm{v}_{hi}|||_{V,K}:=(\|\nabla_{w,m-1}\{\bm{v}_{hi},\bm{0}\}  \|_{0,K}^{2}
+\| \eta^{\frac{1}{2}}\bm{\Pi}_{k}^{B} \bm{v}_{hi}\|_{0,\partial K}^{2})^{\frac{1}{2}}.
\end{align*}
\end{theorem}

\subsection{Iteration scheme}

Due to the nonlinearity of the WG scheme \eqref{WG}, we shall employ  the following  Oseen's iteration algorithm:

Given $\bm{u}_{h}^{0}$, seek $(\bm{u}_{h}^{l},p_{h}^{l})$ with $l=1,2,...$,  such that
\begin{subequations}\label{Oseen}
\begin{align}
a_{h}(\bm{u}_{h}^{l},\bm{v}_{h})+b_{h}(\bm{v}_{h},p_{h}^{l})+c_{h}(\bm{u}_{h}^{l-1} ;\bm{u}_{h}^{l},\bm{v}_{h} )+d_{h}(\bm{u}_{h}^{l-1} ;\bm{u}_{h}^{l},\bm{v}_{h} )
 =&(\bm{f },\bm{v}_{hi}),\\
	b_{h}(\bm{u}_{h}^{l},q_{h})=&0,
\end{align}
\end{subequations}
for $\forall (\bm{v}_{h}, q_h )\in  \bm{V}_{h}^{0}\times Q_{h}^{0}$.

It is not difficult to know that the linear system \eqref{Oseen} is uni-solvent for given $(\bm{u}_{h}^{l-1},p_{h}^{l-1})$ and that it holds
\begin{eqnarray}\label{bound-ul}
|||\bm{u}_h^l||| _{V}
\leq   \frac1 \nu \|\bm{f}\|_{*,h},\quad l=1,2,....
\end{eqnarray}

We have the following  convergence result.
\begin{theorem}\label{th73}
Assume that $(\bm{u}_{h},p_{h})\in \bm{V}_{h}^{0}\times Q_{h}^{0}$ is the solution of  the WG scheme (\ref{WG}). Under the  condition
\begin{align}\label{COND}
 2C_{\widetilde{r}}^{r}C_{r}\alpha \frac{\|\bm{f}\|^{r-2}_{\ast,h}}{ \nu^{r-1}}+\frac{\mathcal{N}_{h}\|\bm{f}\|_{\ast,h}}{\nu^{2}}<1
\end{align}
  the Oseen's iteration scheme \eqref{Oseen}  is convergent in the following sense:
\begin{align}\label{conv-p}
\lim_{l\rightarrow\infty}|||\bm{u}_{h}^{l}-\bm{u}_{h} |||_{V}=0,
\lim_{l\rightarrow\infty}|||p_{h}^{l}-p_{h} |||_{Q}=0.
\end{align}
\end{theorem}
\begin{proof}
Denote 	$e_{u}^{l}:=\bm{u}_{h}^{l}-\bm{u}_{h}$ and $e_{p}^{l}:=p_{h}^{l}-p_{h}$. 
Subtracting \eqref{WG} from \eqref{Oseen} gives
\begin{subequations}\label{7.12}
\begin{align}
a_{h}(e_{u}^{l},\bm{v}_{h})
+b_{h}(\bm{v}_{h},e_{p}^{l})-
c_{h}(\bm{u}_{h} ;\bm{u}_{h},\bm{v}_{h})+c_{h}(\bm{u}_{h}^{l-1} ;\bm{u}_{h}^{l},\bm{v}_{h})&\label{7.12a}\nonumber\\
-d_{h}(\bm{u}_{h} ;\bm{u}_{h},\bm{v}_{h})
+d_{h}(\bm{u}_{h}^{l-1} ;\bm{u}_{h}^{l},\bm{v}_{h} )&=0,\\
b_{h}(e_{u}^{l},q_{h})&=0,\label{7.12b}&
\end{align}
\end{subequations}
for any $(\bm{v}_{h}, q_h )\in  \bm{V}_{h}^{0}\times Q_{h}^{0}$.
Taking $\bm{v}_{h}=e_{u}^{l}$, $q_{h}=e_{p}^{l}$ in \eqref{7.12} and using the definition of $c_{h}(\cdot;\cdot,\cdot)$,   Lemmas \ref{Lemma 2.4}, \ref{Lemma 2.9} and \ref{Lemma 3.1}, and the estimates \eqref{sta} and \eqref{bound-ul}, we have
\begin{align}\label{7.13}
 \nu|||e_{u}^{l}|||^{2}_{V}
 =&c_{h}(\bm{u}_{h} ;\bm{u}_{h},e_{u}^{l})-c_{h}(\bm{u}_{h}^{l-1} ;\bm{u}_{h}^{l},e_{u}^{l} )
+ d_{h}(\bm{u}_{h} ;\bm{u}_{h},e_{u}^{l})-d_{h}(\bm{u}_{h}^{l-1} ;\bm{u}_{h}^{l},e_{u}^{l} )\nonumber\\
 =&c_{h}(\bm{u}_{h} ;\bm{u}_{h},e_{u}^{l})
 -c_{h}(\bm{u}_{h}^{l-1} ;\bm{u}_{h},e_{u}^{l} )-c_{h}(\bm{u}_{h}^{l-1};e_{u}^{l},e_{u}^{l})
 -d_{h}(e_{u}^{l-1} ;\bm{u}_{h},e_{u}^{l})\nonumber\\
\leq&c_{h}(\bm{u}_{h} ;\bm{u}_{h},e_{u}^{l})
 -c_{h}(\bm{u}_{h}^{l-1} ;\bm{u}_{h},e_{u}^{l} )
 -d_{h}(e_{u}^{l-1} ;\bm{u}_{h},e_{u}^{l})\nonumber\\
  \leq &(2C_{\widetilde{r}}^{r}C_{r}\alpha \frac{\|\bm{f}\|^{r-2}_{\ast,h}}{ \nu^{r-2}}+\frac{\mathcal{N}_{h}\|\bm{f}\|_{\ast,h}}{\nu}) |||e_{u}^{l-1}|||_{V}  ||| e_{u}^{l}|||_{V},
\end{align}
 which implies
 \begin{align}\label{7.14}
 |||e_{u}^{l}|||_{V}
   \leq &  \mathcal{M} |||e_{u}^{l-1}|||_{V},
\end{align}
with $\mathcal{M}:= 2C_{\widetilde{r}}^{r}C_{r}\alpha \frac{\|\bm{f}\|^{r-2}_{\ast,h}}{ \nu^{r-1}}+\frac{\mathcal{N}_{h}\|\bm{f}\|_{\ast,h}}{\nu^{2}}$. This further means that
\begin{align}
|||e_{u}^{l}|||_{V}
  \leq  & \mathcal{M}|||e_{u}^{l-1}|||_{V} \leq ...
  \leq  \mathcal{M}^{l}|||e_{u}^{0}|||_{V}.
\end{align}
In view of \eqref{COND}, we know that $0< \mathcal{M}<1$.  Thus, we obtain
\begin{align}\label{7.16}
\lim_{l\rightarrow\infty} |||e_{u}^{l}|||_{V}=\lim_{l\rightarrow\infty} |||\bm{u}_{h}^{l}-\bm{u}_{h} |||_{V}=0.
\end{align}

The thing left is to prove the second convergence relation of \eqref{conv-p}.
From   \eqref{7.12a}  it follows
\begin{align*}
b_{h}(\bm{v}_{h},e^{l}_{p})=&-a_{h}(e_{u}^{l},\bm{v}_{h})
+c_{h}(\bm{u}_{h} ;\bm{u}_{h},\bm{v}_{h})-c_{h}(\bm{u}_{h}^{l-1} ;\bm{u}_{h}^{l},\bm{v}_{h} )\nonumber\\
&+d_{h}(\bm{u}_{h} ;\bm{u}_{h},\bm{v}_{h})
-d_{h}(\bm{u}_{h}^{l-1} ;\bm{u}_{h}^{l},\bm{v}_{h})\nonumber\\
=&-a_{h}(e_{u}^{l},\bm{v}_{h})+\big(c_{h}(\bm{u}_{h} ;\bm{u}_{h},\bm{v}_{h})
 -c_{h}(\bm{u}_{h}^{l-1} ;\bm{u}_{h},\bm{v}_{h} )\big)-c_{h}(\bm{u}_{h}^{l-1};e_{u}^{l},\bm{v}_{h})\nonumber\\
&-\big(d_h(\bm{u}_{h} ;e_{u}^{l},\bm{v}_{h} )+d_h(e_{u}^{l-1} ;e_{u}^{l},\bm{v}_{h} )+d_h(e_{u}^{l-1};\bm{u}_{h},\bm{v}_{h}) \big),
\end{align*}
for all $\bm{v}_{h}\in \bm{V}_{h}^{0} $. By  Lemma  \ref{theoremLBB} we have
$$
 |||e^{l}_{p}|||_{Q}
  \lesssim\sup_{ \bm{v}_{h}\in \bm{V}_{h}^{0}}\frac{b_{h}(\bm{v}_{h},e^{l}_{p})}{|||\bm{v}_{h}|||_{V} }.$$
  The above two results,
   together with \eqref{7.16} and  Lemmas  \ref{Lemma 2.4}, \ref{Lemma 2.9} and \ref{Lemma 3.1}, yield the desired conclusion
\begin{align*}
\lim_{l\rightarrow\infty} |||e^{l}_{p}|||_{Q}
=\lim_{l\rightarrow\infty} |||p_{h}^{l}-p_{h}|||_{Q}=0.
\end{align*}
This completes this proof.

\end{proof}

\section{Numerical experiments}

In this section, we provide some numerical tests to verify the performance of the WG scheme  \eqref{WG} for the Brinkman-Forchheimer model \eqref{BF0} in two dimensions. We adopt the Oseen's iterative algorithm \eqref{Oseen} with the initial guess $\bm{u}_{hi}^{0}=\bm{0}$ and the stop criterion
\begin{equation}
\|\bm{u}_{h}^{l}-\bm{u}_{h}^{l-1}\|_{0}<1e-8.
\end{equation}
in all the  numerical examples, i.e. Examples \ref{EX7.1} to \ref{EX7.4}.

\begin{figure}[!t]
		  \centering
	 \includegraphics[width= 2.7in]{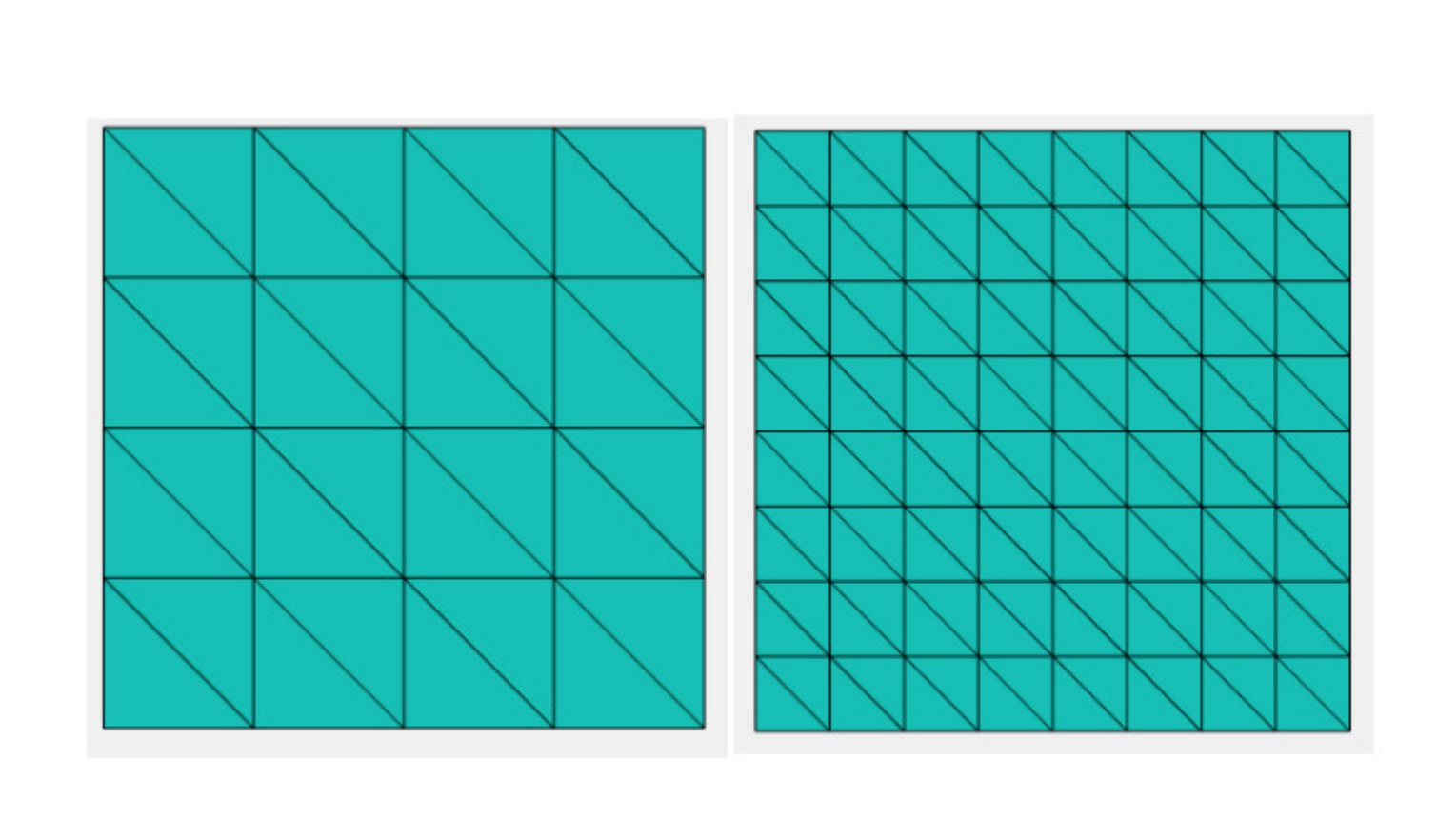}\\
		   \caption{Uniform triangular meshes:    $4\times4$ mesh (left)      and $8\times 8 $   mesh (right).}
		\label{fig1:mesh}
		\end{figure}

\begin{exam} \label{EX7.1}
 Set $\Omega=  [0,1]\times[0,1]$, $\nu=1$,  $\alpha=5$ and $r=10$ in the model \eqref{BF0}.
 The exact solution $(\bm{u},p)$ is given as follows:
\begin{equation}\label{EX2}
\left\{\begin{array}{ll}
u_{1}=10x^{2}(x-1)^{2}y(y-1)(2y-1),        \\
u_{2}=-10x(x-1)(2x-1)y^{2}(y-1)^{2} ,    \\
p=10(2x-1)^{2}(2y-1).
\end{array}\right.
\end{equation}
We compute the scheme \eqref{WG} on uniform triangular meshes  (cf. Figure \ref{fig1:mesh}), with  $m=1,2$, $k=m-1,m$. Numerical results of  $\|\bm{u}-\bm{u}_{hi}\|_0$, $\|\nabla \bm{u}-\nabla_{h}\bm{u}_{hi}\|_0$, $\|p-p_{hi}\|_0$ and $\| \nabla\cdot\bm{u}_{hi} \|_{0,\infty}$.
are listed in Tables \ref{tab1} and \ref{tab2}.

\end{exam}

From the numerical results we have the following observations: 
\begin{itemize}
  \item The convergence rates of $\|\nabla \bm{u}-\nabla_{h}\bm{u}_{hi}\|_0$ and $\|p-p_{hi}\|_0$  for the WG scheme are  $m^{th}$ orders  in the cases of $m=1,2$ and $k=m, m-1$. These are  conformable to the theoretical results in Theorem \ref{Theorem 4.1}.

     \item The convergence rate of $\|\bm{u}-\bm{u}_{hi}\|_0$  is   $(m+1)^{th}$ order,   which is  conformable to the theoretical result in Theorem \ref{Theorem 5.1}.

   \item The results of $\|\nabla_{h}\cdot\bm{u}_{hi}\|_{0,\infty} $ 
   are almost zero. This means that the discrete velocity is globally divergence-free, which is consistent with  Theorem \ref{TH2.2}.
\end{itemize}

\begin{table}[h]
	\small
	\caption{\label{tab1}
		History of convergence results for Example \ref{EX7.1}:  $m=1$
		  }
	\centering
	
	\begin{tabular}{c|c|c|c|c|c|c|c|c}
		\Xhline{1pt}
		
		\multirow{2}{*}{$k$} &
		\multirow{2}{*}{$mesh$} &
		
		\multicolumn{2}{c|}{$\frac{\|\bm{u}-\bm{u}_{hi}\|_0}{\|\bm{u}\|_0}$} &
		\multicolumn{2}{c|}{$\frac{\|\nabla \bm{u}-\nabla_{h}\bm{u}_{hi}\|_0}{\|\nabla \bm{u}\|_0}$} &
		
		\multicolumn{2}{c| }{$\frac{\|p-p_{hi}\|_0}{\|p\|_0}$} &
		\multicolumn{1}{c }{
			\footnotesize  $\|\nabla_{h}\cdot\bm{u}_{hi}\|_{0,\infty}$
		} \\
		\cline{3-9}
		
		& &Error &Rate  &Error &Rate  &Error &Rate &Error \\
		\hline
		
		\multirow{6}{*}{$0$}

&$4\times4$     &5.9583e-01   &   -   &5.1516e-01   &   -   &2.8667e-01   &   -   &4.0593e-16\\
&$8\times8$     &1.5876e-01   &1.91   &2.7301e-01   &0.92   &1.4424e-01   &0.99   &2.3028e-16\\
&$16\times16$   &4.1525e-02   &1.93   &1.3851e-01   &0.98   &7.2201e-02   &1.00   &3.1127e-16\\
&$32\times32$   &1.0641e-02   &1.96   &6.9420e-02   &1.00   &3.6100e-02   &1.00   &8.4459e-17\\
&$64\times64$   &2.6985e-03   &1.98   &3.4723e-02   &1.00   &1.8048e-02   &1.00   &2.7905e-17\\
&$128\times128$ &6.9479e-04   &1.96   &1.7364e-02   &1.00   &9.0235e-03   &1.00   &5.6257e-17\\

		\Xhline{1pt}
		\multirow{6}{*}{$1$}
&$4\times4$     &5.6714e-01   &   -   &5.1165e-01   &   -   &2.8667e-01   &   -   &3.4694e-18\\
&$8\times8$     &1.5224e-01   &1.90   &2.7237e-01   &0.91   &1.4425e-01   &0.99   &3.2092e-17\\
&$16\times16$   &3.9918e-02   &1.93   &1.3841e-01   &0.98   &7.2213e-02   &1.00   &1.3010e-18\\
&$32\times32$   &1.0236e-02   &1.96   &6.9404e-02   &1.00   &3.6110e-02   &1.00   &2.9328e-17\\
&$64\times64$   &2.5908e-03   &1.98   &3.4720e-02   &1.00   &1.8054e-02   &1.00   &8.2115e-17\\
&$128\times128$ &6.5160e-04   &1.99   &1.7363e-02   &1.00   &9.0267e-03   &1.00   &1.4732e-17\\
			
		\Xhline{1pt}
	\end{tabular}
	
\end{table}

\begin{table}[h]
	\small
	\caption{\label{tab2}
		History of convergence results for Example \ref{EX7.1}:  $m=2$
		  }
	\centering
	
	\begin{tabular}{c|c|c|c|c|c|c|c|c}
		\Xhline{1pt}
		
		\multirow{2}{*}{$k$} &
		\multirow{2}{*}{$mesh$} &
		
		\multicolumn{2}{c|}{$\frac{\|\bm{u}-\bm{u}_{hi}\|_0}{\|\bm{u}\|_0}$} &
		\multicolumn{2}{c|}{$\frac{\|\nabla \bm{u}-\nabla_{h}\bm{u}_{hi}\|_0}{\|\nabla \bm{u}\|_0}$}&
		
		\multicolumn{2}{c| }{$\frac{\|p-p_{hi}\|_0}{\|p\|_0}$} &
		\multicolumn{1}{c }{
			\footnotesize  $\|\nabla_{h}\cdot\bm{u}_{hi}\|_{0,\infty}$
		} \\
		\cline{3-9}
		
		& &Error &Rate  &Error &Rate  &Error &Rate &Error \\
		\hline
		
		\multirow{6}{*}{$1$}
&$4\times4$     &5.9915e-02   &   -   &1.3040e-01   &   -   &3.3330e-02   &   -   &2.2560e-14\\
&$8\times8$     &7.6055e-03   &2.98   &3.4961e-02   &1.90   &8.3117e-03   &2.00   &7.0083e-16\\
&$16\times16$   &9.4731e-04   &3.01   &8.9617e-03   &1.96   &2.0761e-03   &2.00   &1.9606e-15\\
&$32\times32$   &1.1827e-04   &3.00   &2.2616e-03   &1.99   &5.1888e-04   &2.00   &7.4921e-16\\
&$64\times64$   &1.4789e-05   &3.00   &5.6761e-04   &1.99   &1.2970e-04   &2.00   &2.5543e-17\\
&$128\times128$ &1.8494e-06   &3.00   &1.4215e-04   &2.00   &3.2423e-05   &2.00   &3.5312e-16\\

		\Xhline{1pt}
		\multirow{6}{*}{$2$}
&$4\times4$     &5.6852e-02   &   -   &1.3016e-01   &   -   &3.3281e-02   &   -   &1.1310e-15\\
&$8\times8$     &7.3762e-03   &2.95   &3.4894e-02   &1.90   &8.2971e-03   &2.00   &1.7295e-15\\
&$16\times16$   &9.3175e-04   &2.98   &8.9496e-03   &1.96   &2.0724e-03   &2.00   &2.5093e-15\\
&$32\times32$   &1.1713e-04   &2.99   &2.2591e-03   &1.99   &5.1795e-04   &2.00   &4.0441e-17\\
&$64\times64$   &1.4693e-05   &3.00   &5.6704e-04   &1.99   &1.2947e-04   &2.00   &1.0278e-15\\
&$128\times128$ &1.8403e-06   &3.00   &1.4201e-04   &2.00   &3.2366e-05   &2.00   &1.1529e-16\\

		\Xhline{1pt}
	\end{tabular}
	
\end{table}

\begin{figure}[htbp!]
\centering
\subfigure[ velocity (Taylor-Hood)     ]
{\includegraphics[width=4.5cm]{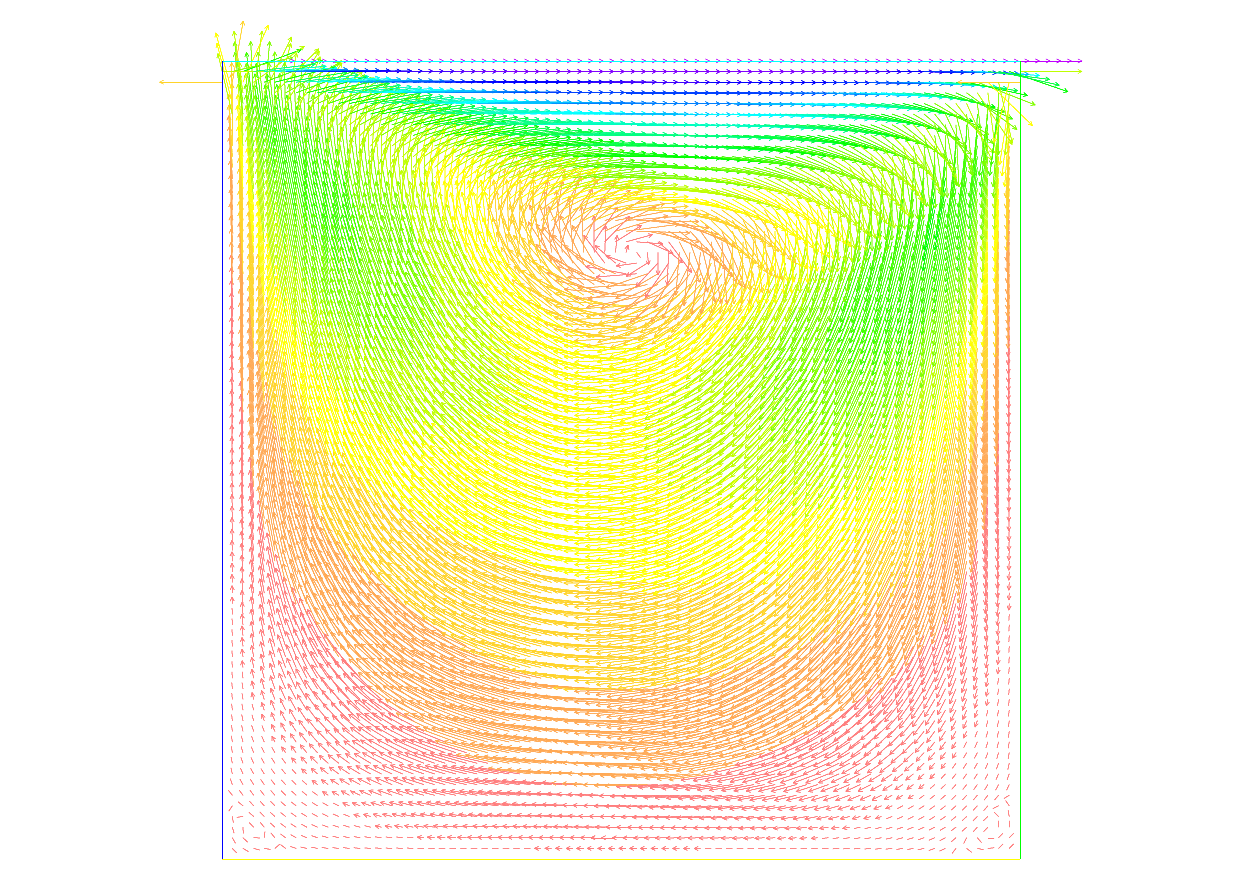}}
\quad
\subfigure[  pressure  (Taylor-Hood) ]
{\includegraphics[width=4.5cm]{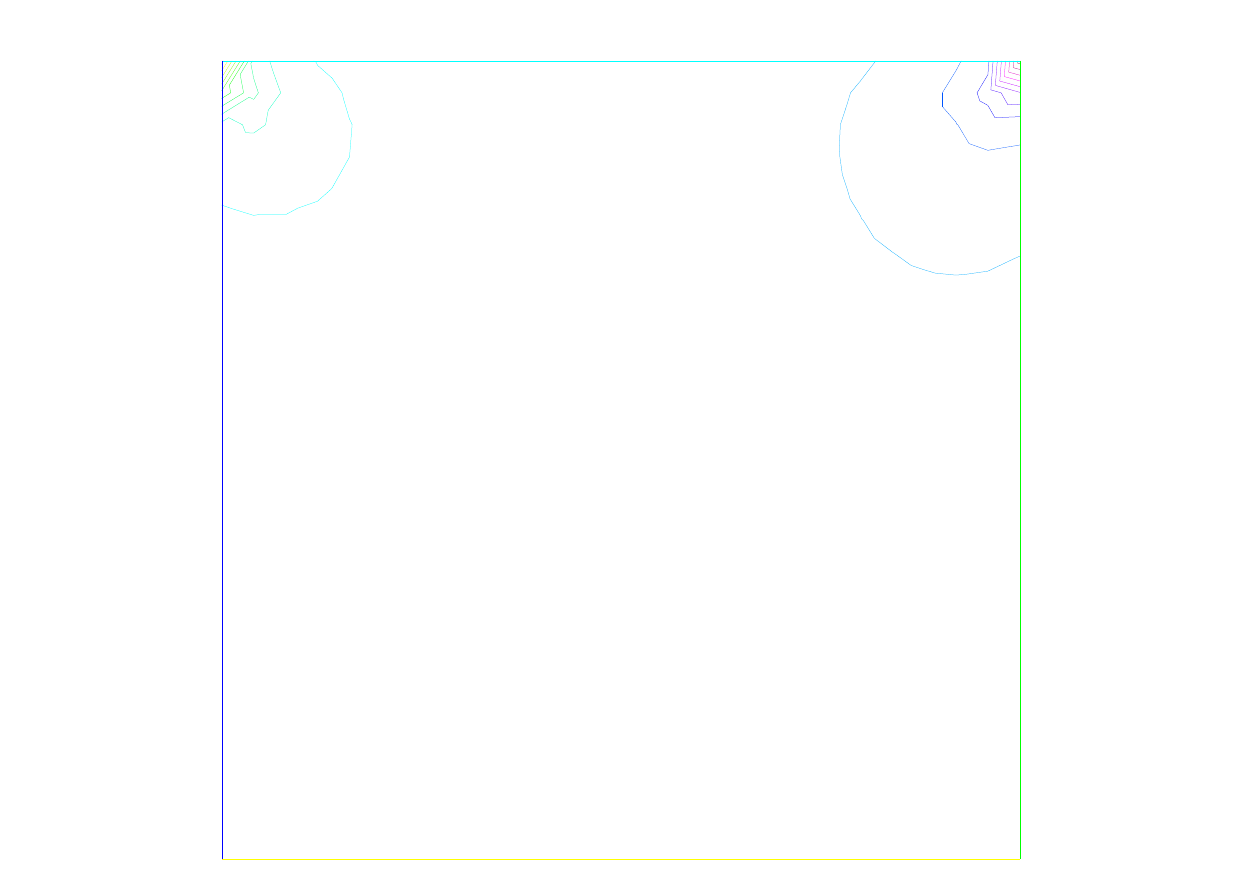}}\\

\subfigure[  velocity (WG) ]
{\includegraphics[width=4.5cm]{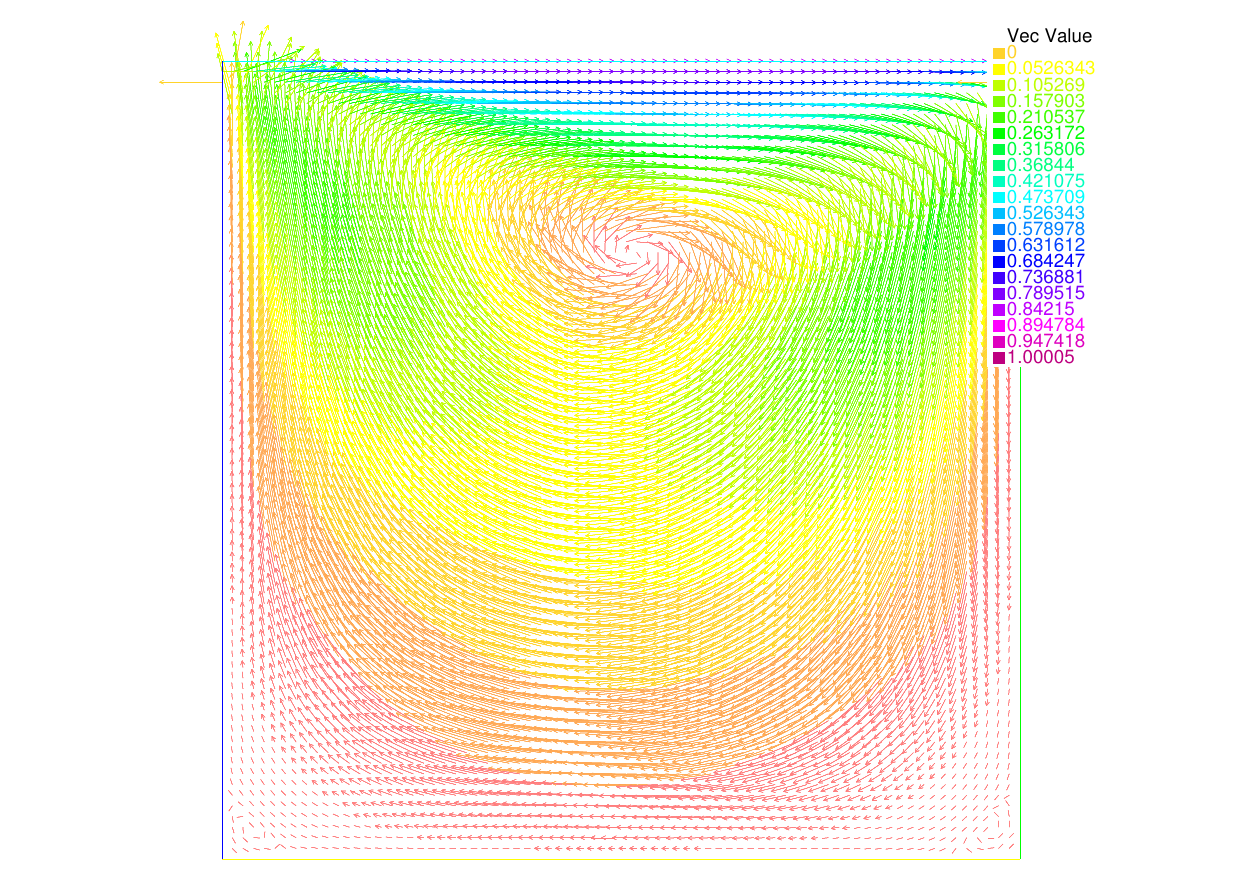}}
\quad
\subfigure[  pressure (WG) ]
{\includegraphics[width=4.5cm]{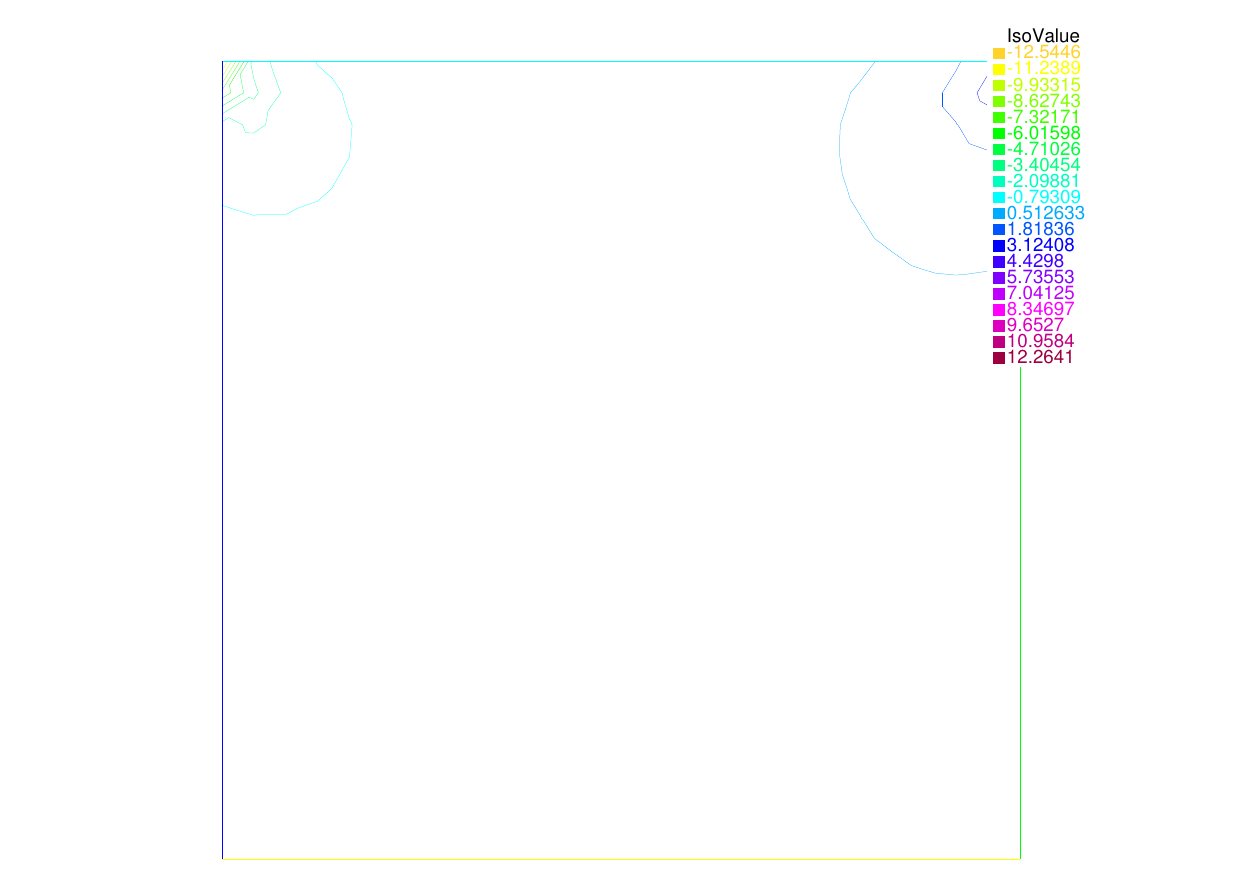}}\\
\caption{ The velocity streamlines  and pressure contours for Example \ref{EX7.2}: $\alpha=0$}
\label{fig21:21}
\end{figure}

\begin{figure}[htbp!]
\centering
\subfigure[ velocity: $\alpha=1$]
{\includegraphics[width=4 cm]{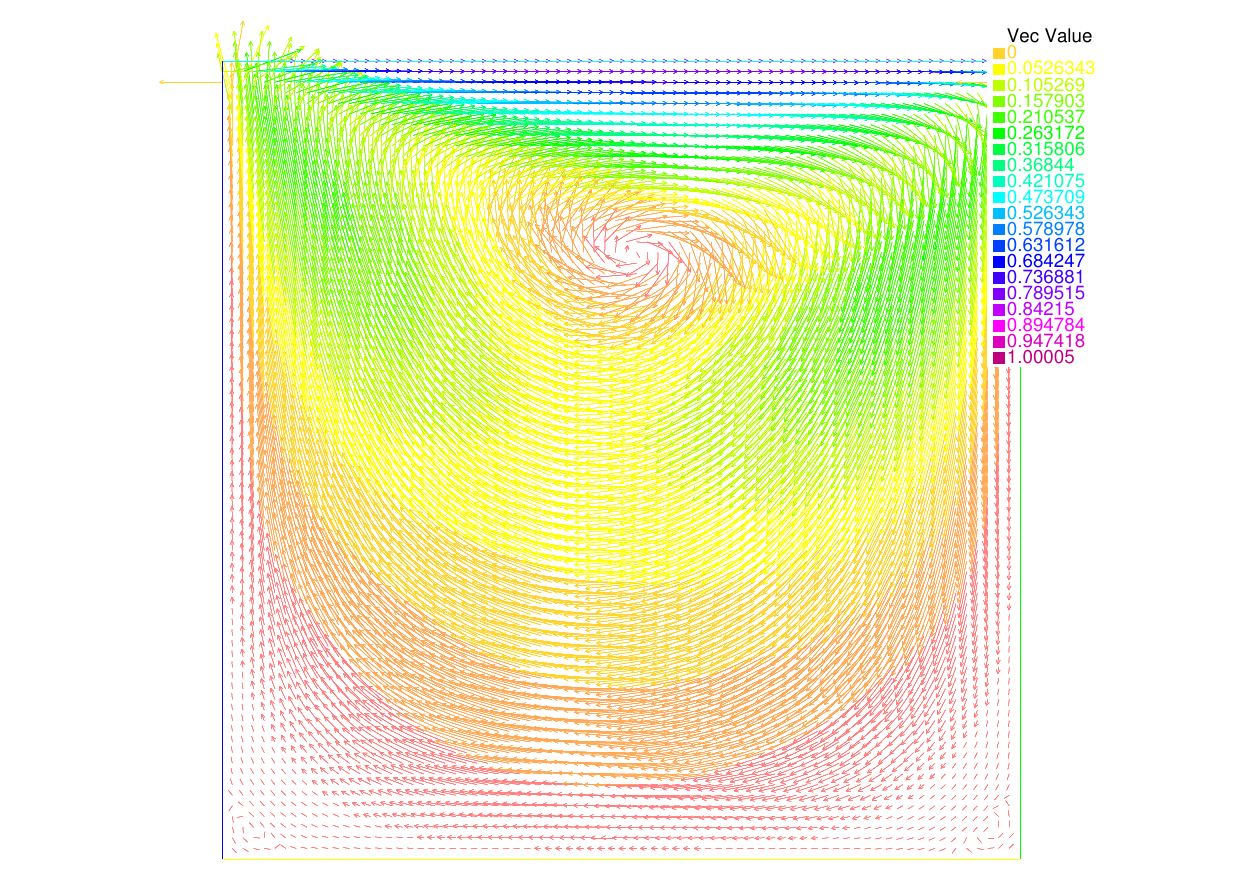}}
\subfigure[ velocity: $\alpha=50 $]
{\includegraphics[width=4 cm]{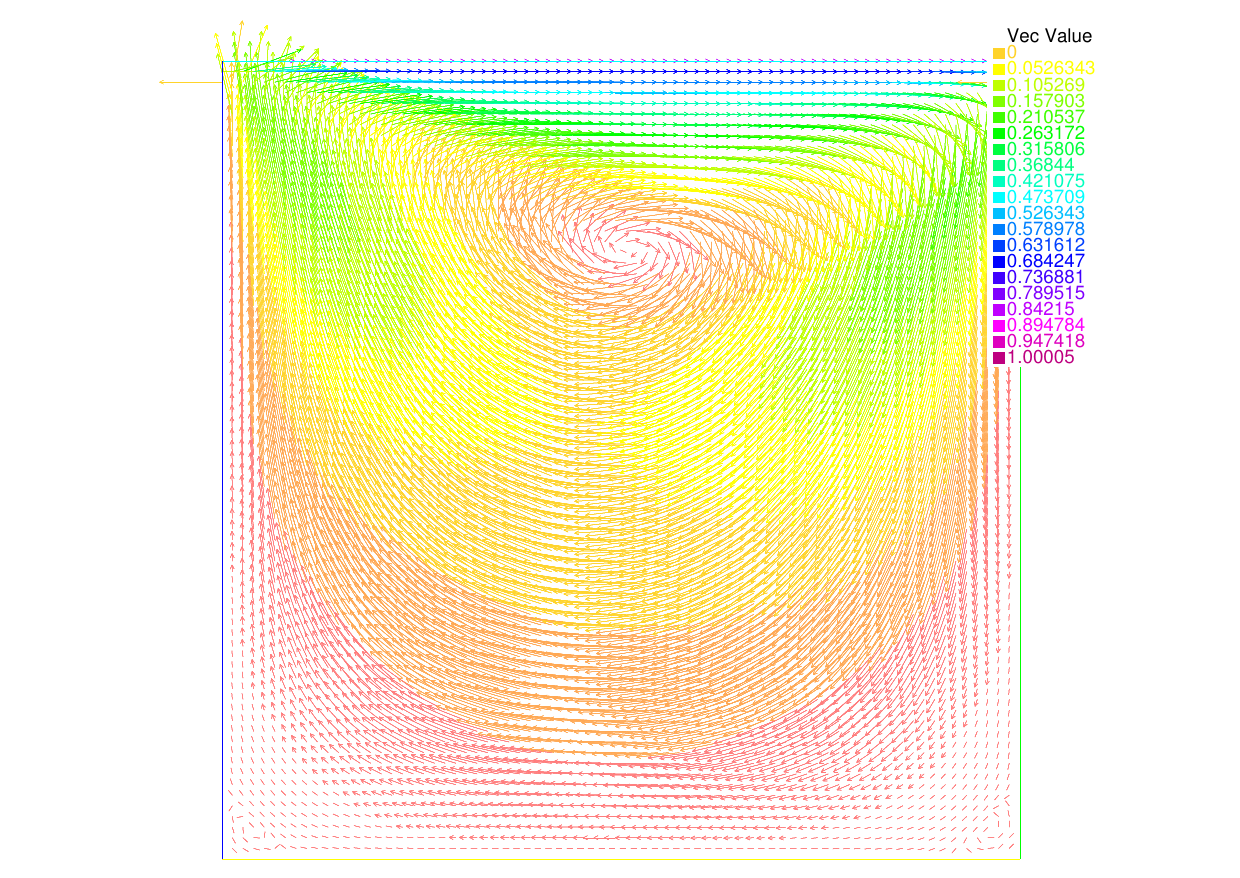}}%
\subfigure[ velocity: $\alpha=100 $]
{\includegraphics[width=4 cm]{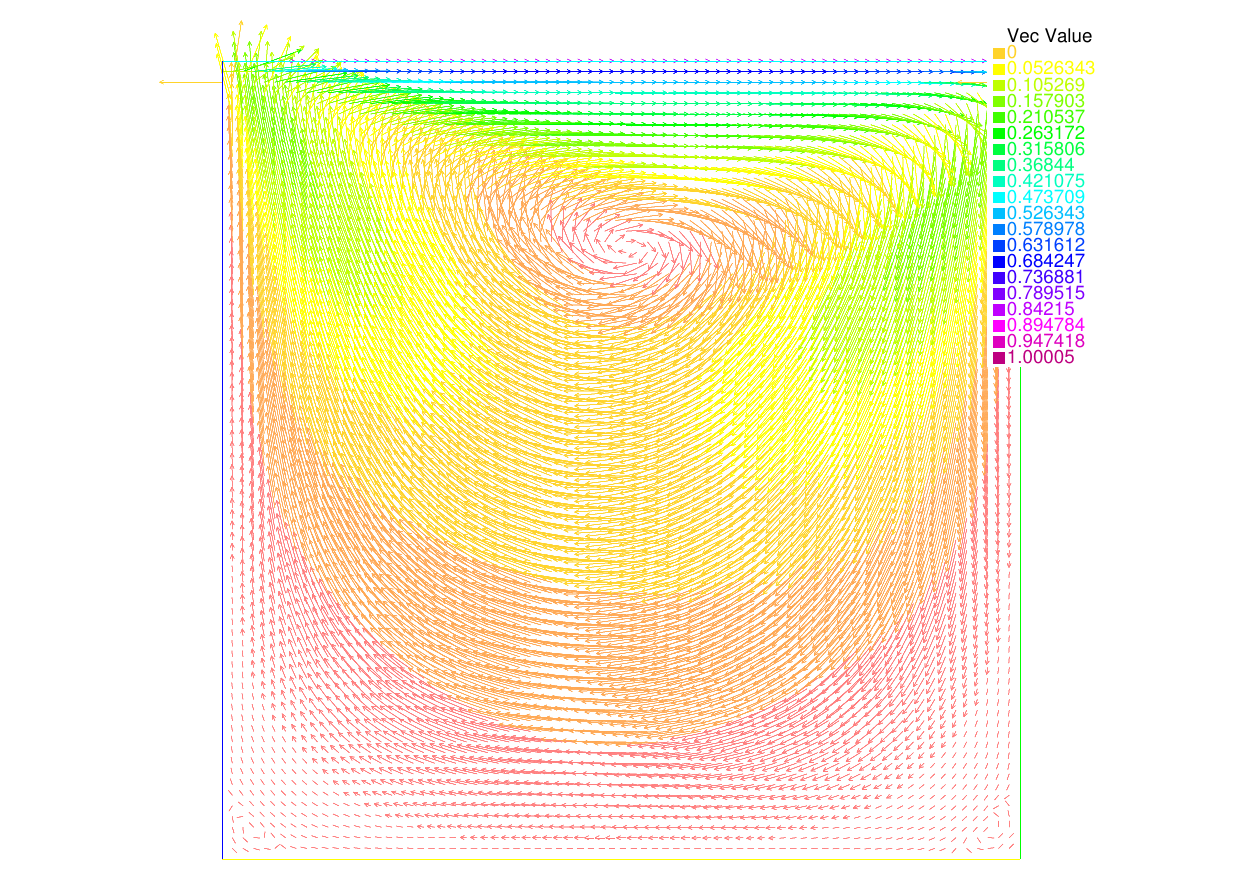}}\\
\subfigure[  pressure: $\alpha=1 $]
{\includegraphics[width=3.9 cm]{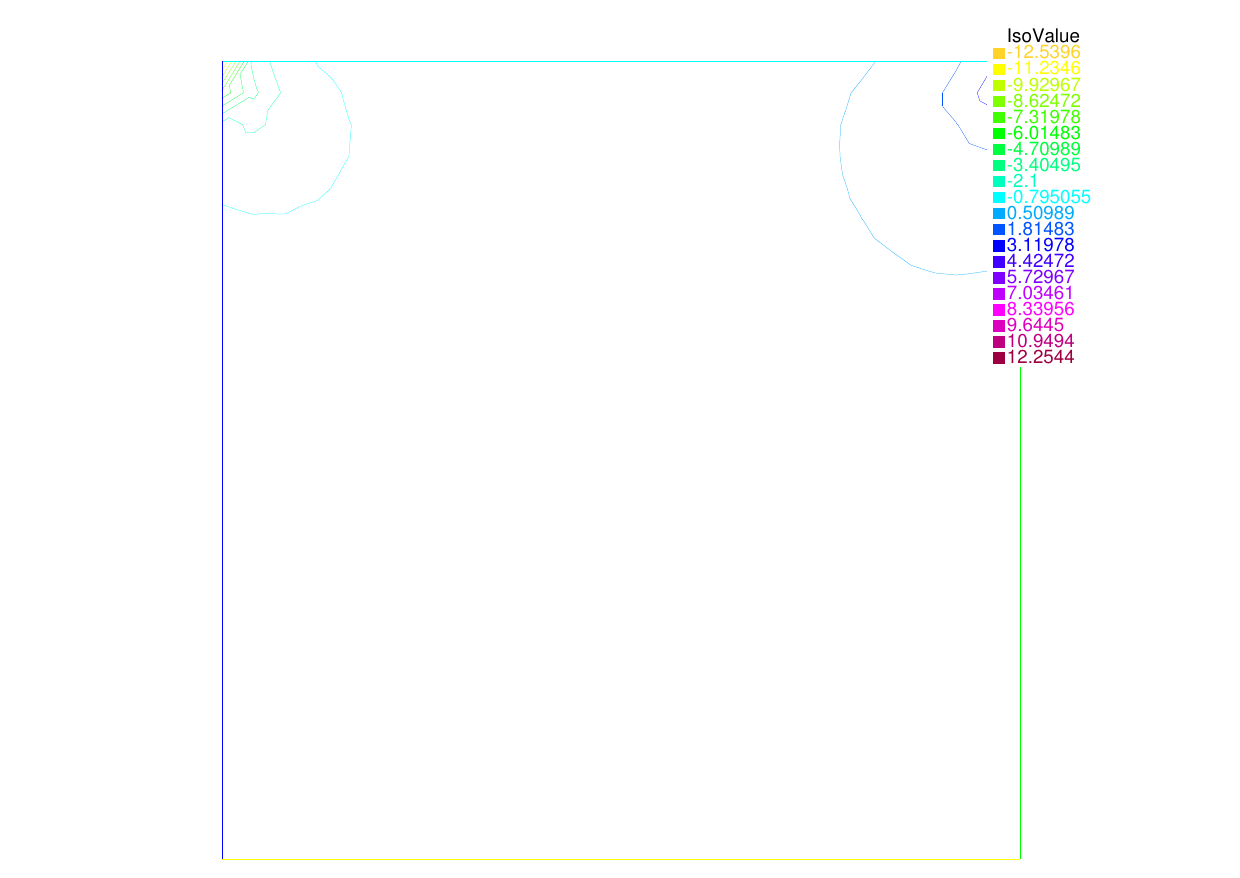}}
\subfigure[  pressure: $\alpha=50 $]
{\includegraphics[width=3.9 cm]{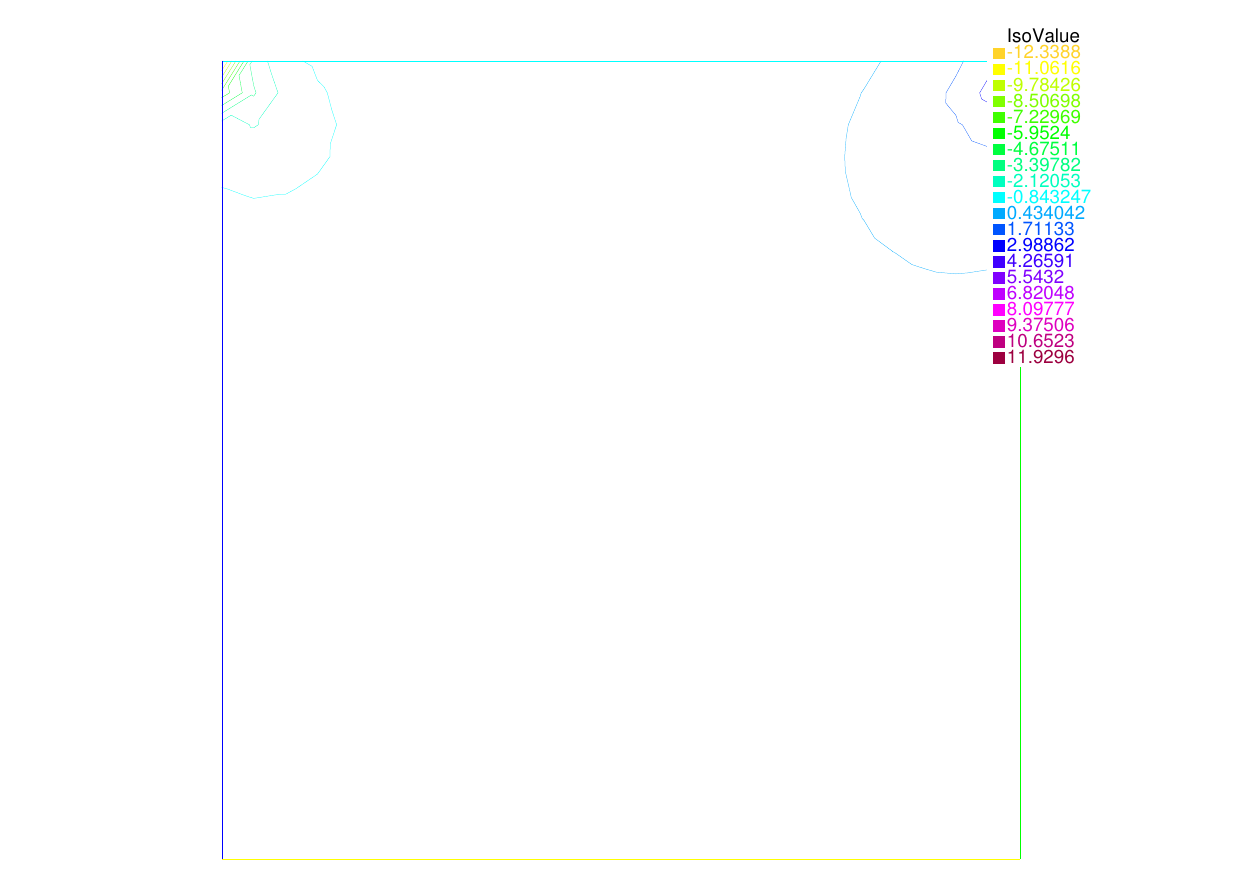}}
\subfigure[  pressure: $\alpha=100 $]
{\includegraphics[width=3.9 cm]{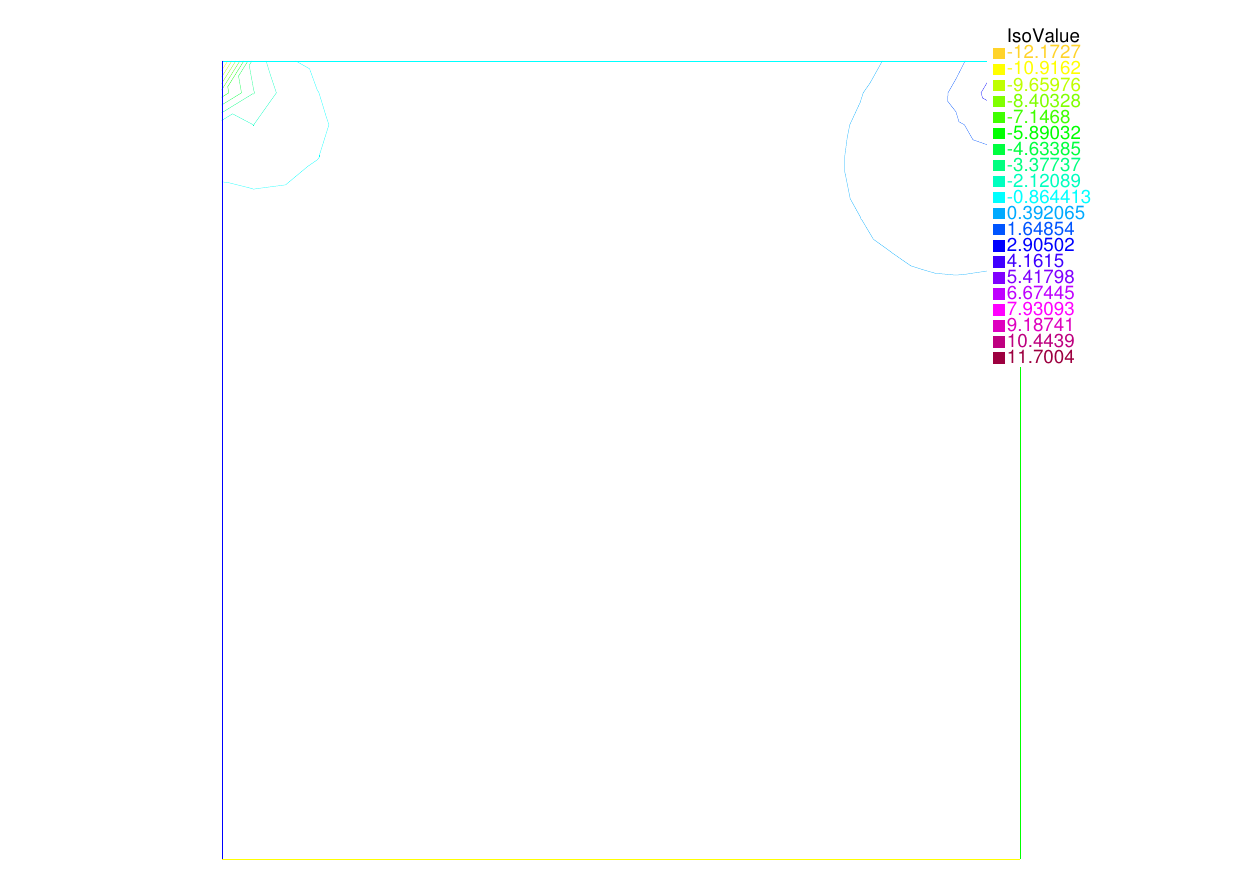}}
\caption{The   velocity streamlines and pressure contours  for Example \ref{EX7.2}:  $r=5$ and   $\alpha=1, 50, 100$}
\label{fig22:22}
\end{figure}

\begin{figure}[htbp!]
\centering
\subfigure[velocity: $r=3$]
{\includegraphics[width=4cm]{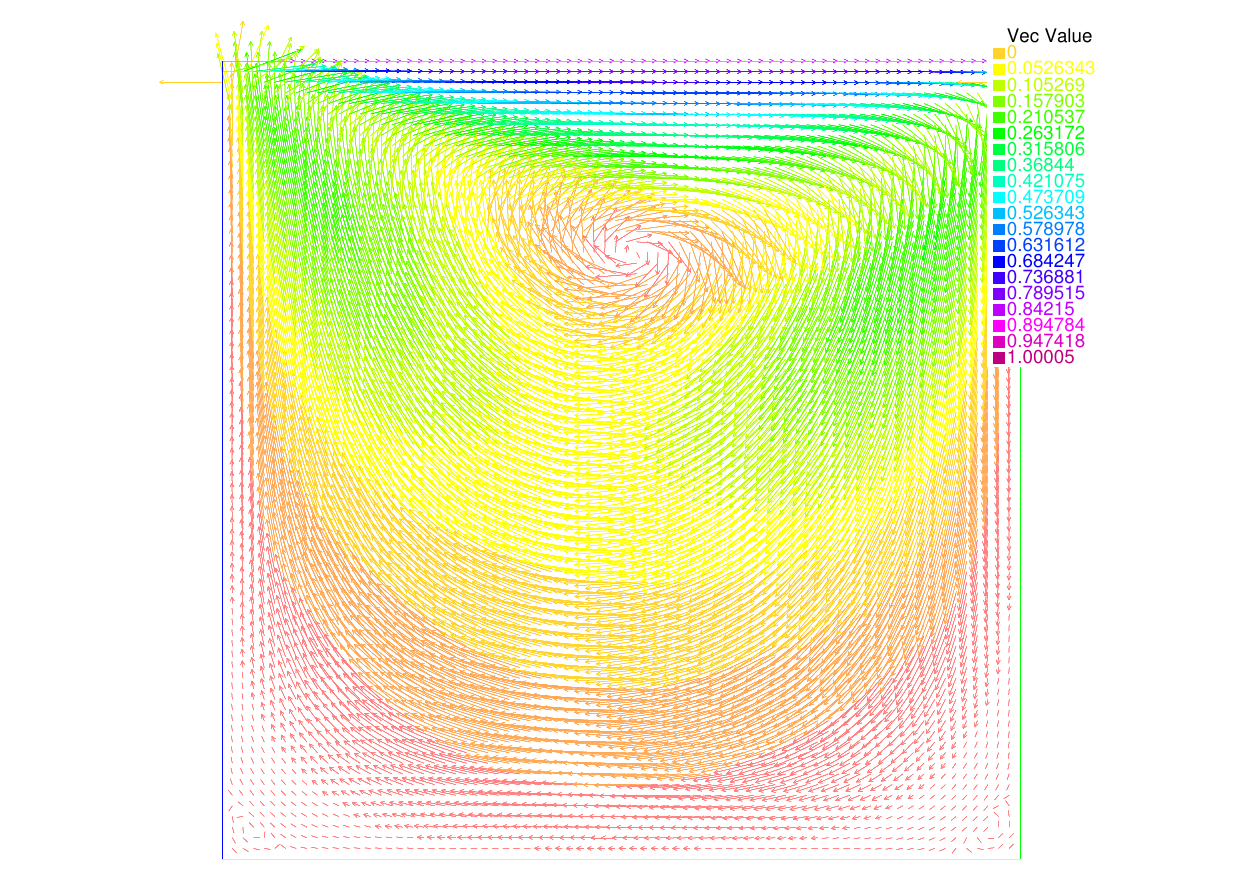}}
\subfigure[velocity: $r=5$]
{\includegraphics[width=4cm]{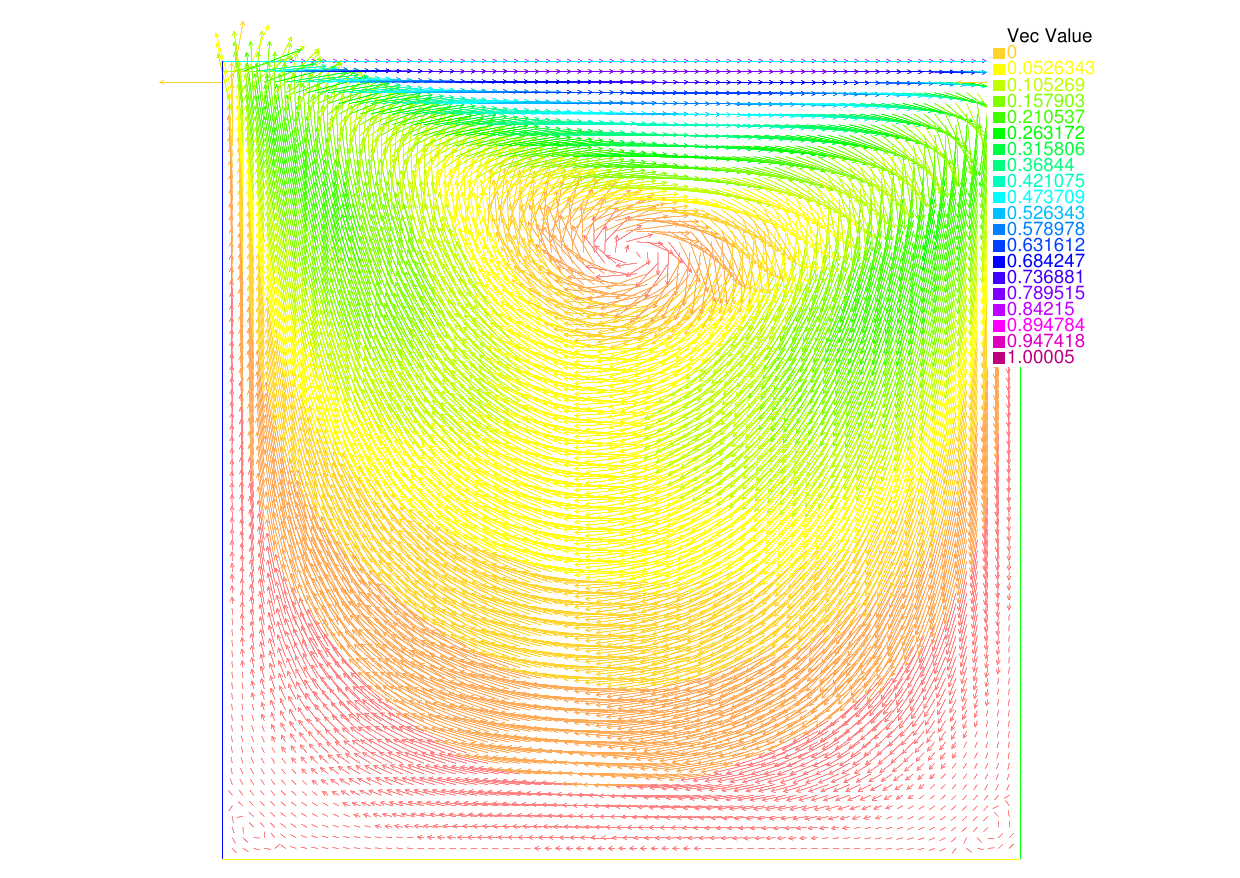}}
\subfigure[velocity: $r=50$]
{\includegraphics[width=4cm]{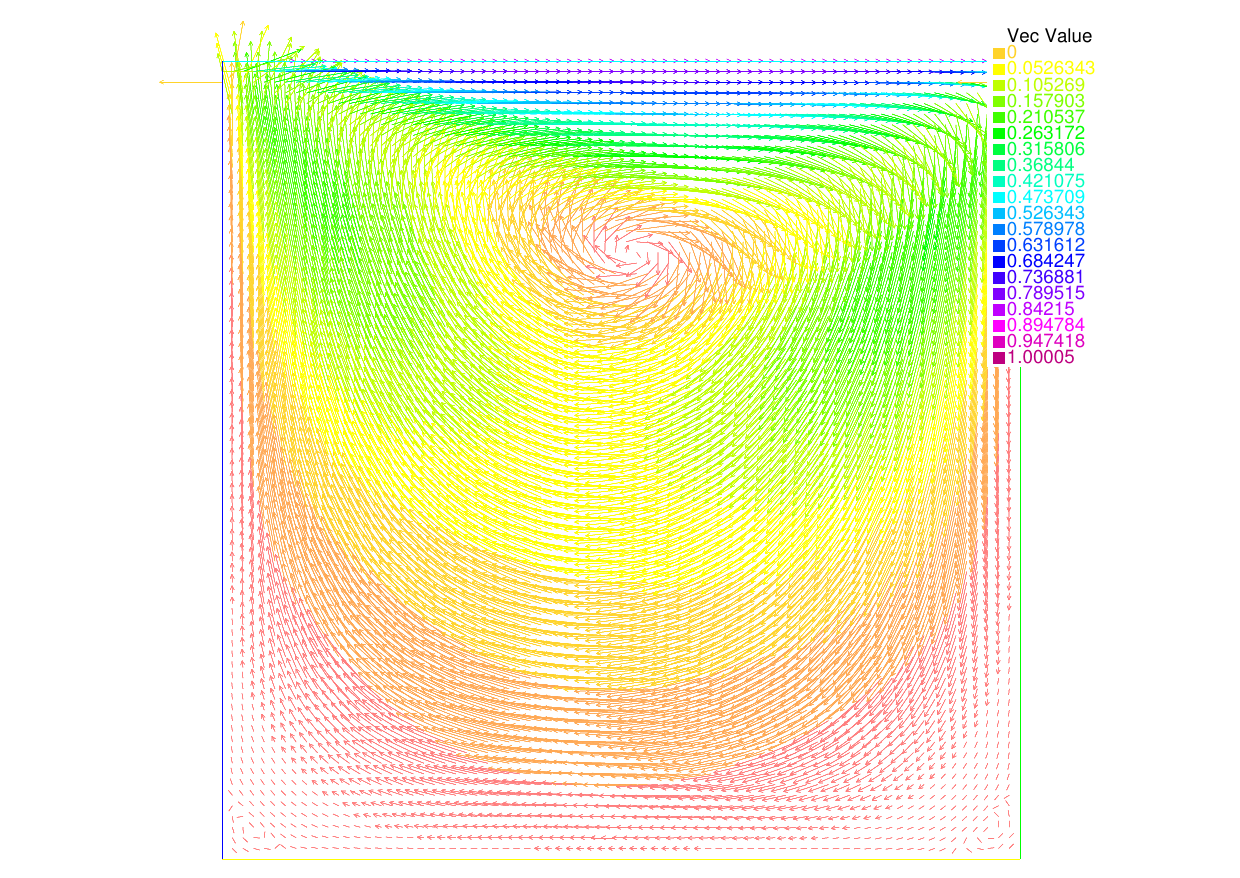}}
\quad
\subfigure[pressure: $r=3$]
{\includegraphics[width=4cm]{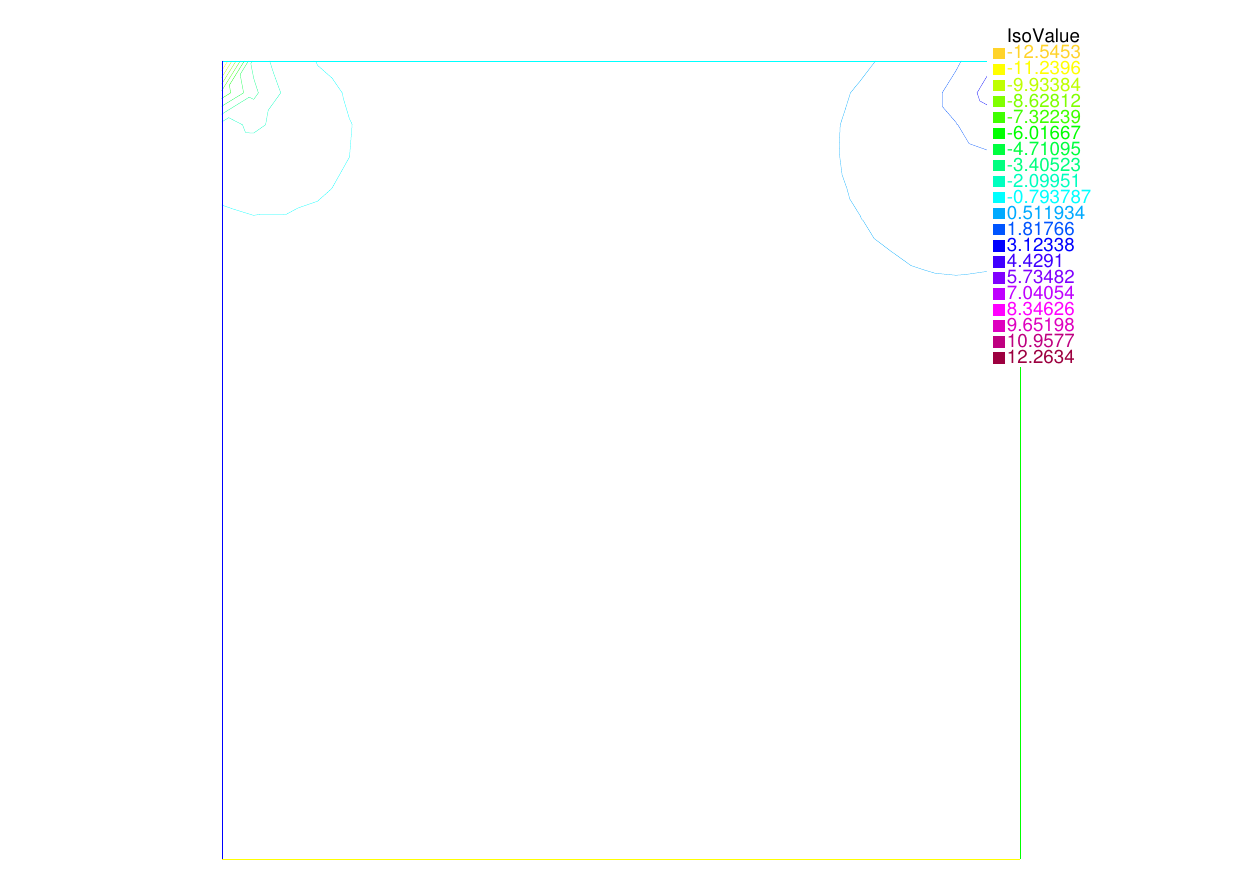}}
\subfigure[ pressure: $r=5$]
{\includegraphics[width=4cm]{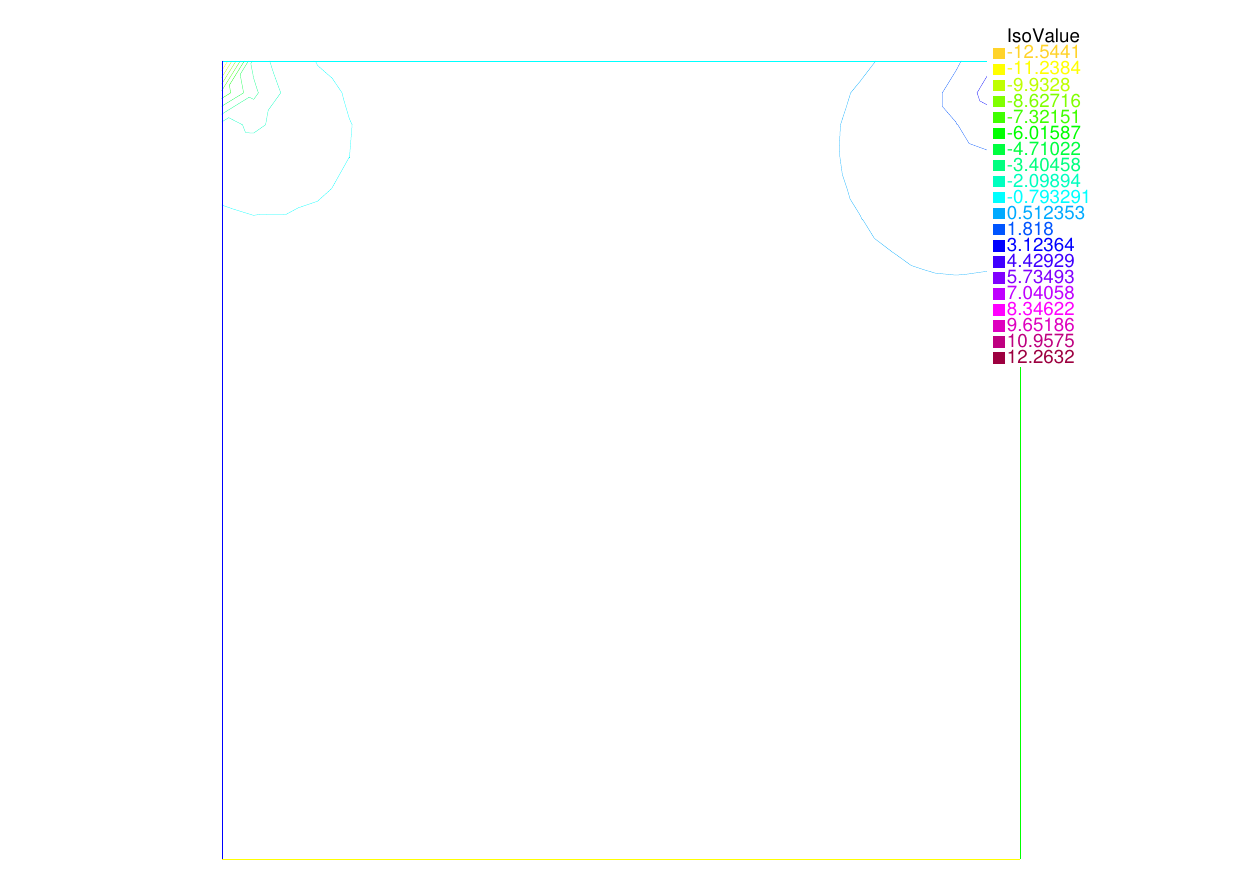}}
\subfigure[ pressure: $r=50$]
{\includegraphics[width=4cm]{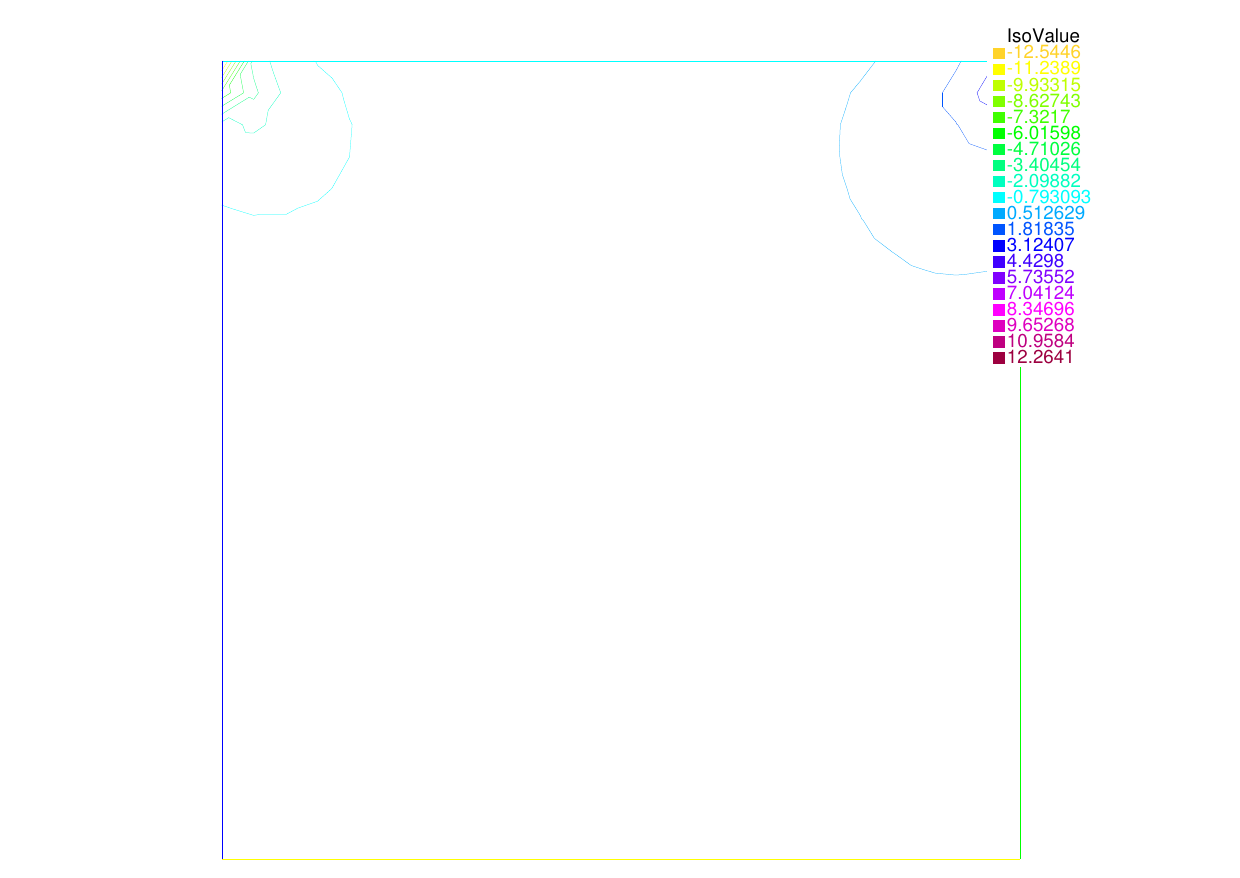}}
\caption{The   velocity streamlines and pressure contours  for Example \ref{EX7.2}:  $\alpha=5$ and    $r=3,5,50$}
\label{fig23:23}
\end{figure}

\begin{exam}[The lid-driven cavity flow problem]\label{EX7.2}
This problem is used to test the influence of damping parameters $\alpha$  and $r$ on the solution of the WG scheme. Take   $\Omega$ =$[0, 1]\times[0, 1]$, $\nu=0.1$  and $\bm{f}=\bm{0}$. The   boundary conditions are as follows:
$$\bm{u}|_{x=0}=\bm{u}|_{x=1}=\bm{u}|_{y=0}=\bm{0}, \quad \bm{u}|_{y=1}=(1,0)^T.$$ 
We compute the WG scheme \eqref{WG} with  $m=k=2$ on the $25\times 25$ uniform triangular mesh  (cf. Figure \ref{fig1:mesh}) in the following cases:
\begin{itemize}
\item [ I]. $\alpha=0$, i.e. the case of the Navier-Stokes  equations;

\item [ II].   $r=5$ and $\alpha=1, 50, 100$;

\item [ III]. $\alpha=5$ and $r=3, 5, 50$.
\end{itemize}
The     velocity streamlines and the   pressure contours are displayed in Figures \ref{fig21:21}, \ref{fig22:22} and \ref{fig23:23}.
  As a comparison,  the  referenced  numerical solutions obtained  with  the Taylor-Hood element  
  are also shown   for  $\alpha=0$; see (a) and (b) in Figure \ref{fig21:21}.

 From  Figure \ref{fig22:22} we can see that the shape and   size  of the vortex  change evidently, which means that  the damping effect becomes greater for the velocity  as the damping parameter  $\alpha$ increases.   We can also  see that the  pressure approximation is not significantly affected by $\alpha$.
On the other hand, as shown in
Figure  \ref{fig23:23},  the velocity  and pressure approximations are not  significantly effected by the number  $r$.
 \end{exam}

\begin{figure}[htbp!]
\centering
\subfigure
{\includegraphics[width=12cm]{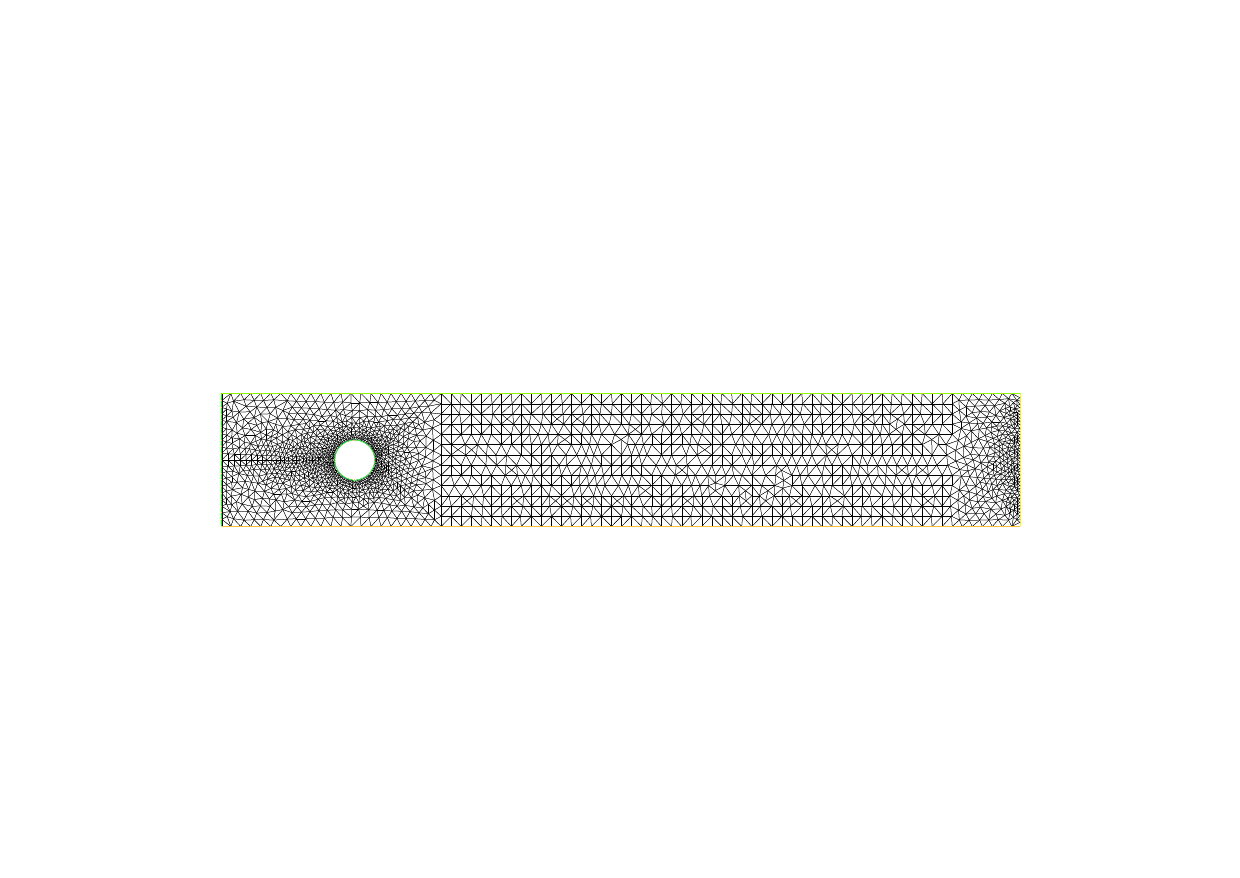}}
\caption{The domain and  finite element mesh for Example \ref{EX7.3}}
\label{fig2:mesh2}
\end{figure}

\begin{figure}[htbp!]
\centering
\setlength{\abovecaptionskip}{0.cm}
\subfigure[velocity (Taylor-Hood)]
{\includegraphics[width=4cm]{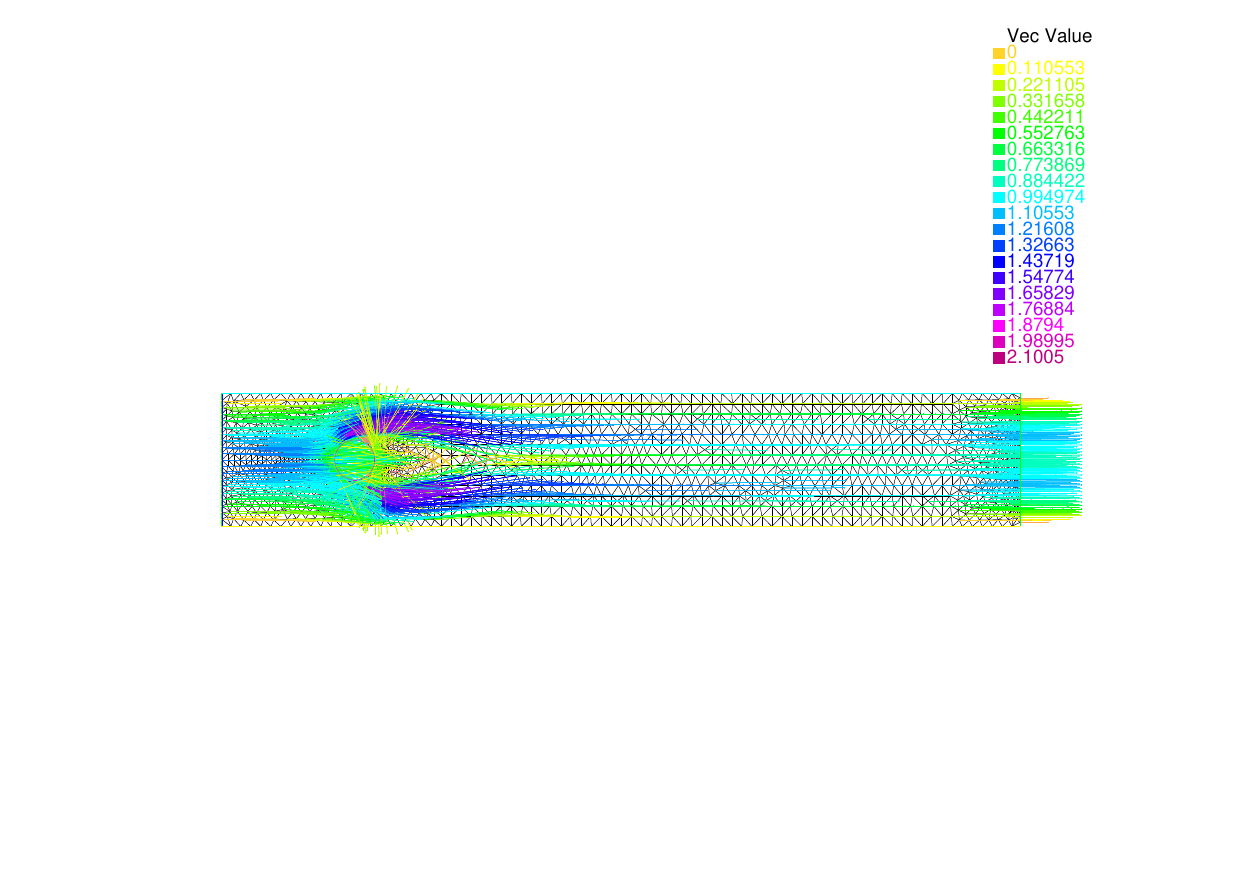}}
\subfigure[vorticity (Taylor-Hood) ]
{\includegraphics[width=4cm]{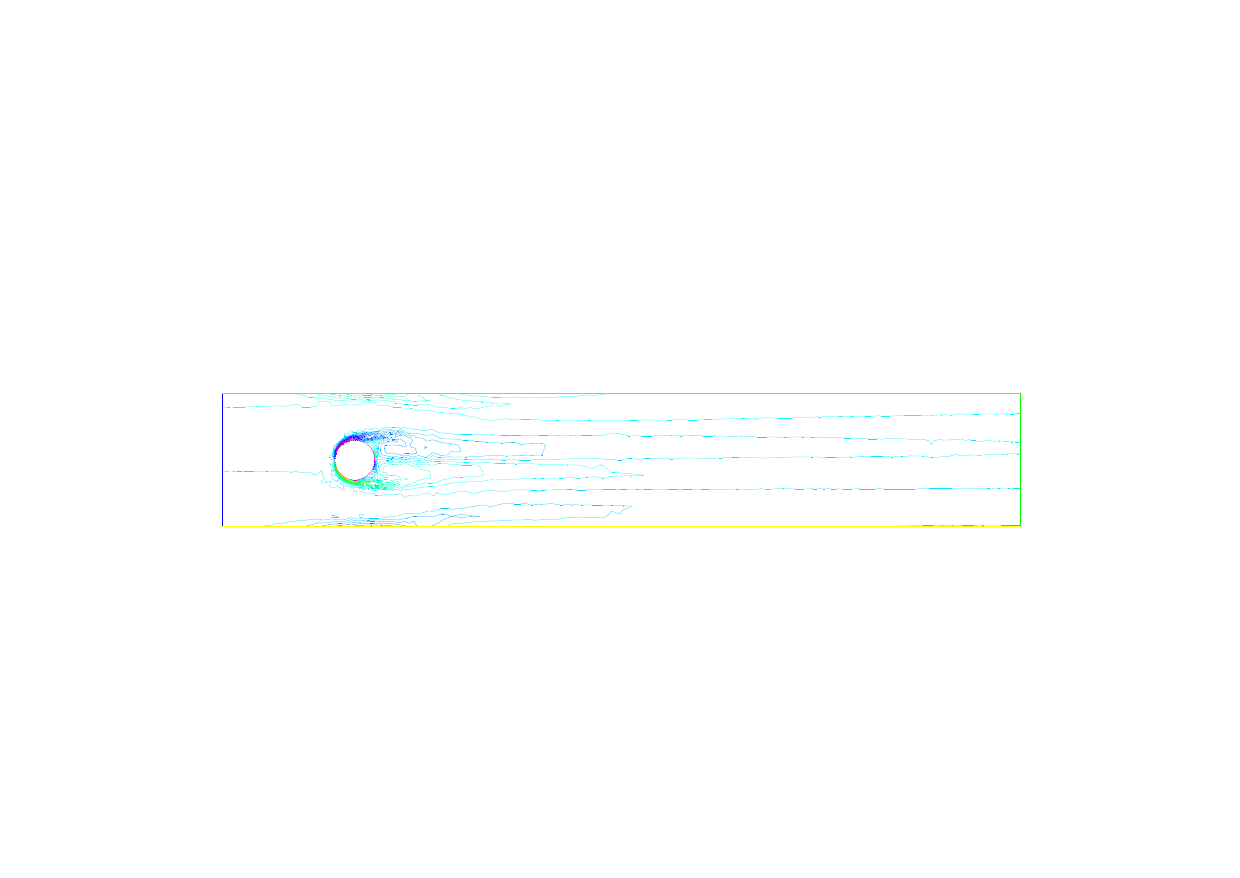}}
\subfigure[pressure (Taylor-Hood) ]
{\includegraphics[width=4cm]{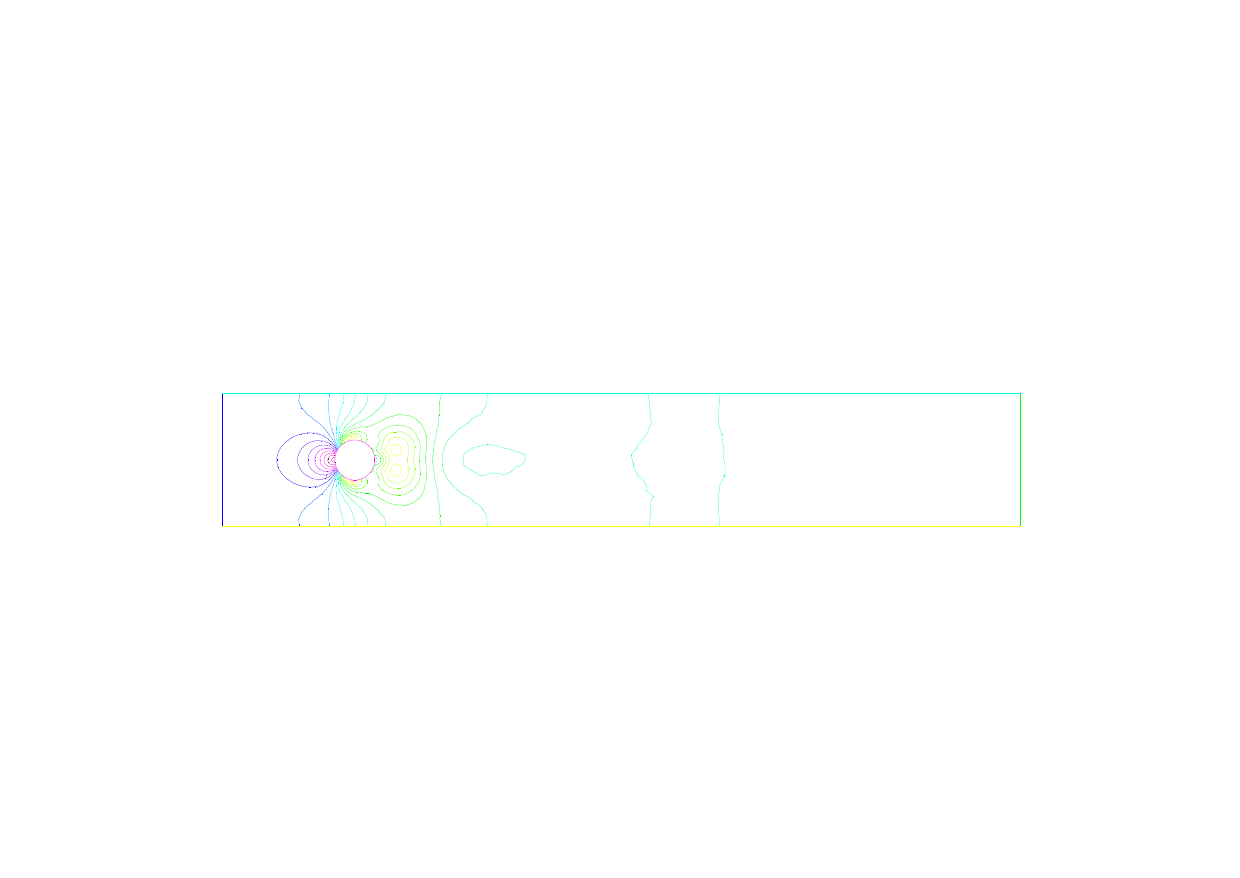}}\\
\subfigure[velocity (WG)]
{\includegraphics[width=4cm]{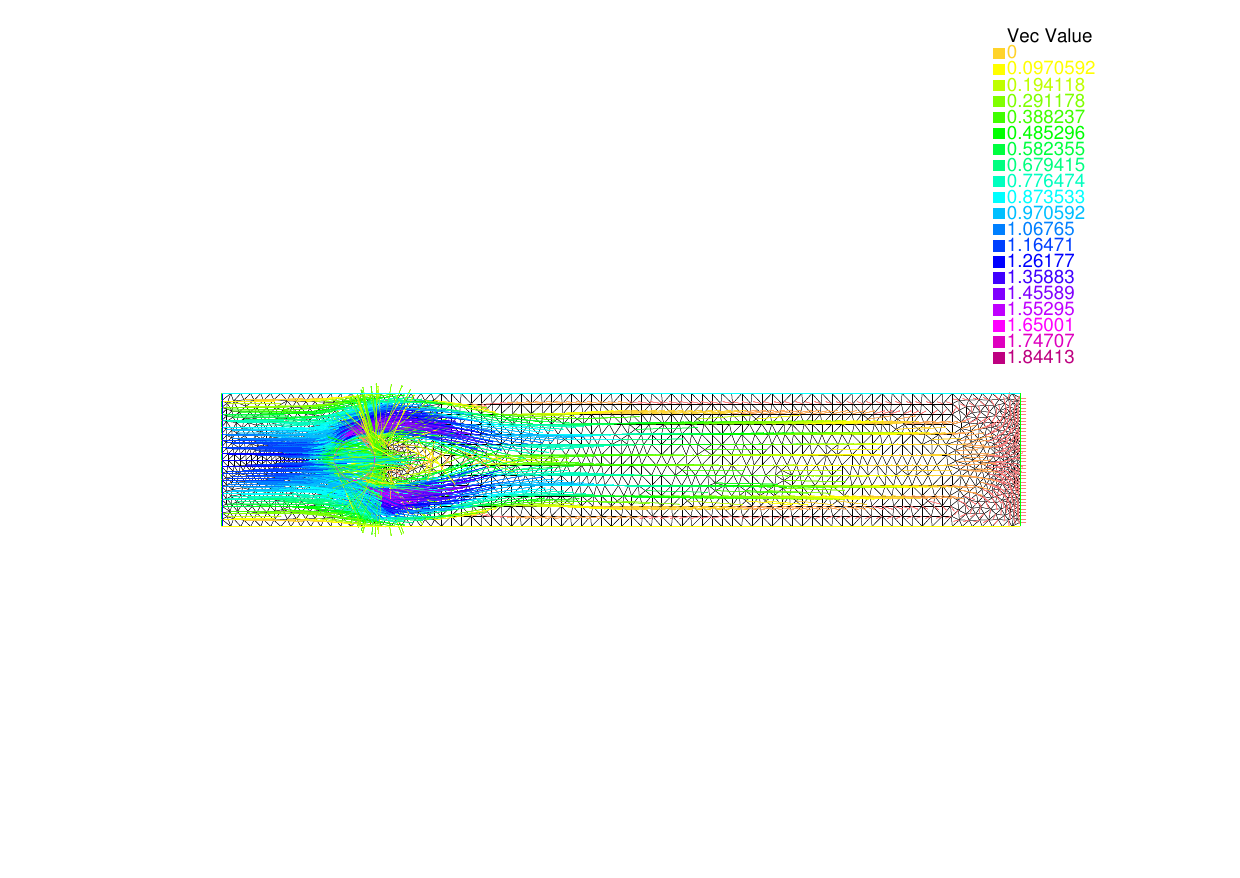}}
\subfigure[vorticity (WG)]
{\includegraphics[width=4cm]{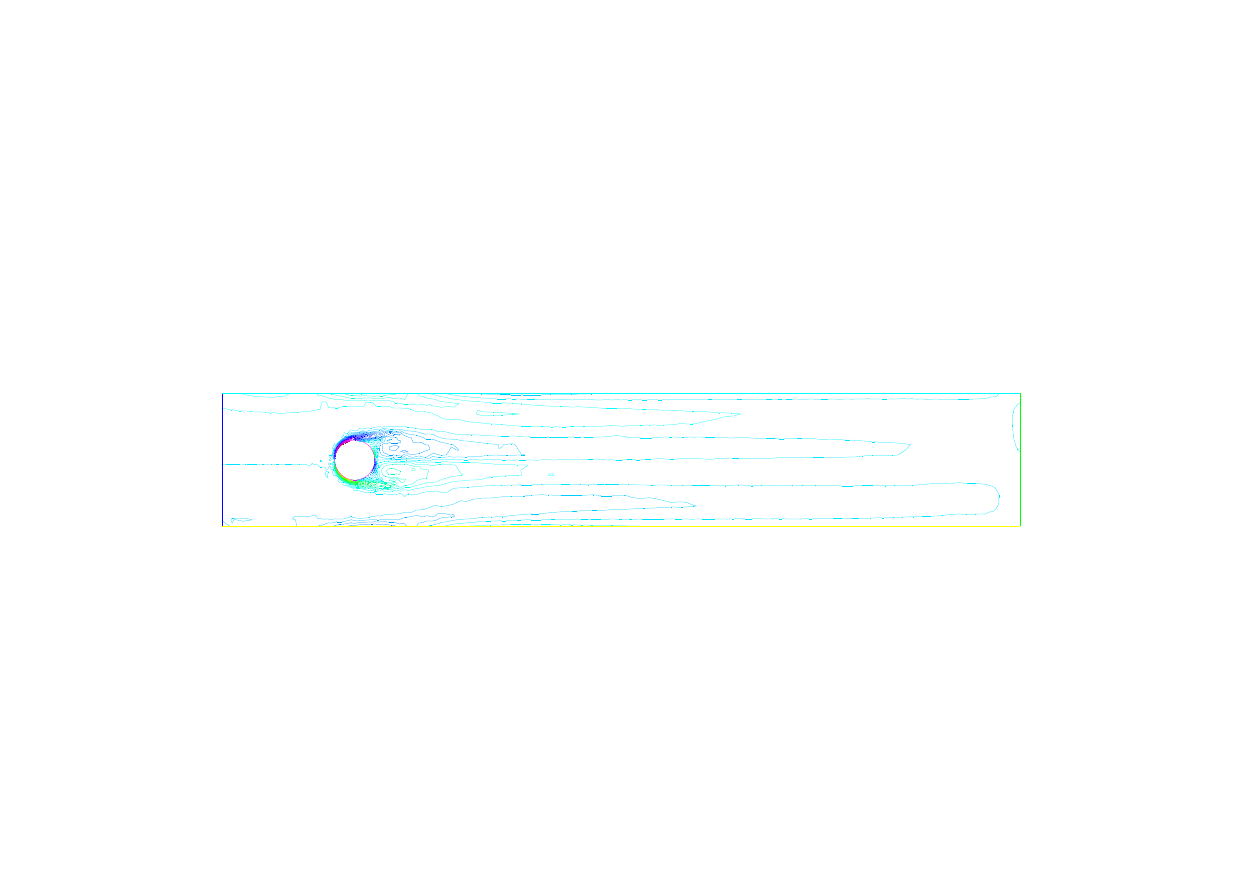}}
\subfigure[pressure (WG)]
{\includegraphics[width=4cm]{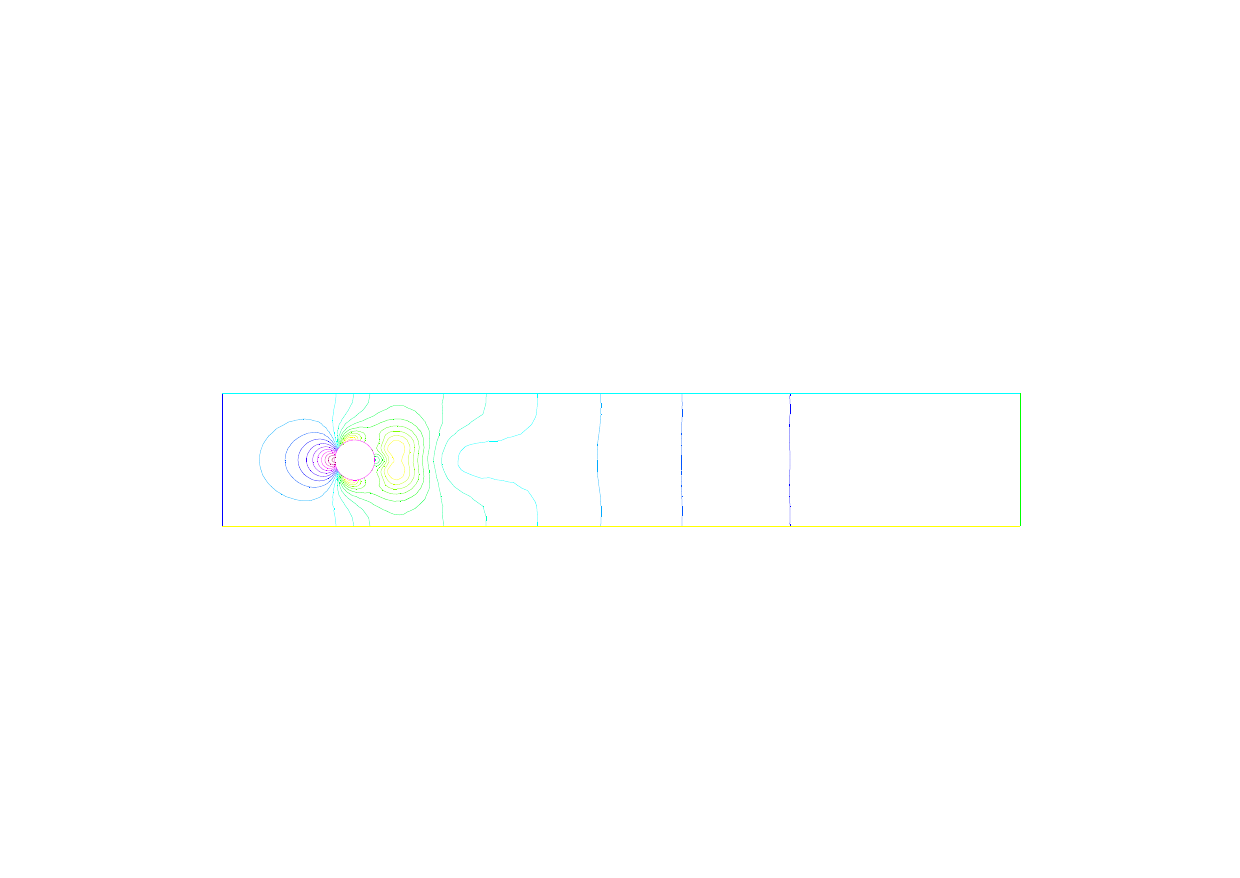}}
\caption{ The velocity streamlines, vortex lines  and pressure contours for Example \ref{EX7.3}: $\alpha=0$}
\label{fig21:2111}
\end{figure}

\begin{figure}[htbp!]
\centering
\subfigure[velocity: $\alpha=0.1$]
{\includegraphics[width=4cm]{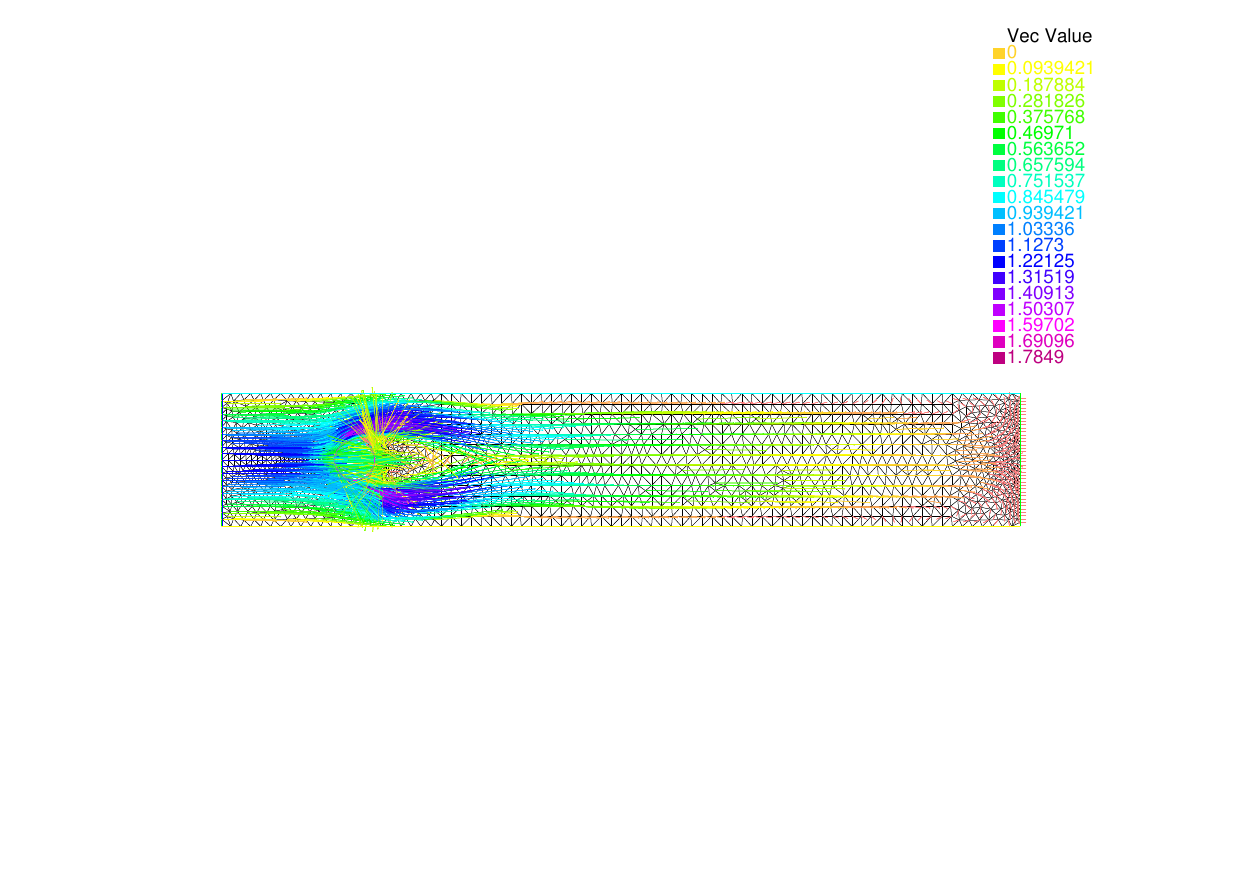}}
\subfigure[velocity: $\alpha=1$]
{\includegraphics[width=4cm]{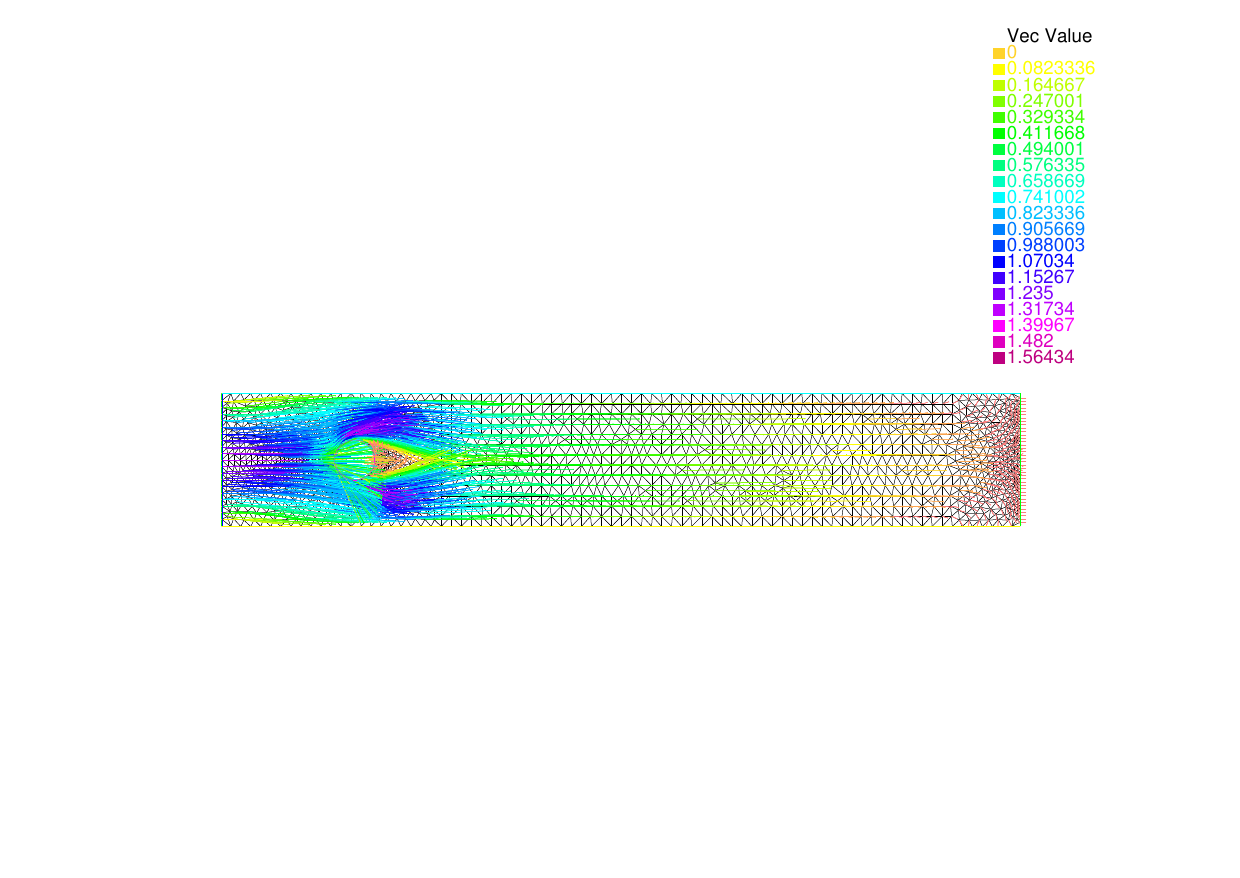}}
\subfigure[velocity: $\alpha=10$]
{\includegraphics[width=4cm]{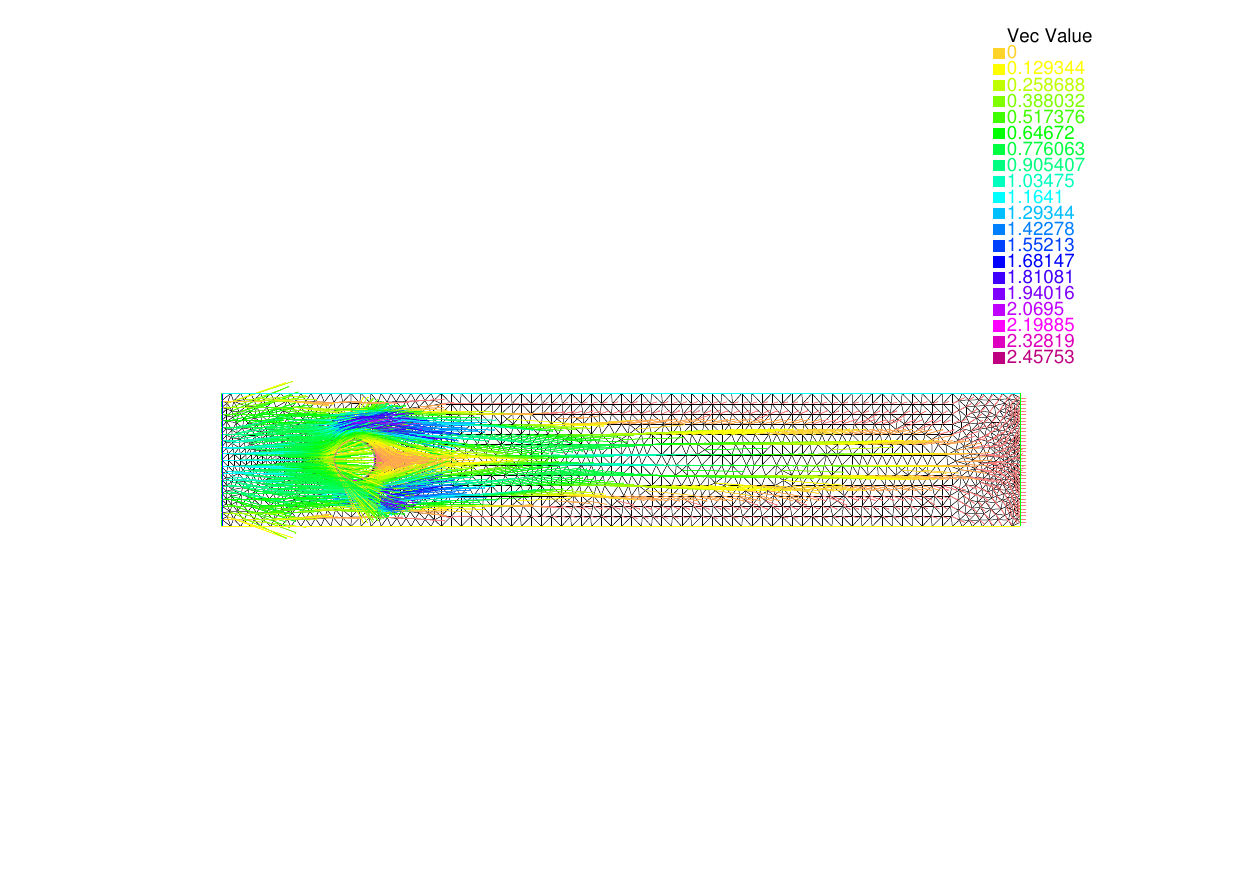}}
\quad
\quad
\subfigure[vorticity: $\alpha=0.1$]
{\includegraphics[width=4cm]{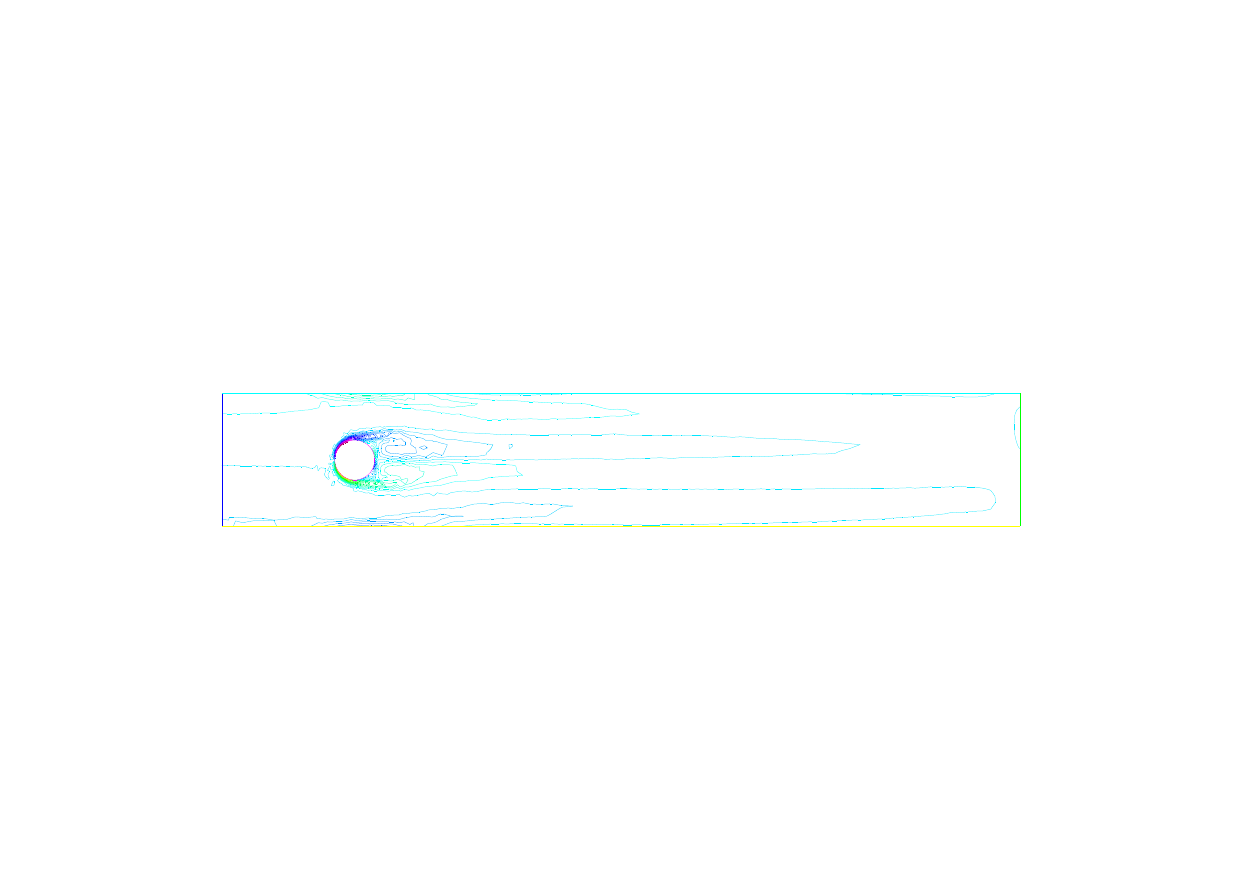}}
\subfigure[vorticity: $\alpha=1$]
{\includegraphics[width=4cm]{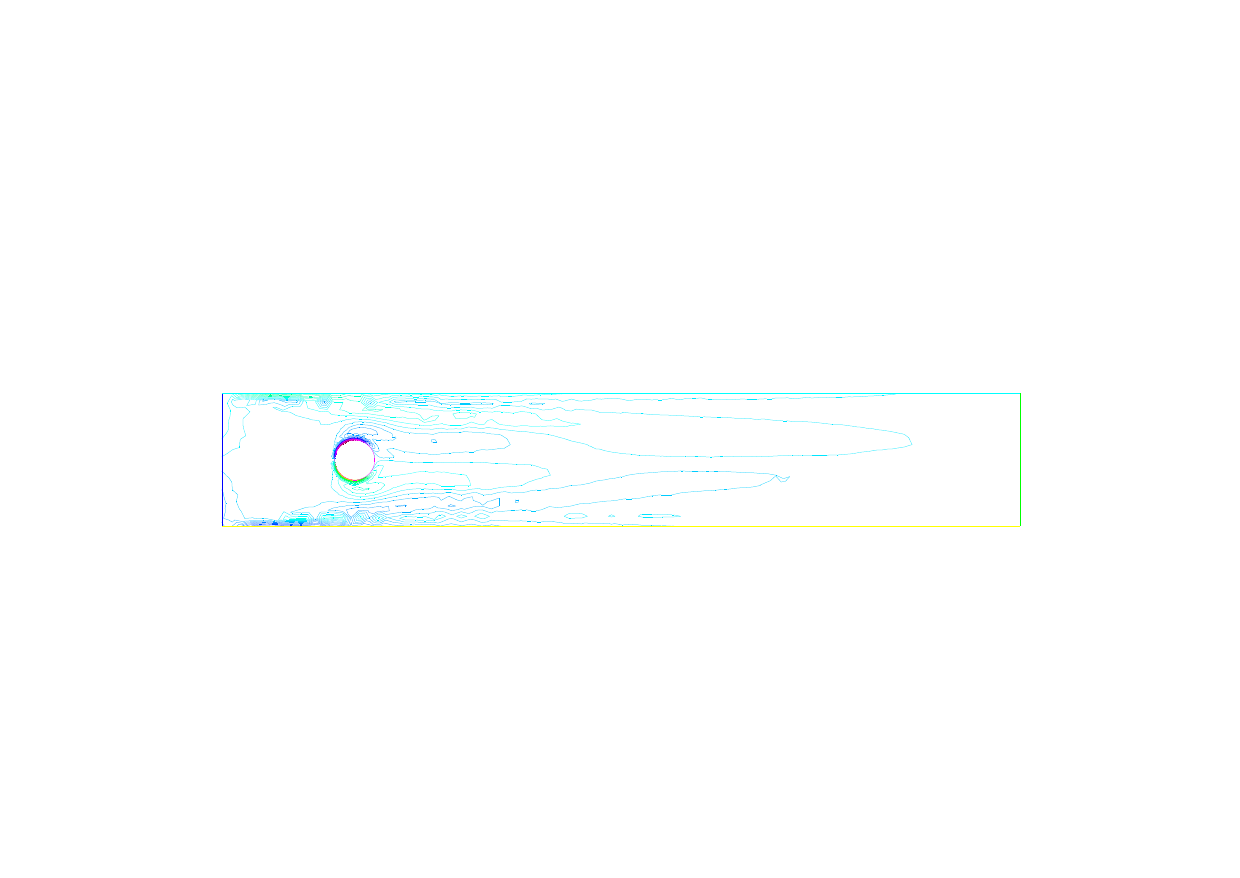}}
\subfigure[vorticity: $\alpha=10$]
{\includegraphics[width=4cm]{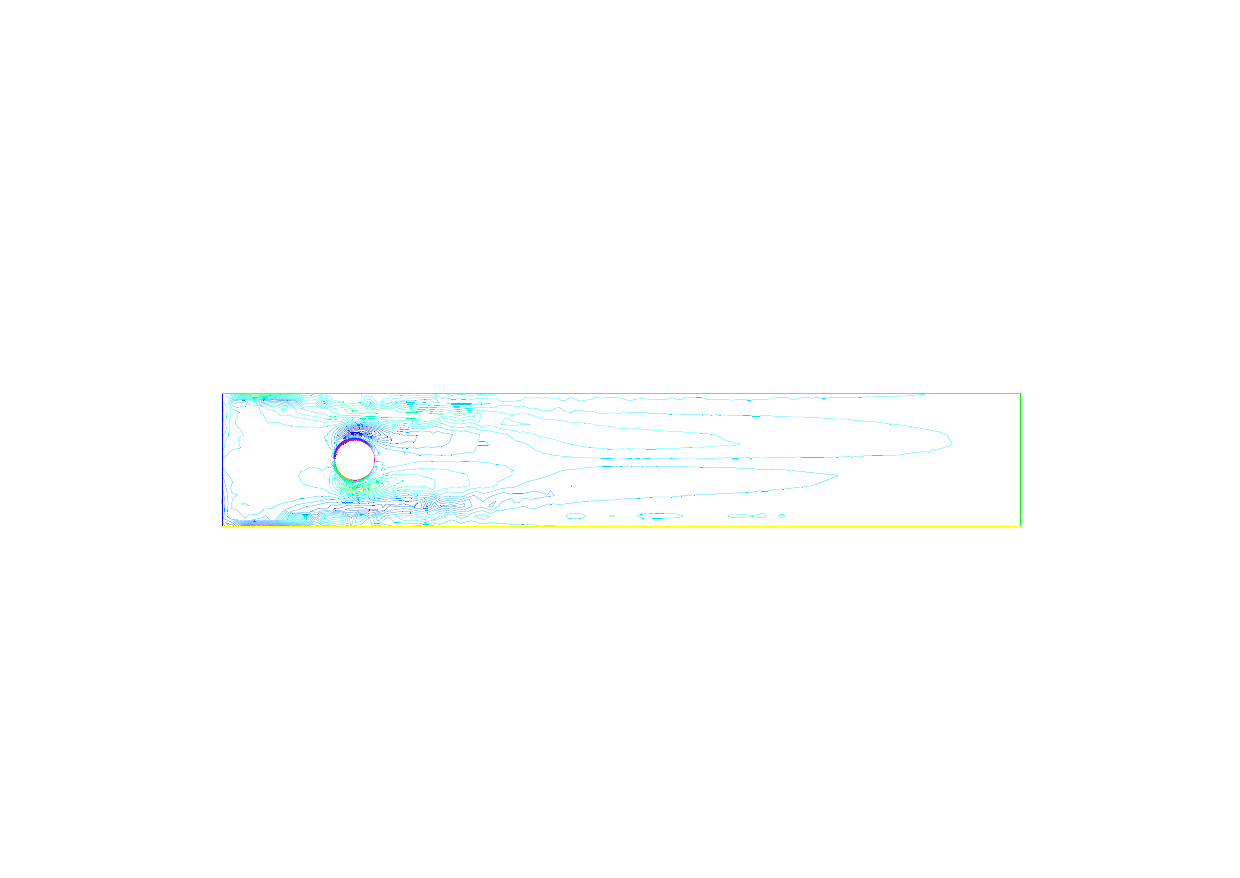}}\hspace{-10mm}
\quad
\subfigure[pressure: $\alpha=0.1$]
{\includegraphics[width=4cm]{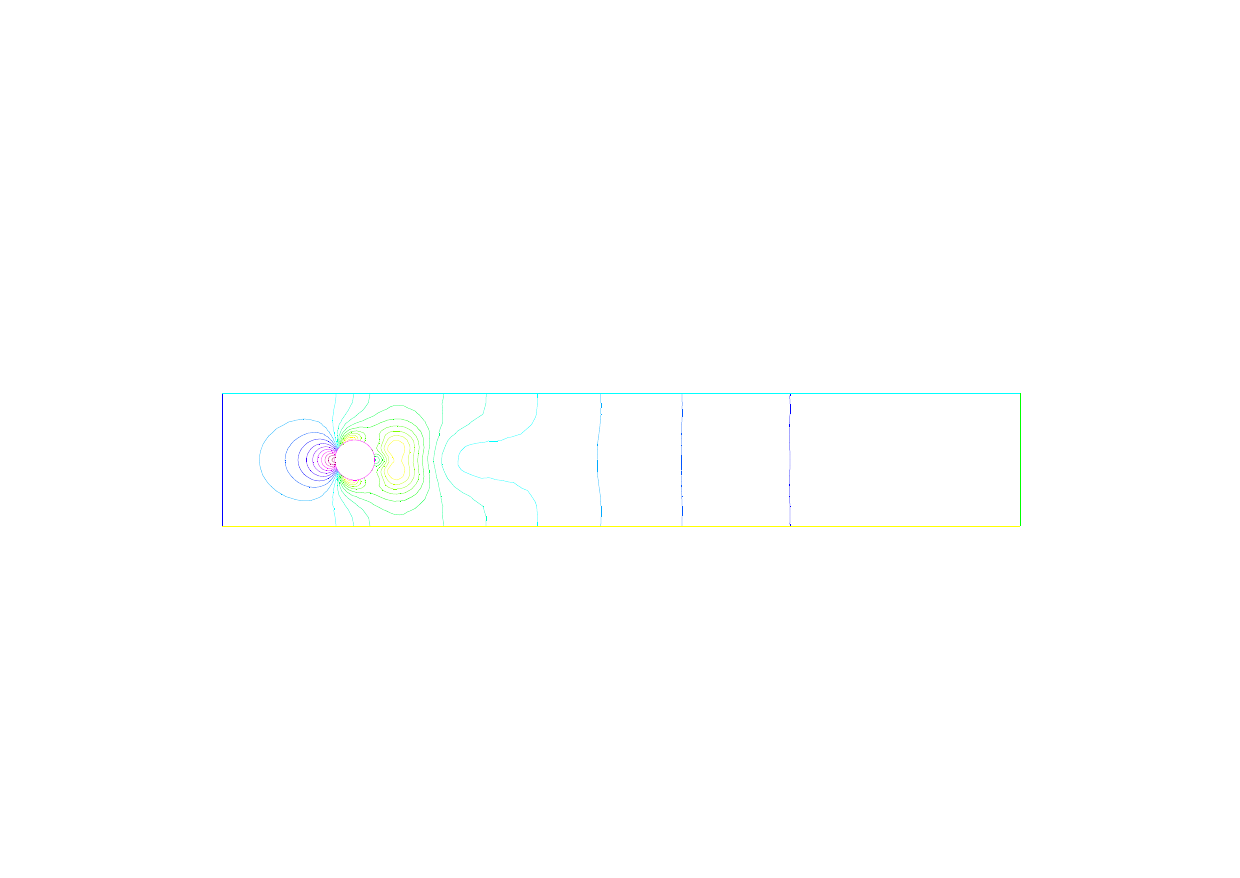}}
\subfigure[pressure: $\alpha=1$]
{\includegraphics[width=4cm]{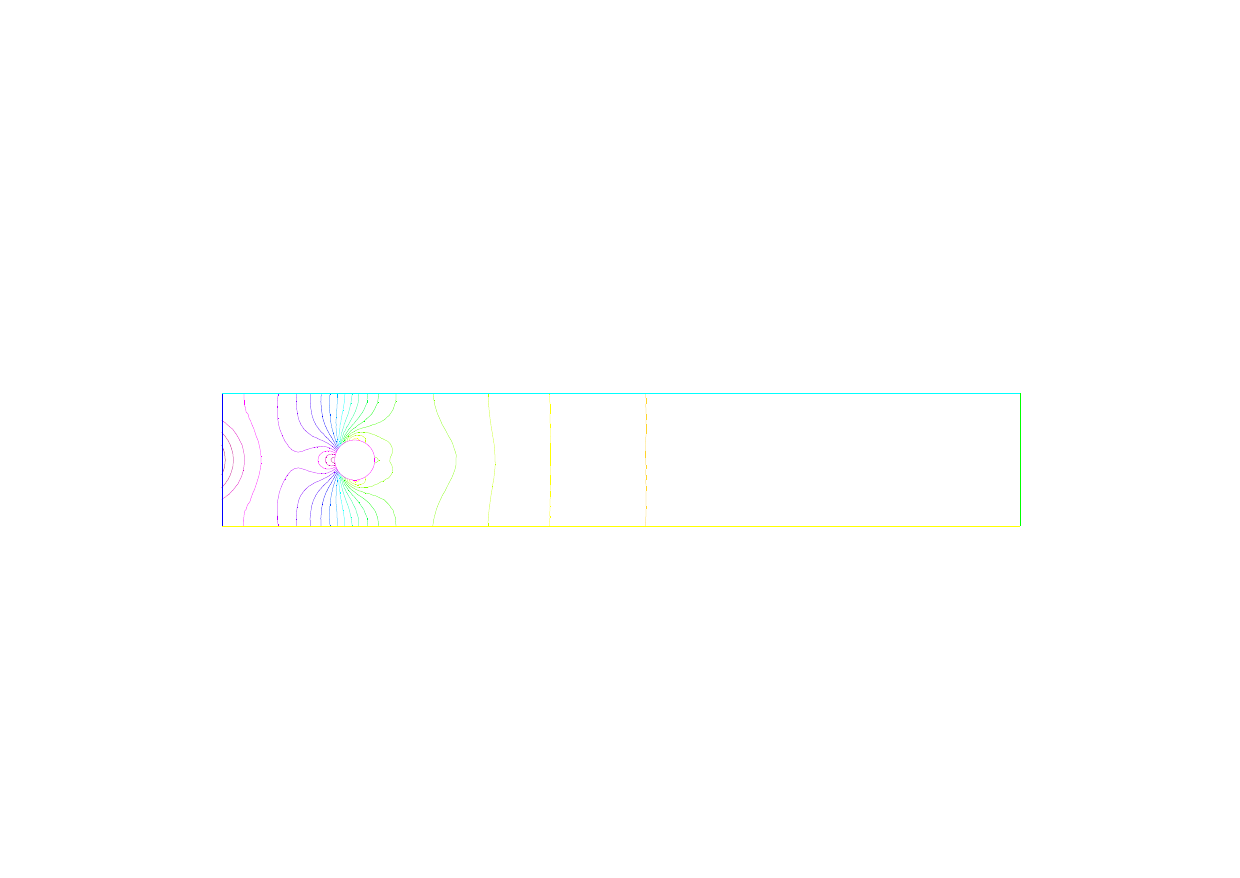}}
\subfigure[pressure: $\alpha=10$]
{\includegraphics[width=4cm]{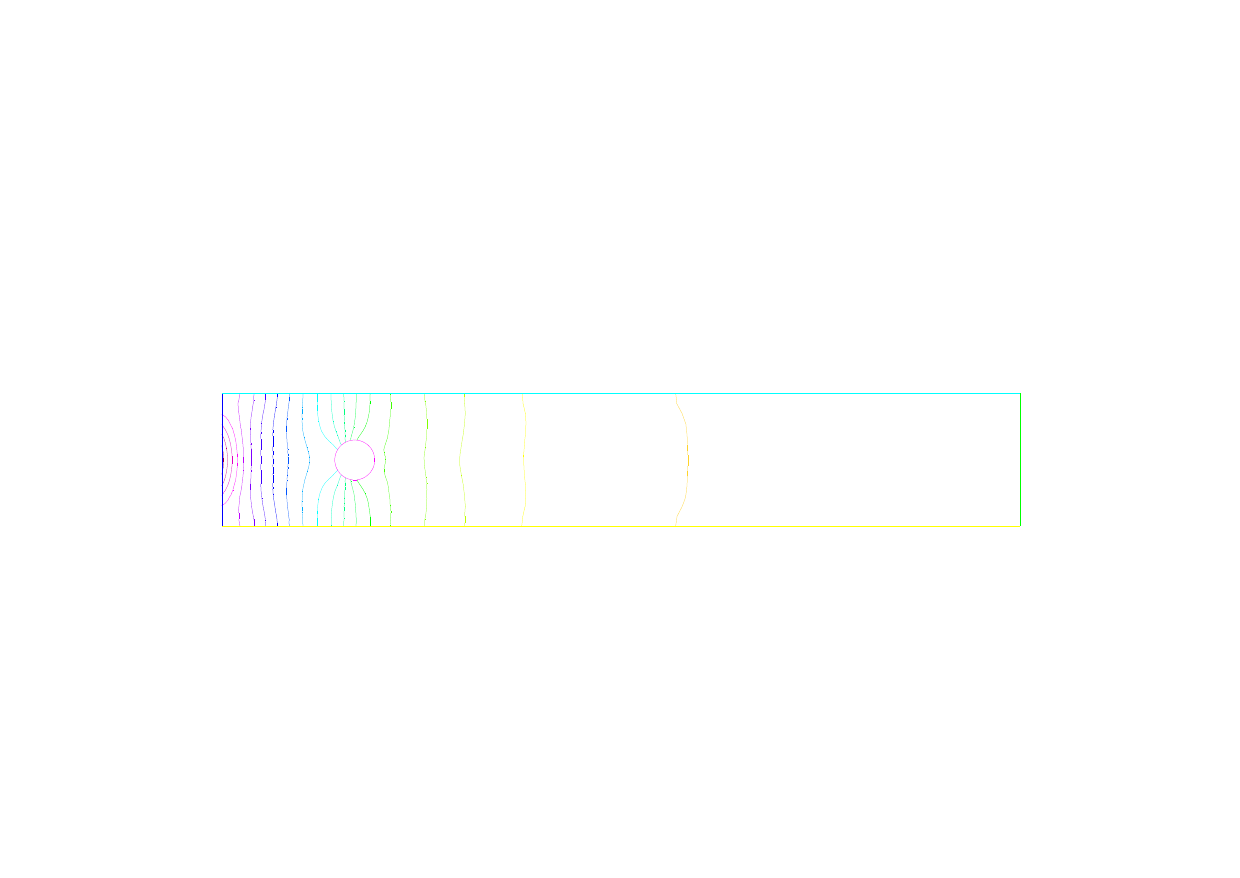}}
\caption{The velocity streamlines, vortex lines  and pressure contours for Example \ref{EX7.3}: $r=3.5$ and $\alpha=0.1, 1, 10$ }
\label{fig31:31}
\end{figure}

\begin{figure}[htbp!]
\centering
\setlength{\abovecaptionskip}{0.cm}
\subfigure[velocity: $r=3$]
{\includegraphics[width=4cm]{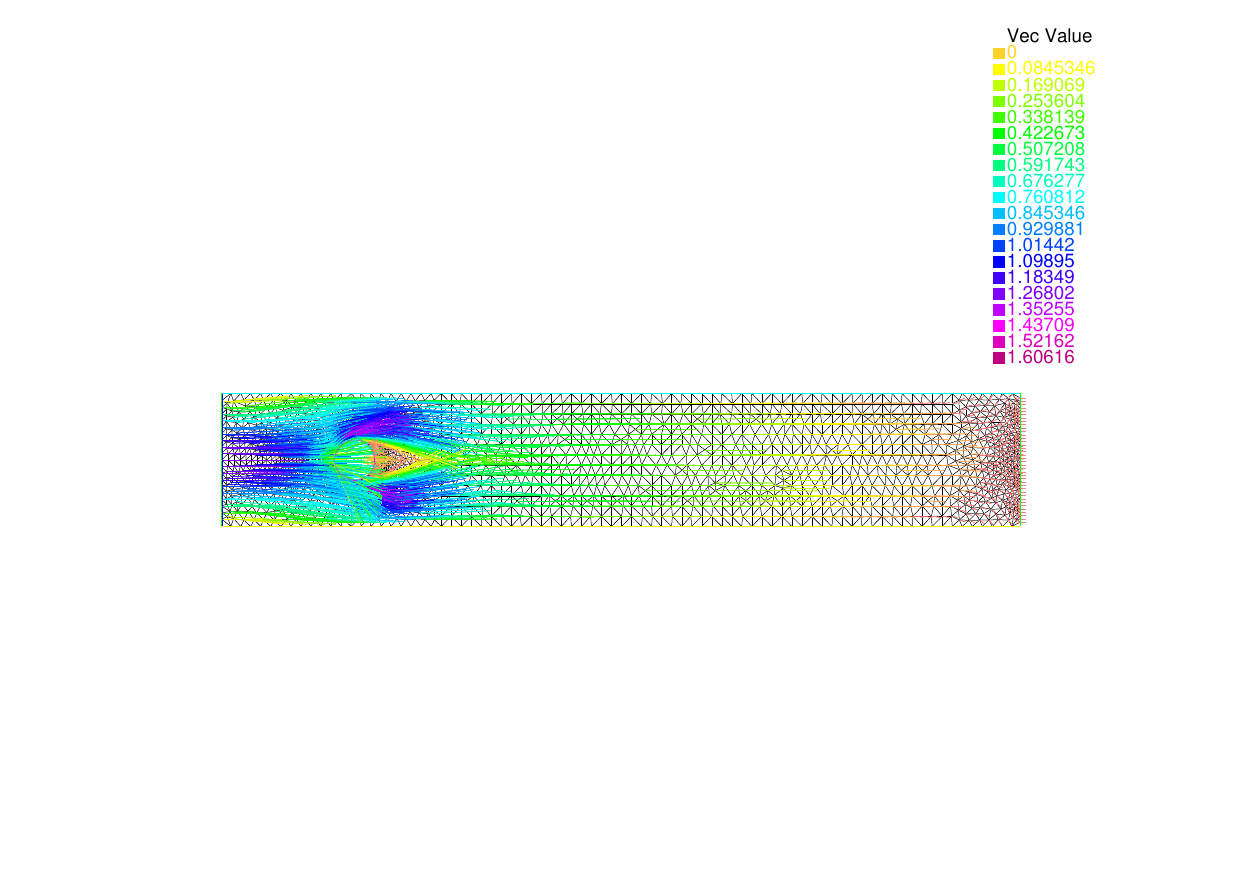}}
\subfigure[velocity: $r=4$]
{\includegraphics[width=4cm]{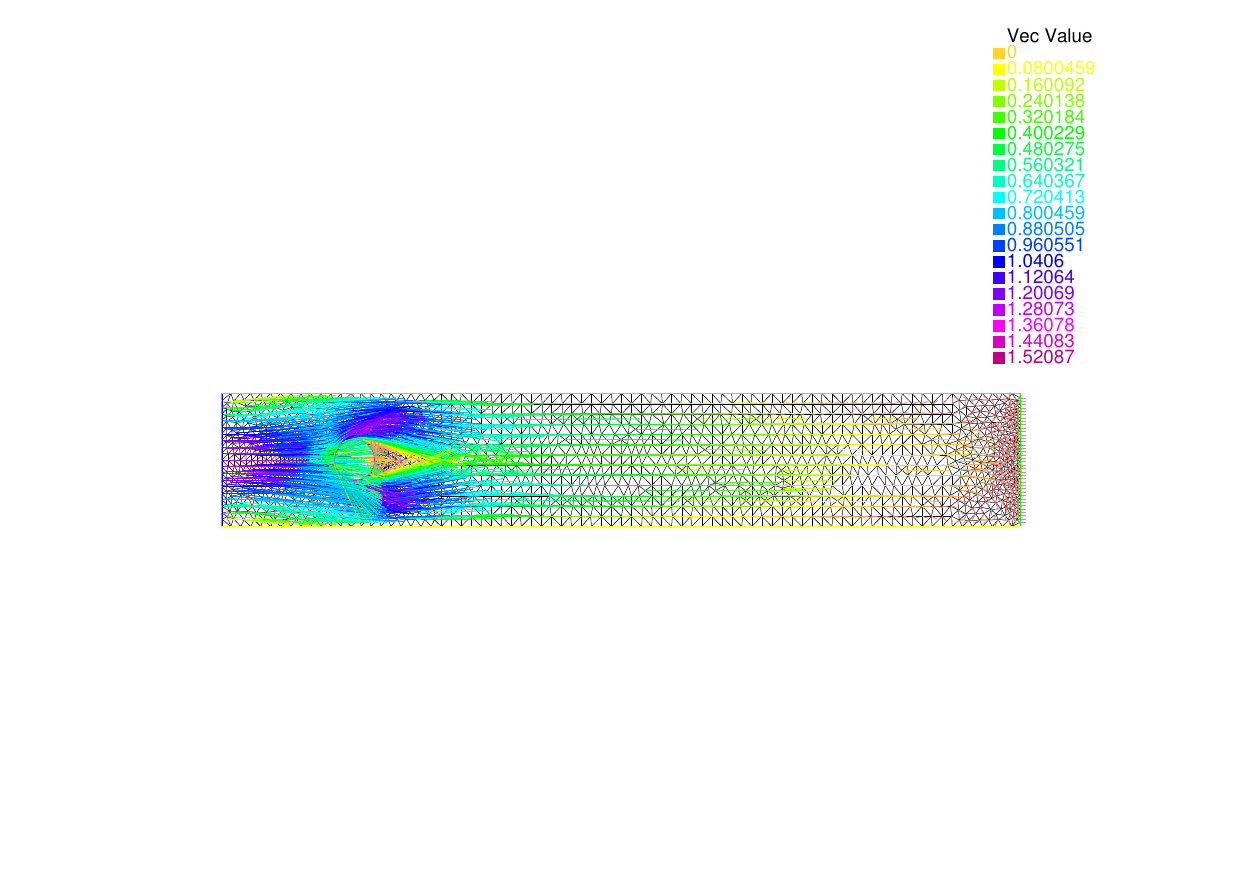}}
\subfigure[velocity: $r=5$]
{\includegraphics[width=4cm]{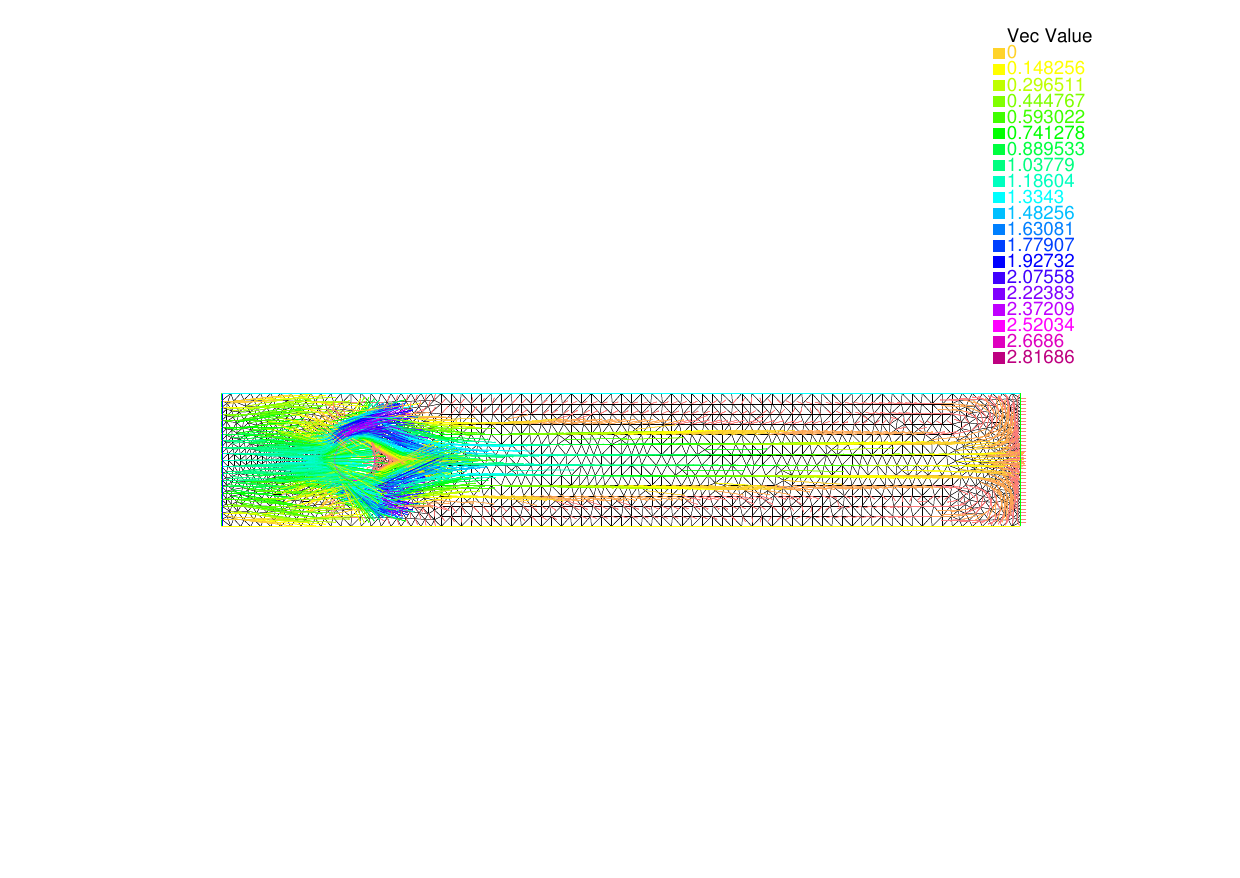}}
\quad
\subfigure[vorticity: $r=3$]
{\includegraphics[width=4cm]{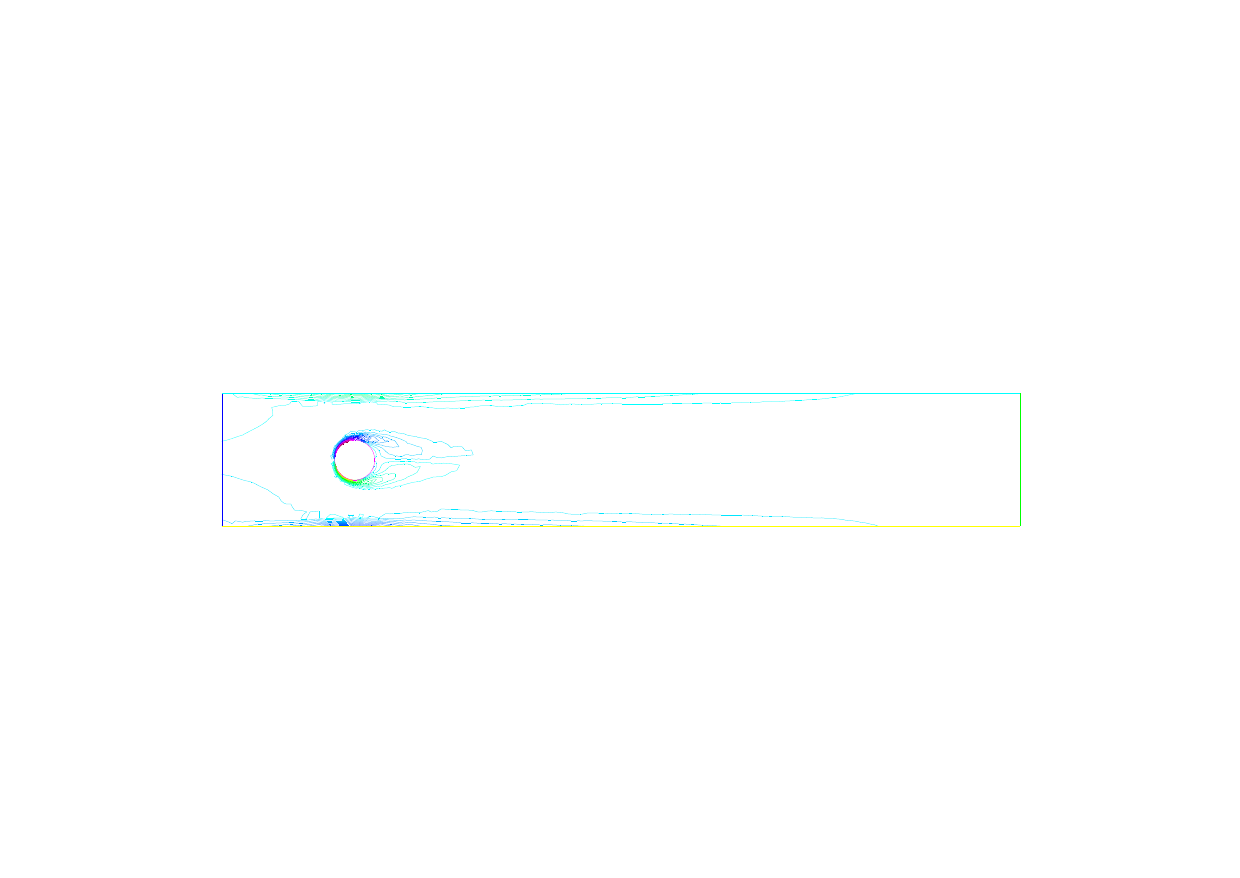}}
\subfigure[vorticity: $r=4$]
{\includegraphics[width=4cm]{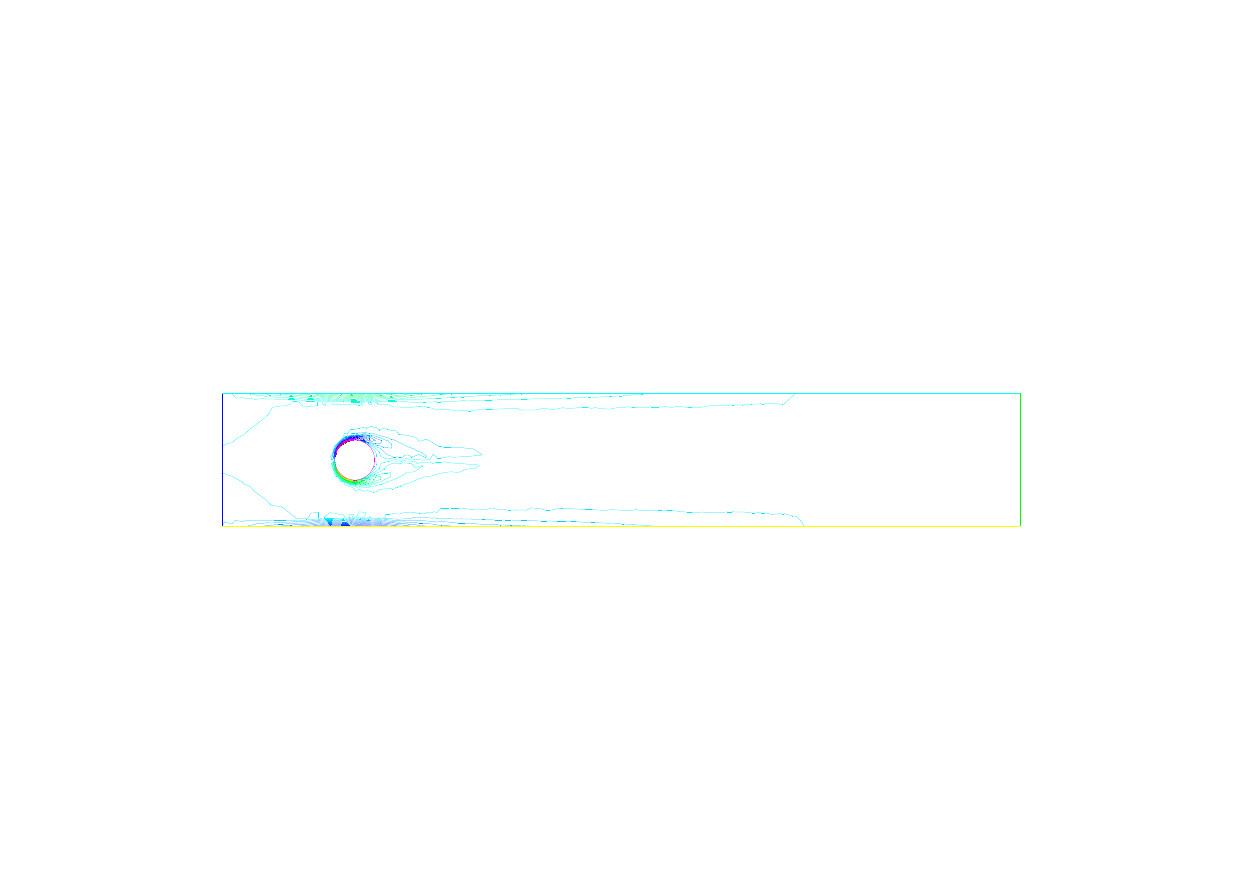}}
\subfigure[vorticity: $r=5$]
{\includegraphics[width=4cm]{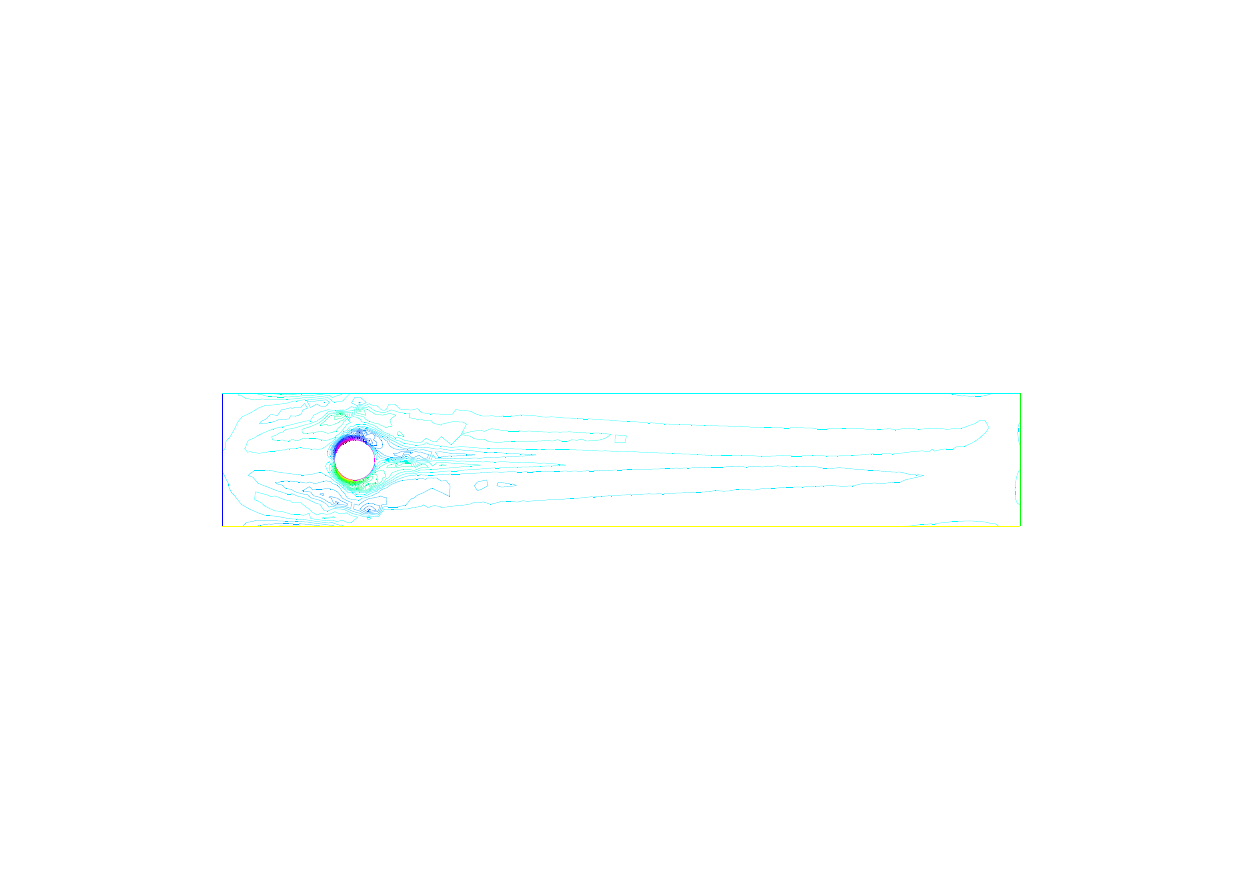}}
\quad
\subfigure[pressure: $r=3$]
{\includegraphics[width=4cm]{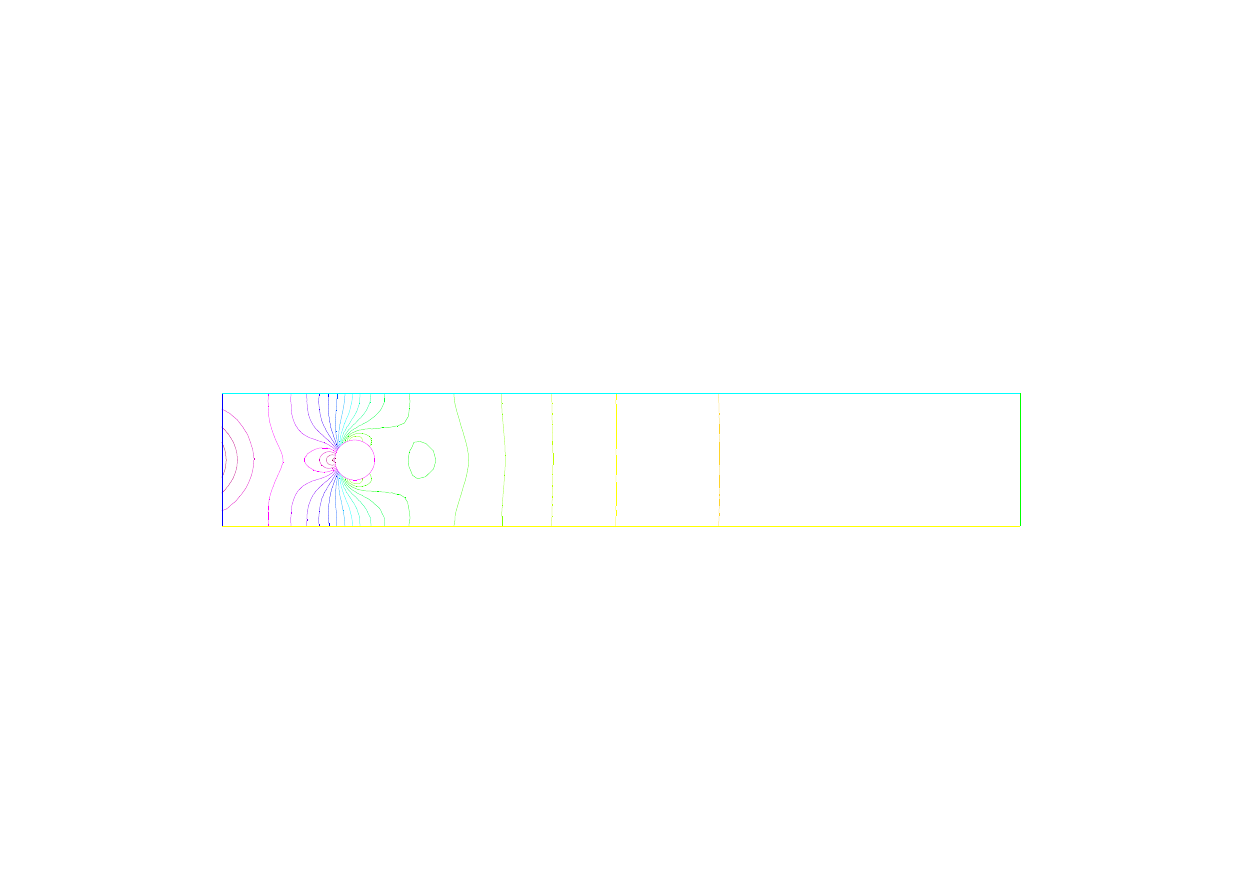}}
\subfigure[pressure: $r=4$]
{\includegraphics[width=4cm]{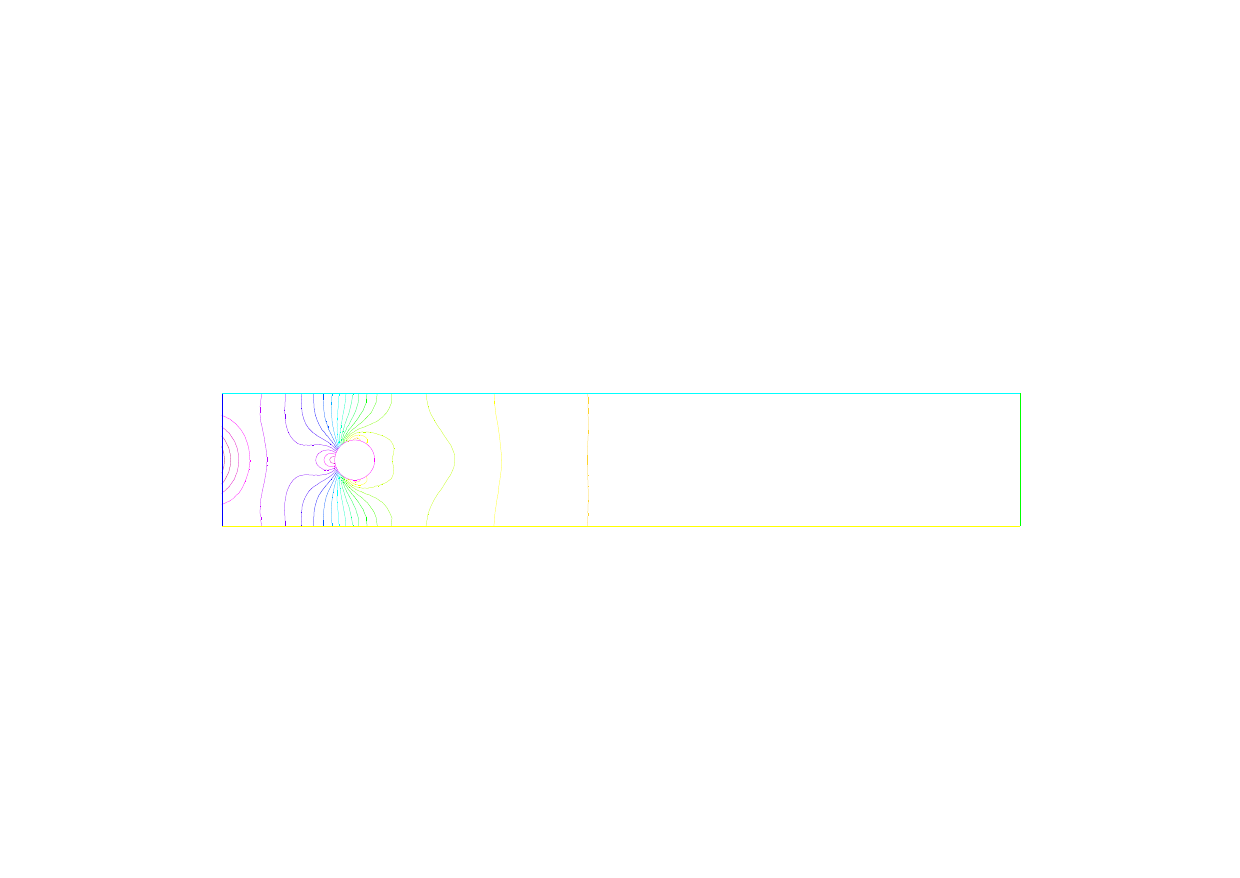}}
\subfigure[pressure: $r=5$]
{\includegraphics[width=4cm]{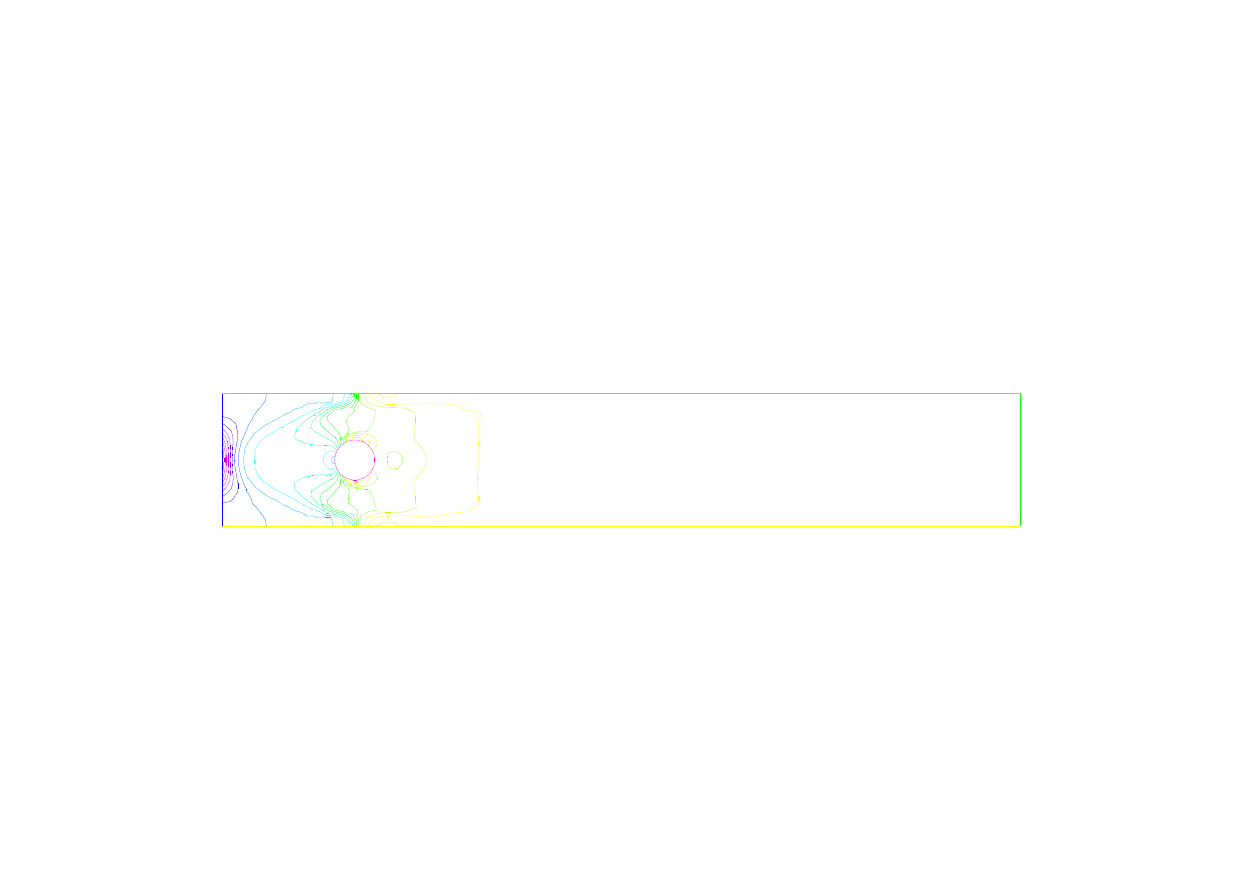}}
\caption{The velocity streamlines, vortex lines  and pressure contours for Example \ref{EX7.3}: $\alpha=1$ and  $r=3,4,5$}
\label{fig32:32}
\end{figure}

\begin{exam}[The problem of flow around a   circular  cylinder]\label{EX7.3}
The flow around a circular cylinder   is examined with the Brinkman-Forchheimer model \eqref{BF0} and the the WG  method. We take
  $\Omega$ =$[0, 6]\times[0, 1] \setminus O_{d}(1,0.5)$, $\nu=0.002$   and $\bm{f}=\bm{0}$, where $O_{d}(1,0.5)$ is a disk with center $(1,0.5)$  and diameter $d=0.3$; see  Figure \ref{fig2:mesh2} for the domain and its finite element mesh. The boundary conditions are as follows:
$$\bm{u}|_{y=0}=\bm{u}|_{y=1}=\bm{u}|_{\partial O_{d}}=\bm{0}, \quad \bm{u}|_{x=0}=(6y(1-y), 0 )^T, $$
$$ \left(-p\bm{I}+\nu\nabla \bm{u}\right){  \bm{n}}\big|_{x=6}=0, $$
where $\bm{I}$ and $ \bm{n}$ are the unit matrix and the outward unit normal vector, respectively.
We compute the WG scheme \eqref{WG} with  $m=k=2$ in the following cases:
\begin{itemize}
\item [ I]. $\alpha=0$, i.e. the case of the Navier-Stokes  equations;

\item [ II].   $r=3.5$ and $\alpha=0.1, 1, 10$;

\item [ III]. $\alpha=1$ and $r=3, 4, 5$.
\end{itemize}
The obtained velocity, vorticity and pressure approximations are shown   in Figures \ref{fig21:2111}, \ref{fig31:31} and \ref{fig32:32}, respectively.  As a comparison,  the  referenced  numerical solutions  obtained  with the Taylor-Hood  element 
are also shown for $\alpha=0$; see (a), (b)  and (c) in Figure \ref{fig21:2111}.
We can see that our method is effective and the damping effect is gradually enhanced as the parameters $\alpha$ and $r$   increase.

\end{exam}

\begin{figure}[htbp!]
\centering
\subfigure
{\includegraphics[width=12cm]{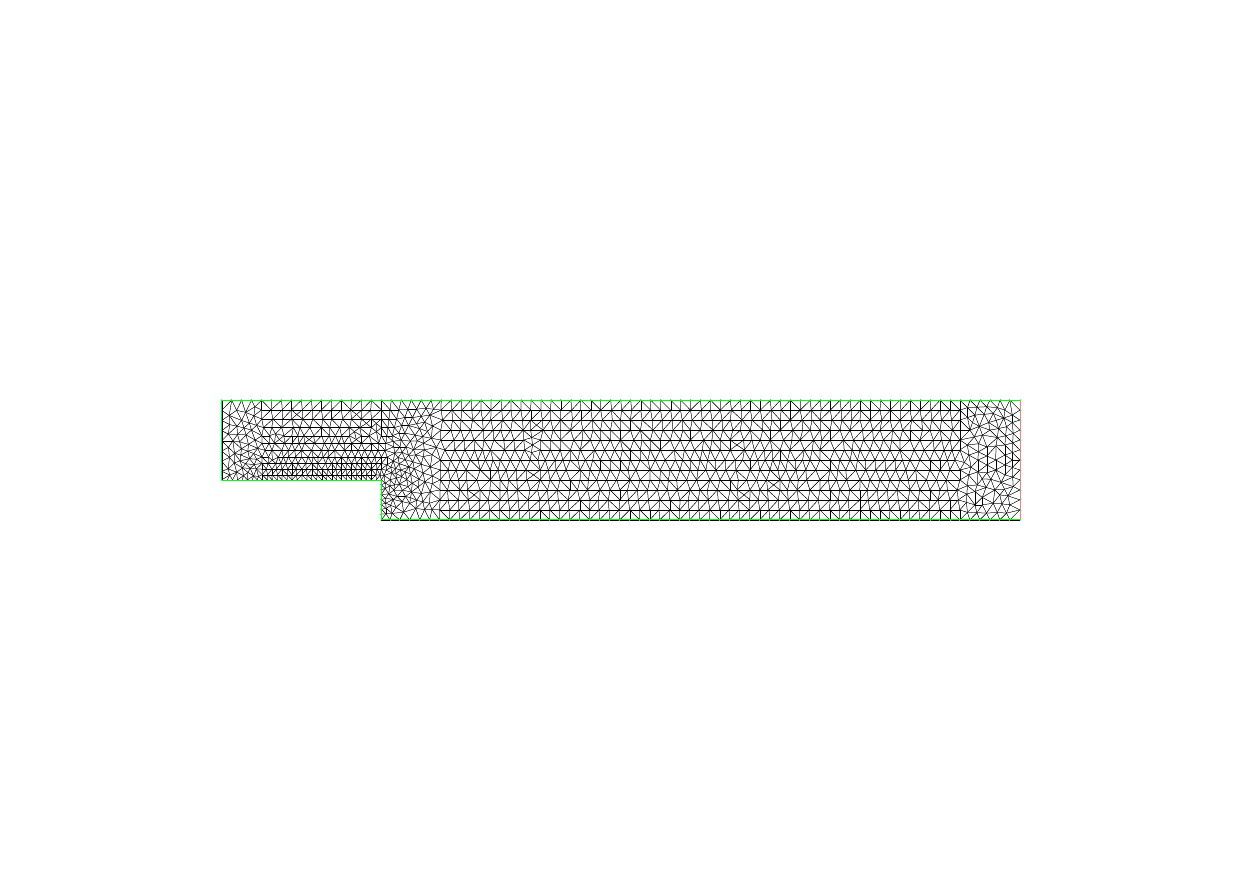}}
\caption{The domain and  finite element mesh for Example \ref{EX7.4}}
\label{fig4:mesh3}
\end{figure}

\begin{figure}[htbp!]
\centering
\setlength{\abovecaptionskip}{0.cm}
\subfigure[ $u_1$ (Taylor-Hood) ]
{\includegraphics[width=4cm]{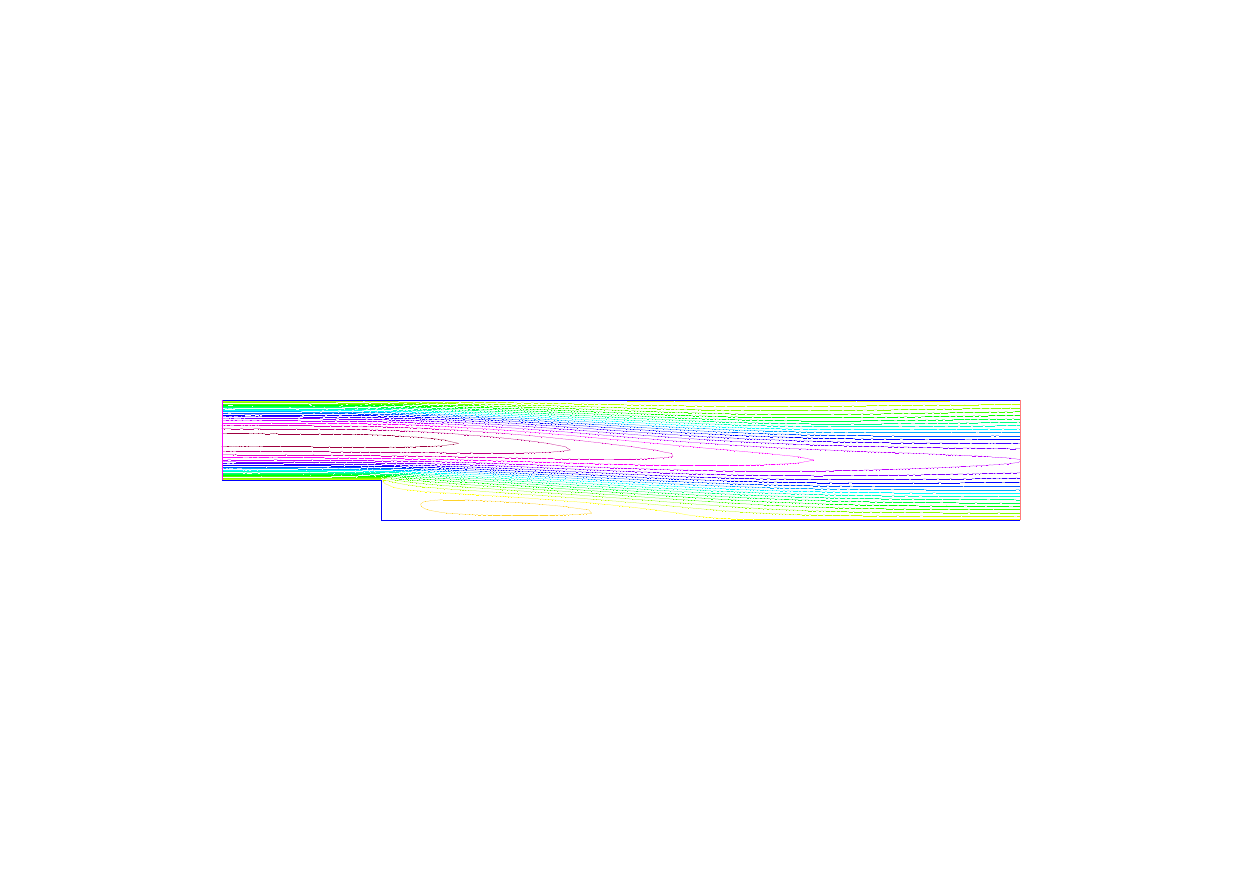}}
\subfigure[  $u_2$ (Taylor-Hood)]
{\includegraphics[width=4cm]{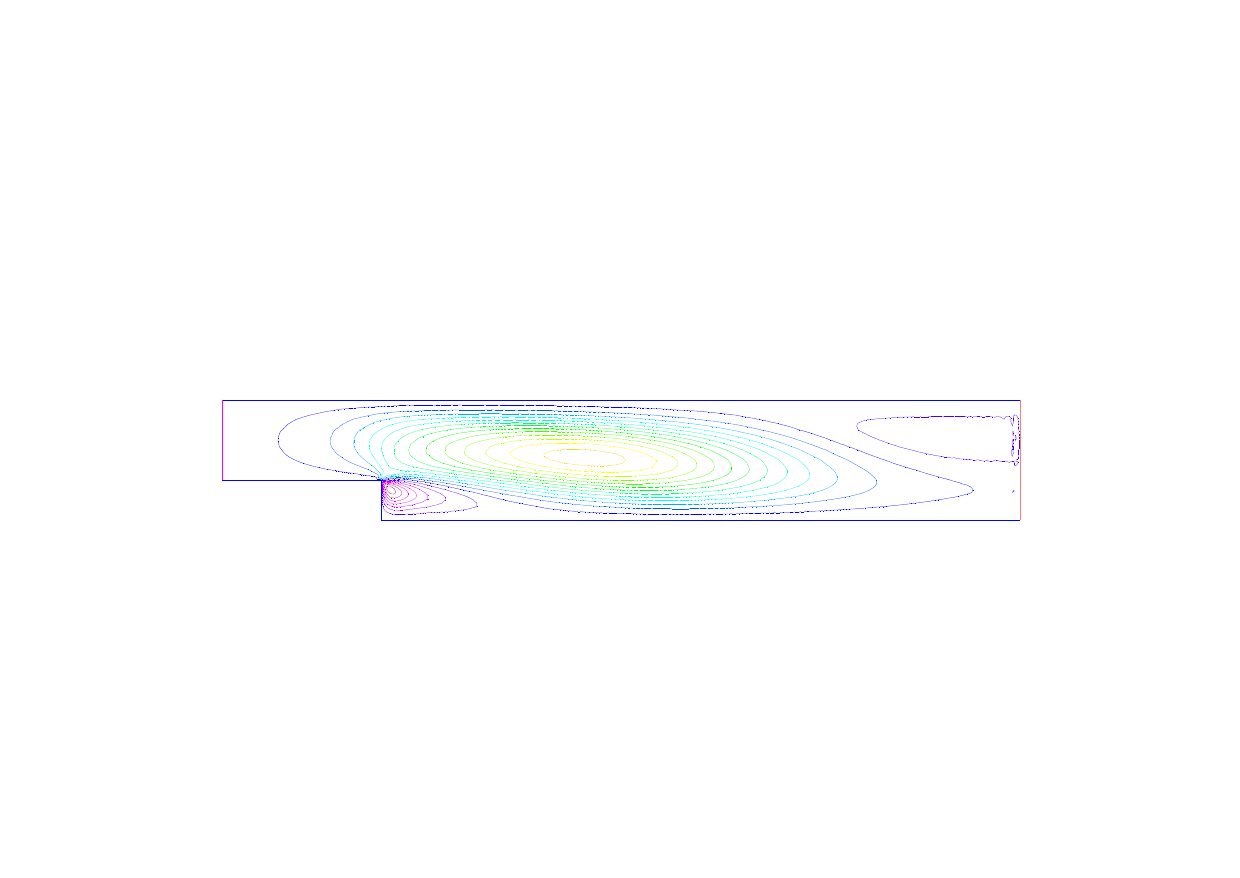}}
\subfigure[pressure (Taylor-Hood)]
{\includegraphics[width=4cm]{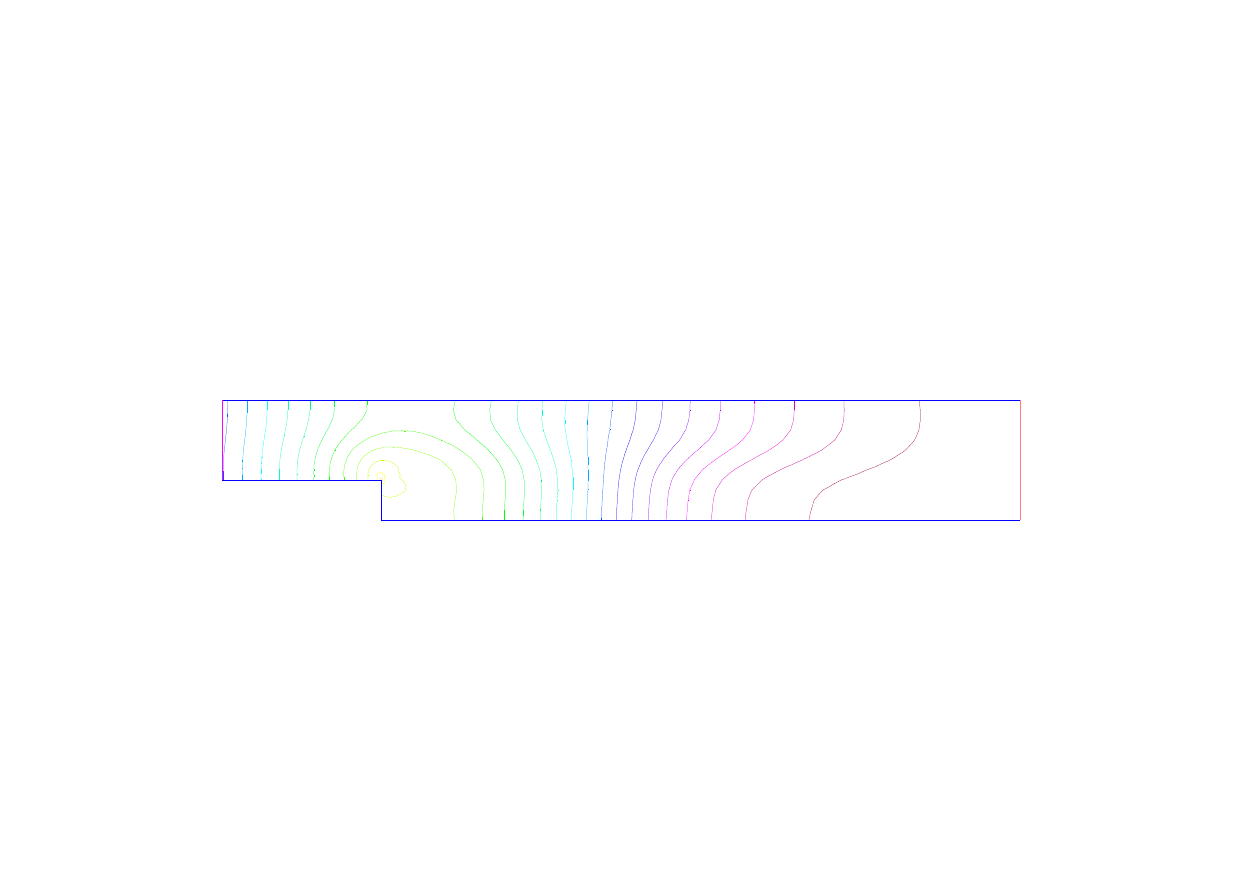}}\\
\subfigure[  $u_1$ (WG)]
{\includegraphics[width=3.9cm]{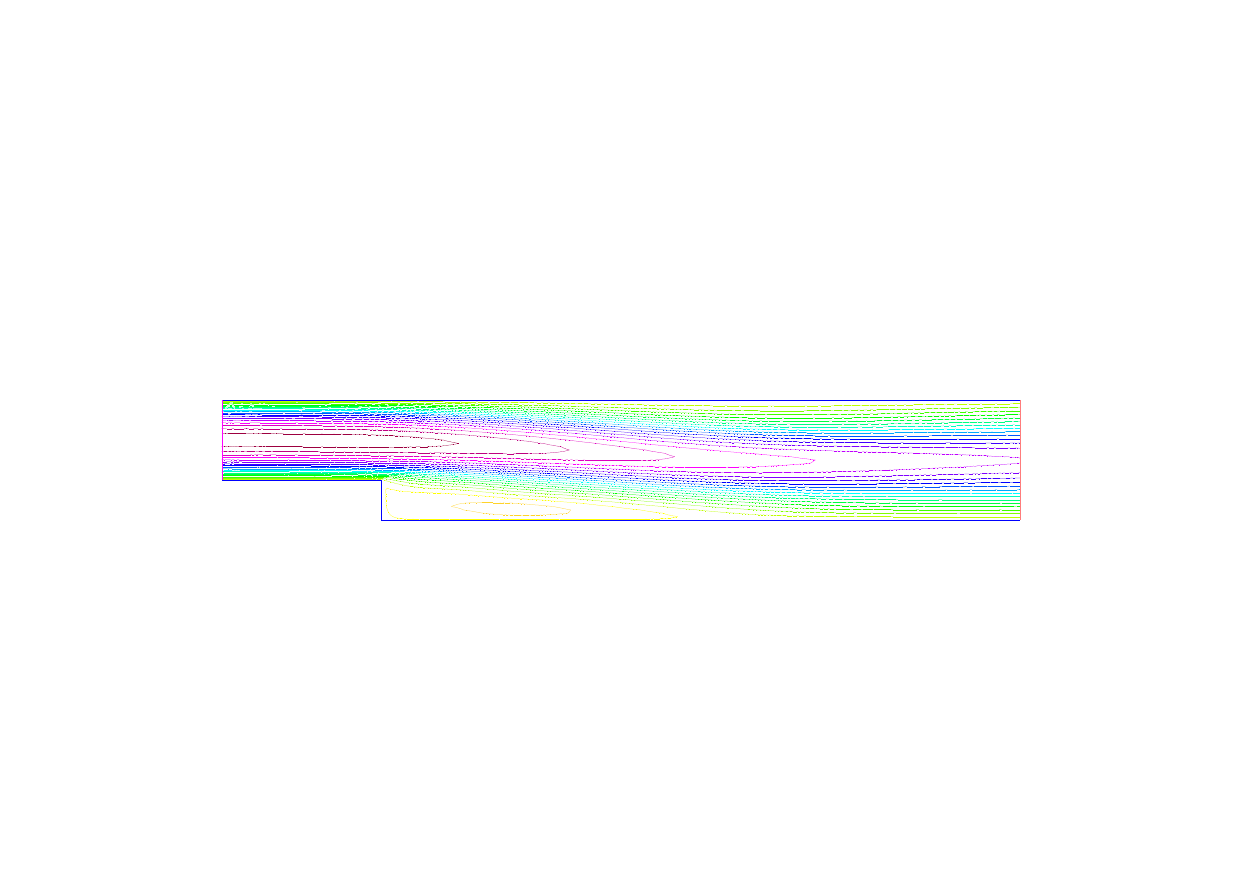}}
\subfigure[ $u_2$ (WG)]
{\includegraphics[width=3.9cm]{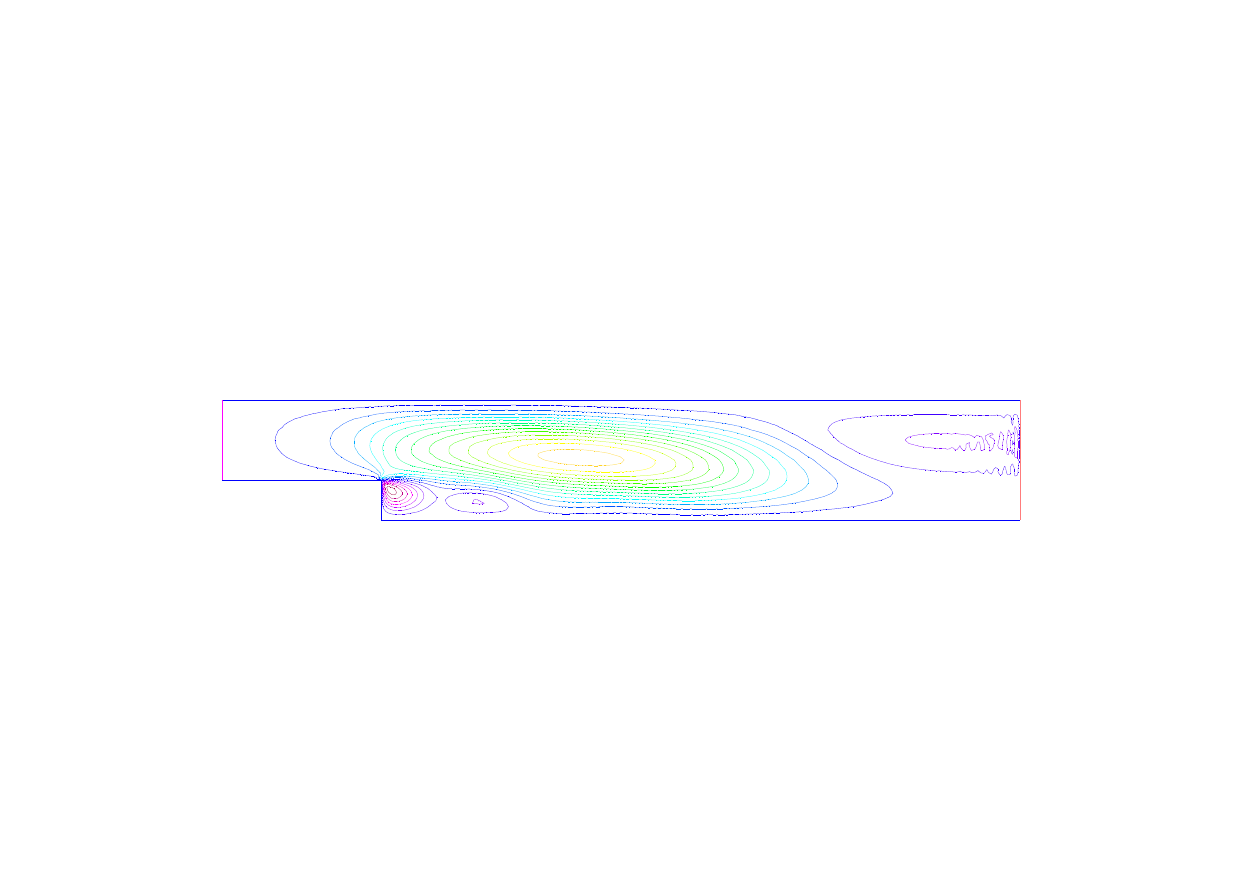}}
\subfigure[pressure (WG)]
{\includegraphics[width=3.9cm]{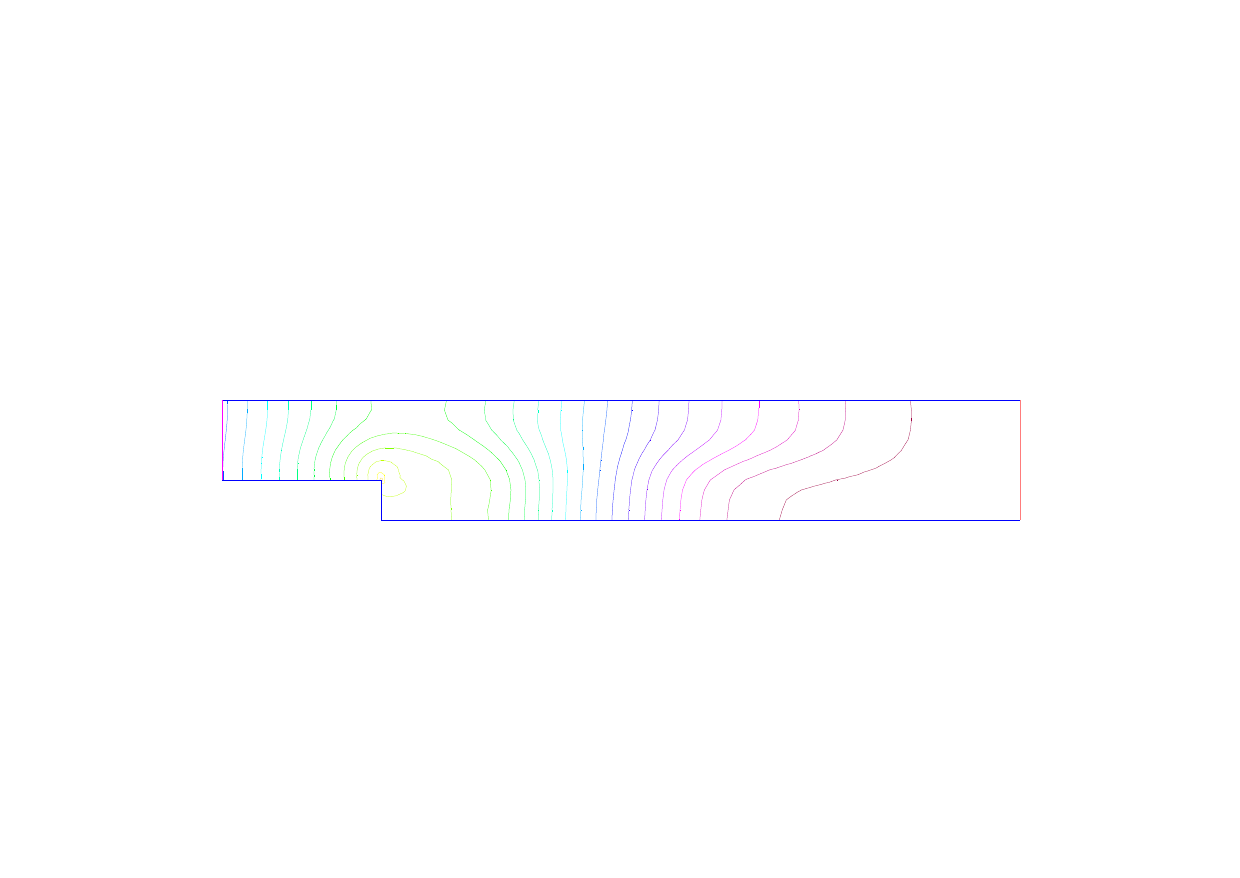}}
 \caption{ The velocity $\bm u_h=(u_1,u_2)^T$   and pressure contours for Example \ref{EX7.4}: $\alpha=0$}
\label{fig41:40}
\end{figure}

\begin{figure}[htbp!]
 \centering
\subfigure[  $u_1$: $\alpha=0.01$]
{\includegraphics[width=4cm]{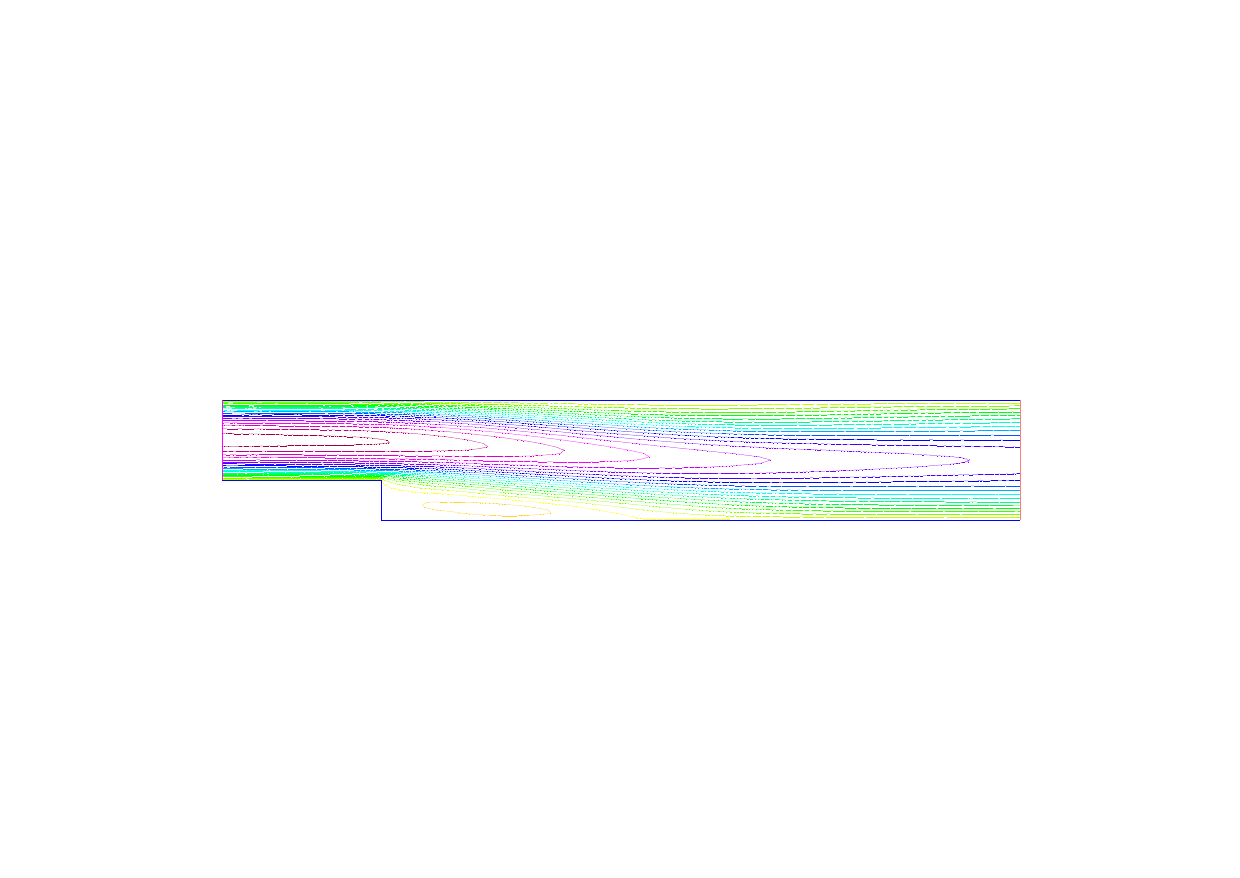}}
\subfigure[ $u_1$: $\alpha=0.1$]
{\includegraphics[width=4cm]{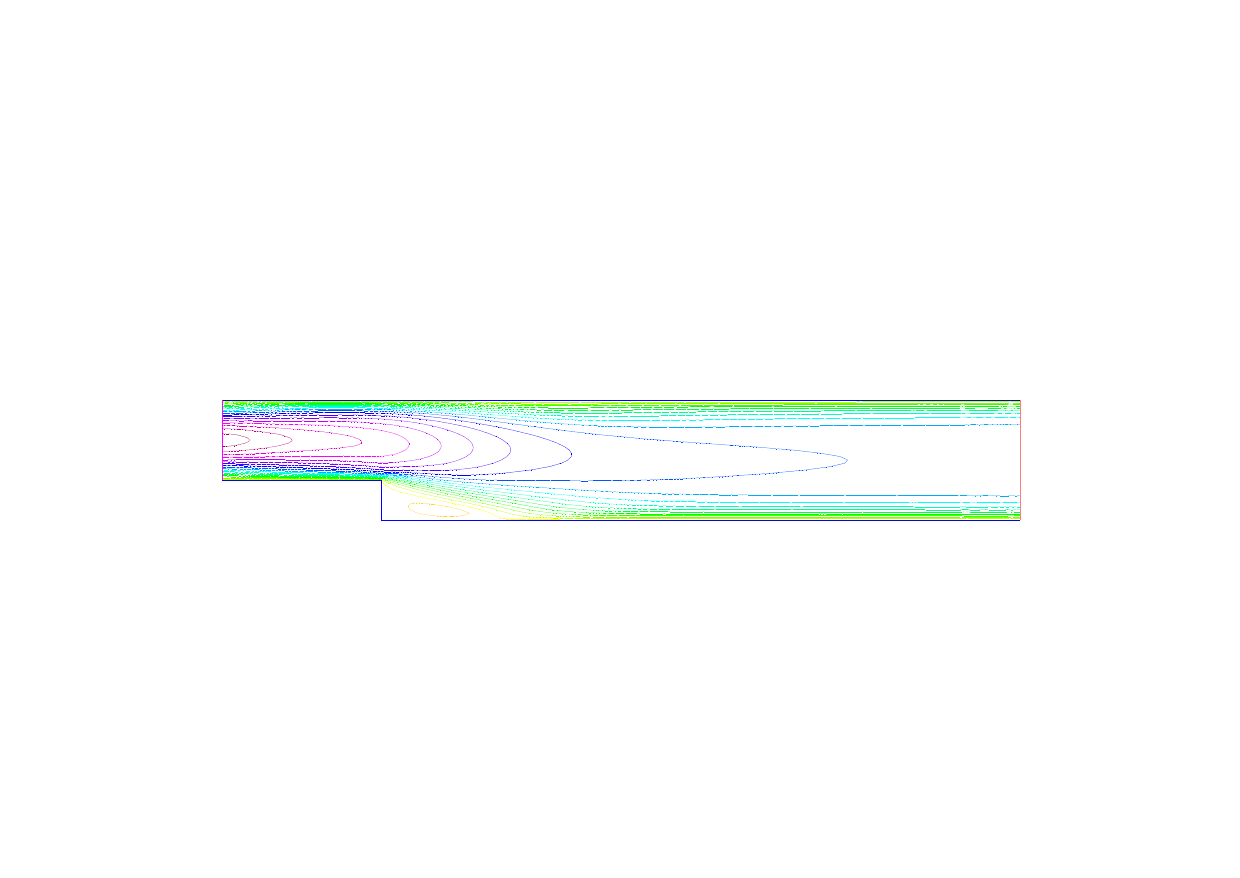}}
\subfigure[  $u_1$: $\alpha=1$]
{\includegraphics[width=4cm]{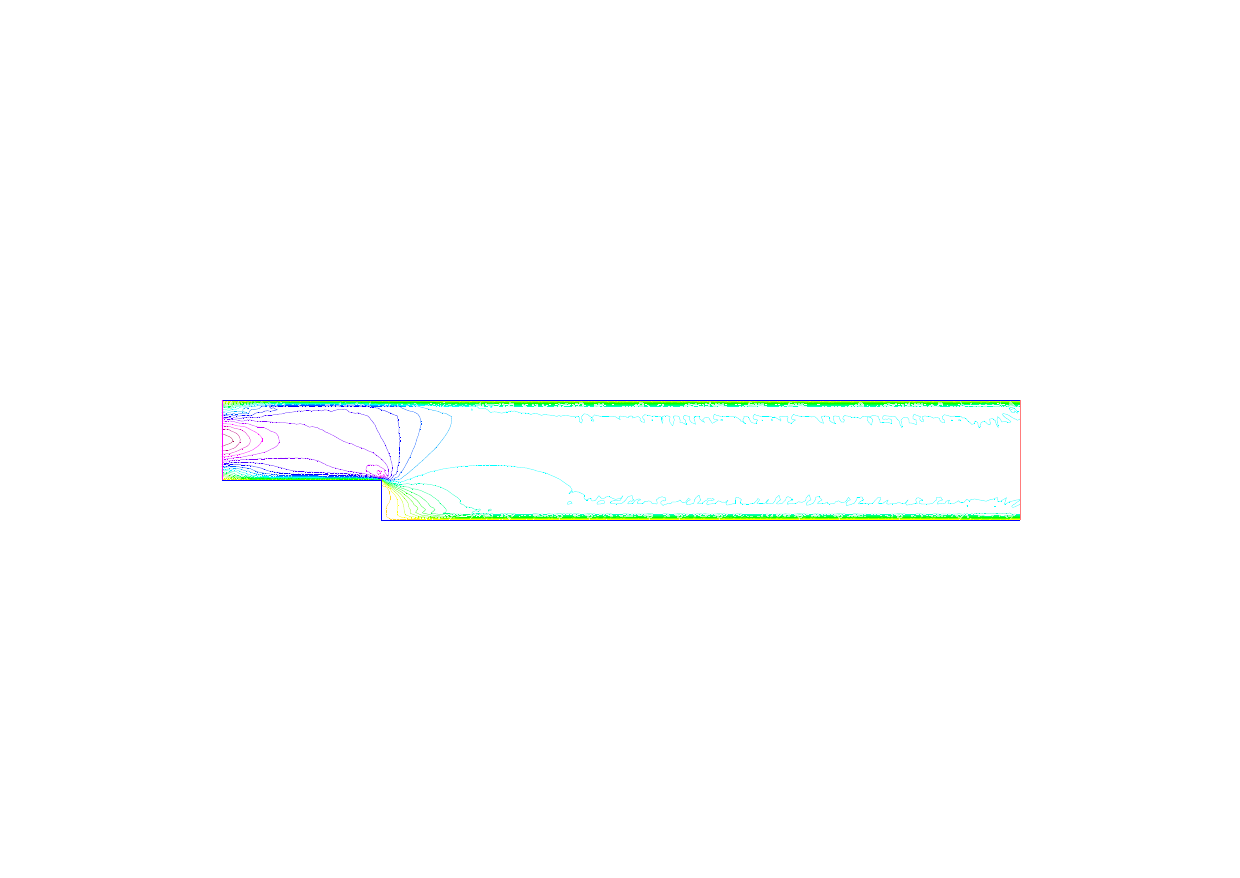}}\\

\subfigure[  $u_2$: $\alpha=0.01$]
{\includegraphics[width=4cm]{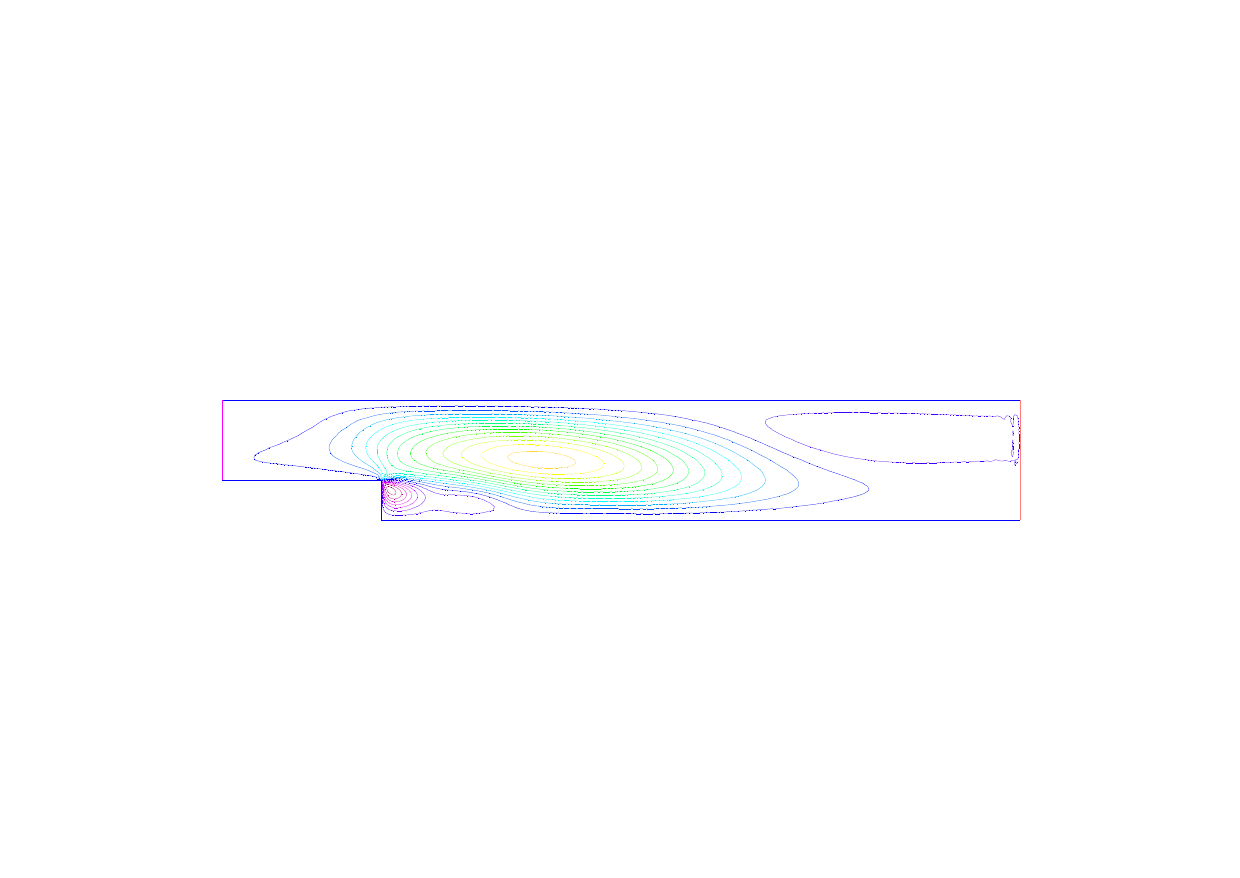}}
\subfigure[  $u_2$: $\alpha=0.1$]
{\includegraphics[width=4cm]{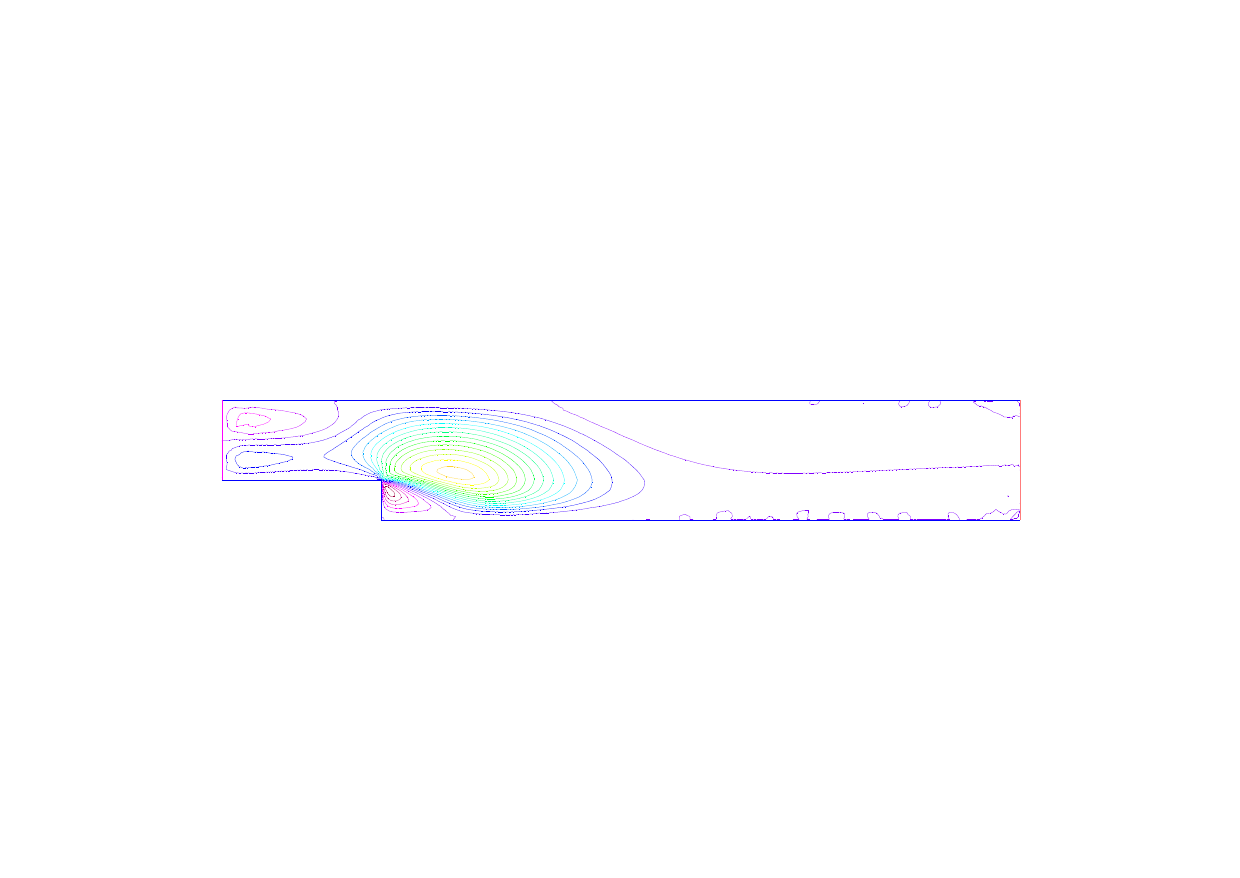}}
\subfigure[  $u_2$: $\alpha=1$]
{\includegraphics[width=4cm]{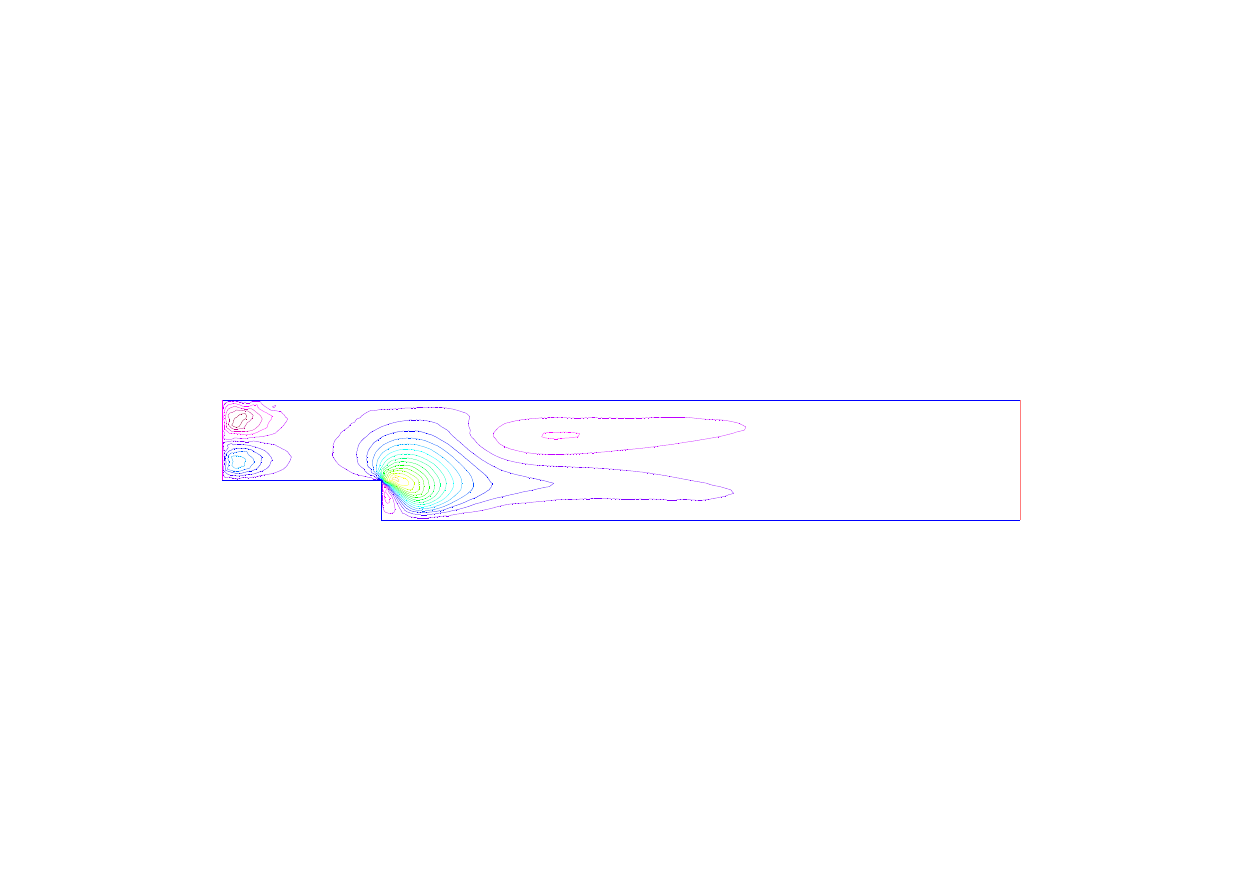}}

\subfigure[pressure: $\alpha=0.01$]
{\includegraphics[width=4cm]{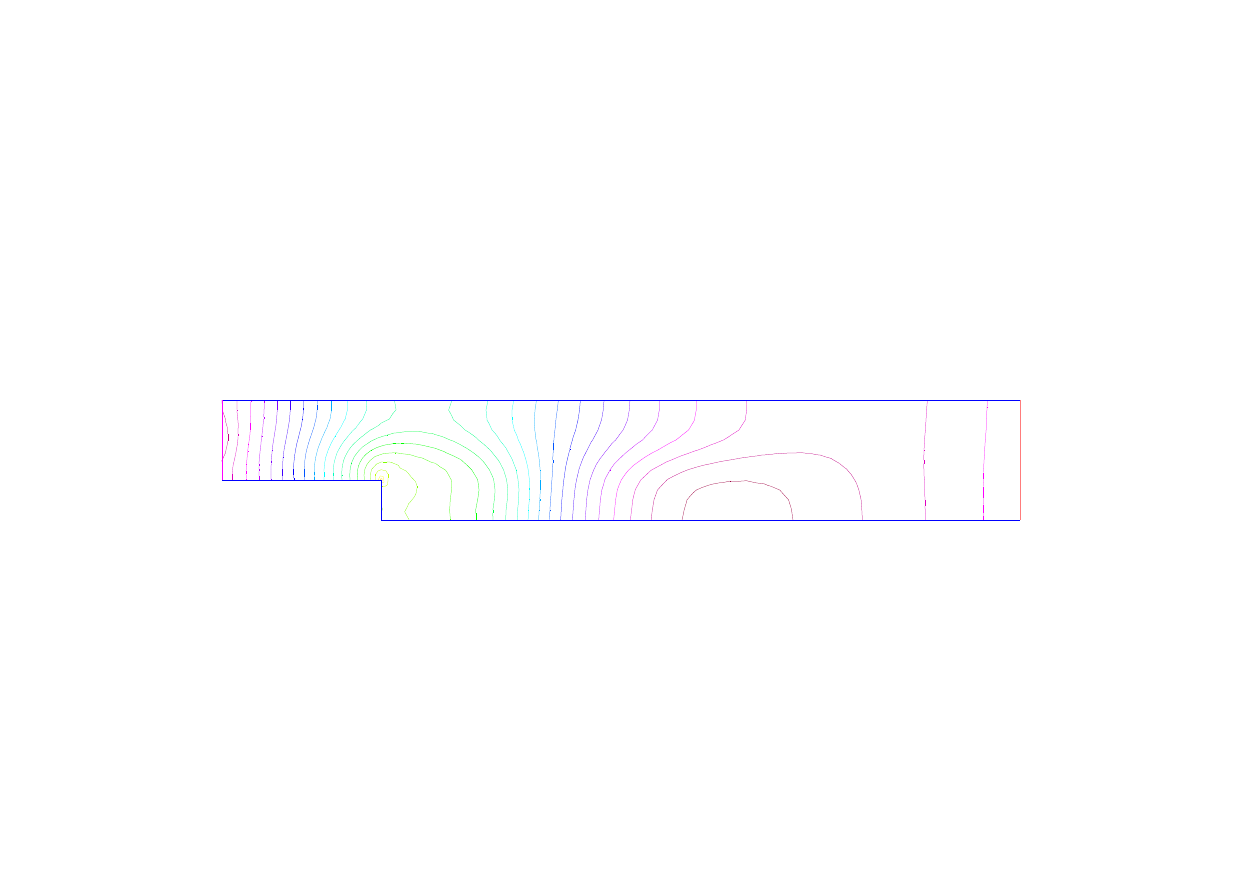}}
\subfigure[pressure: $\alpha=0.1$]
{\includegraphics[width=4cm]{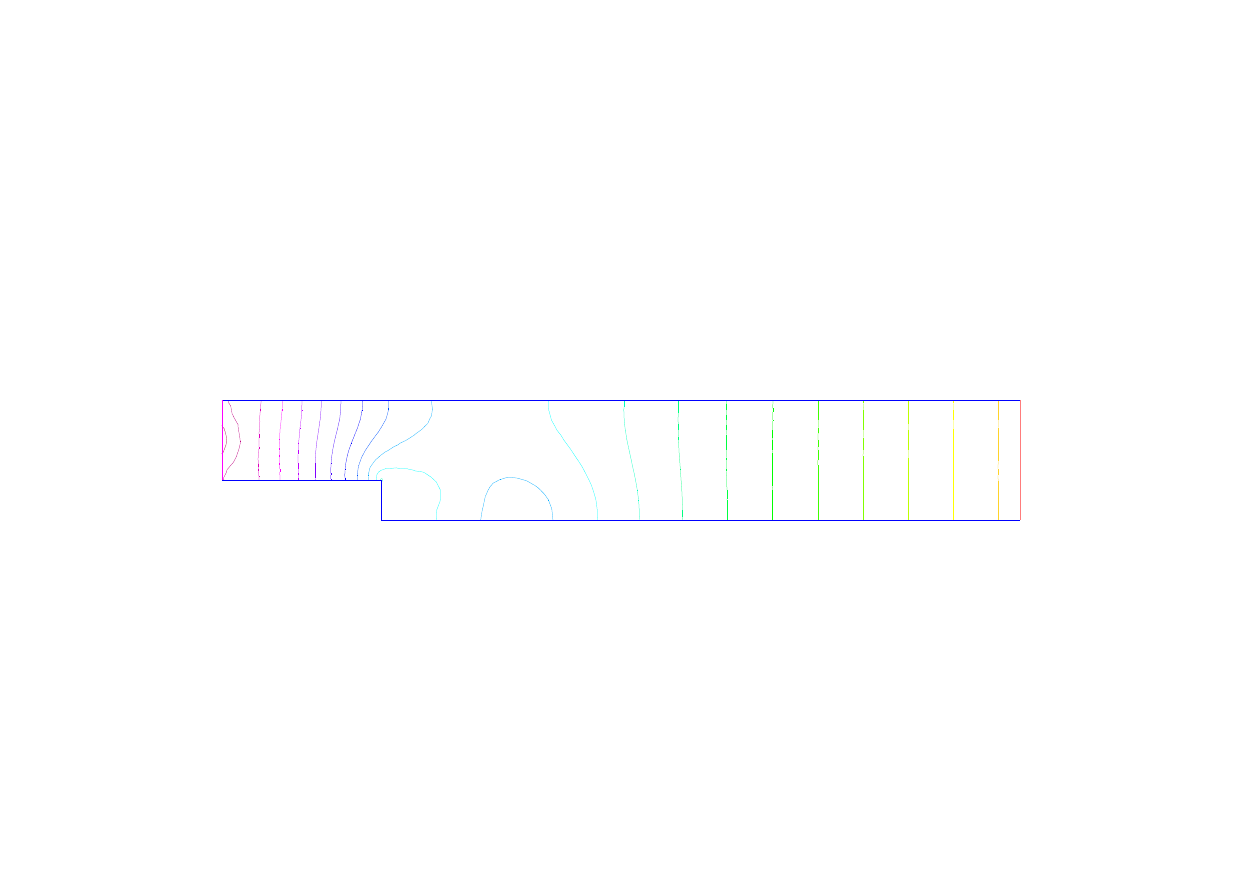}}
\subfigure[pressure: $\alpha=1$]
{\includegraphics[width=4cm]{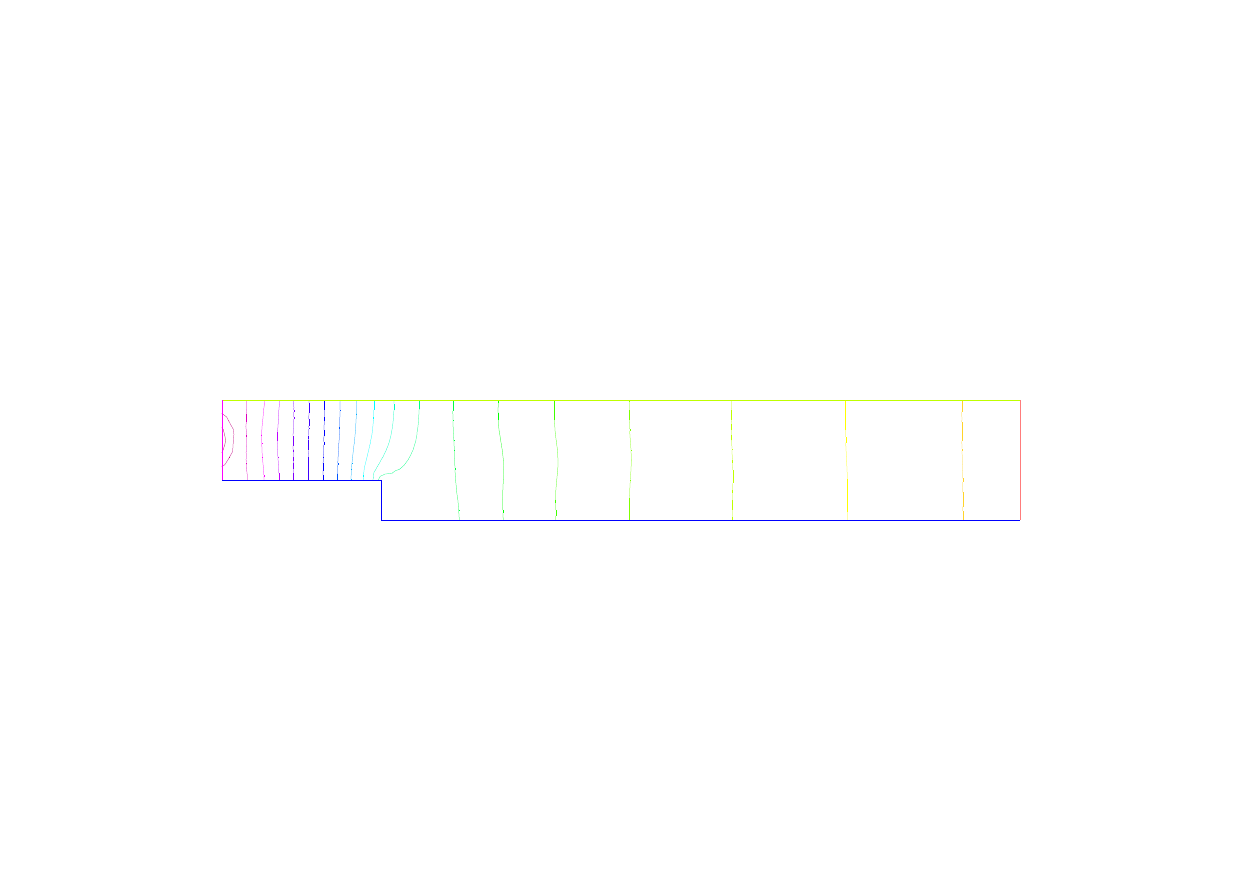}}
\caption{The velocity  $\bm u_h=(u_1,u_2)^T$   and pressure contours for Example \ref{EX7.4}: $r=3.5$ and $\alpha=0.01, 0.1, 1$  }
\label{fig41:41}
\end{figure}

\begin{figure}[htbp!]
\centering
\subfigure[  $u_1$: $r=5$]
{\includegraphics[width=4cm]{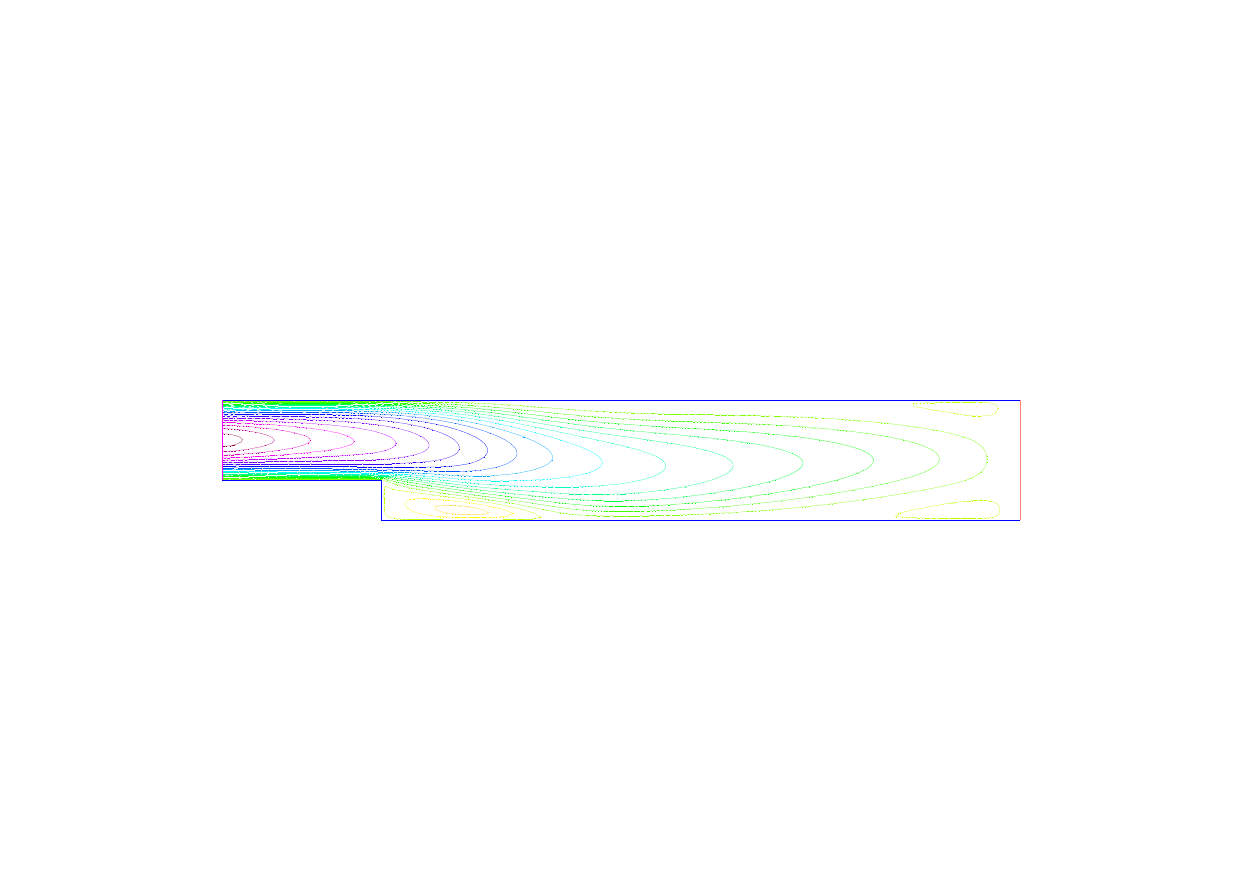}}
\subfigure[ $u_1$: $r=10$]
{\includegraphics[width=4cm]{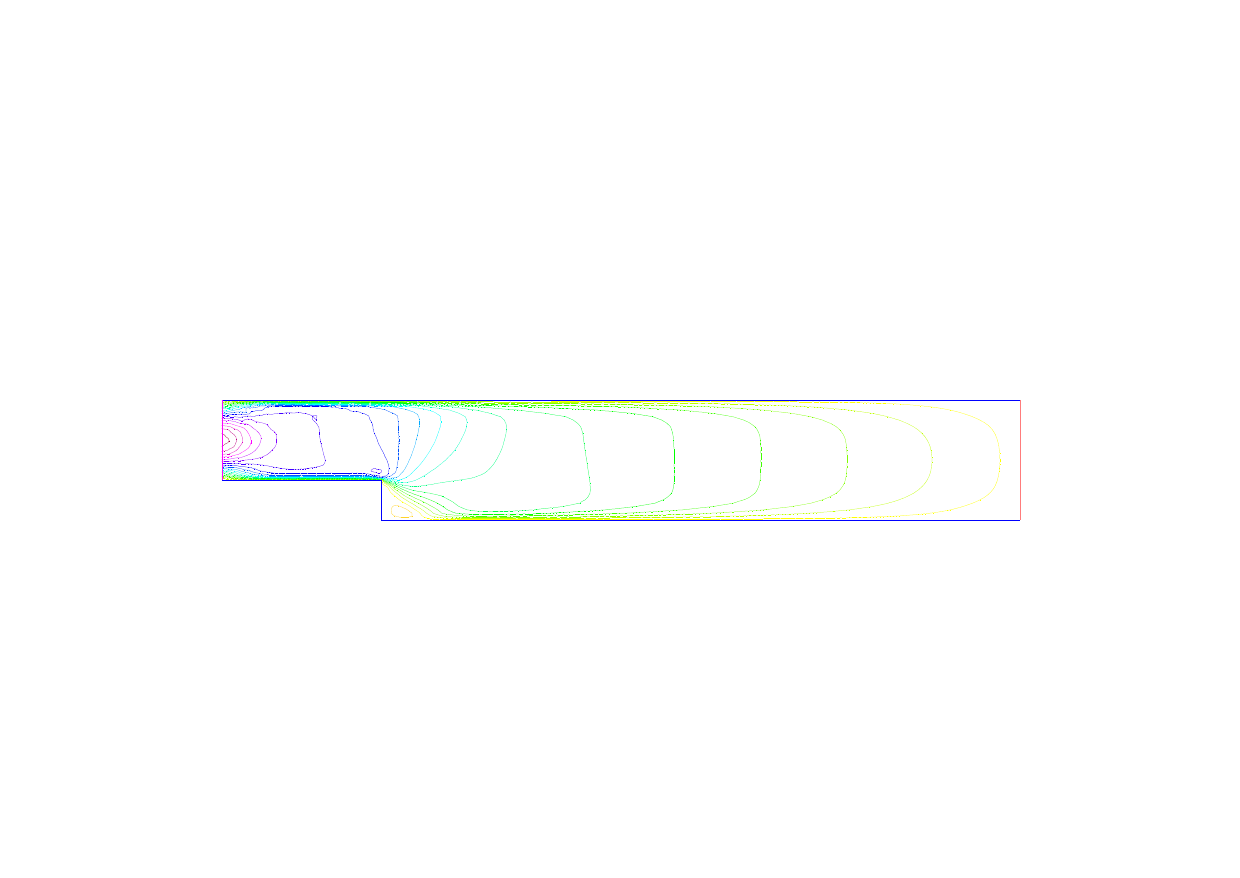}}
\subfigure[  $u_1$: $r=50$]
{\includegraphics[width=4cm]{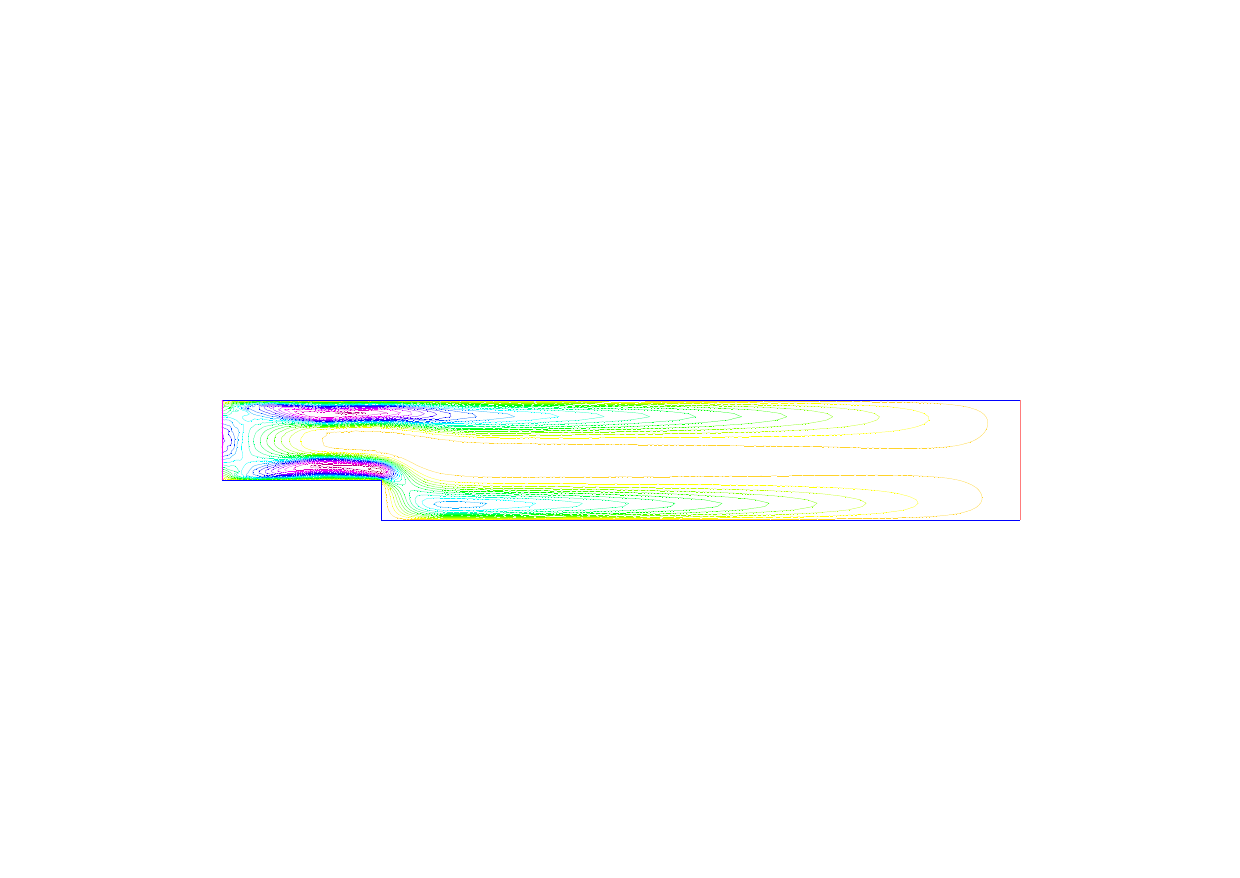}}\\

\subfigure[  $u_2$: $r=5$]
{\includegraphics[width=4cm]{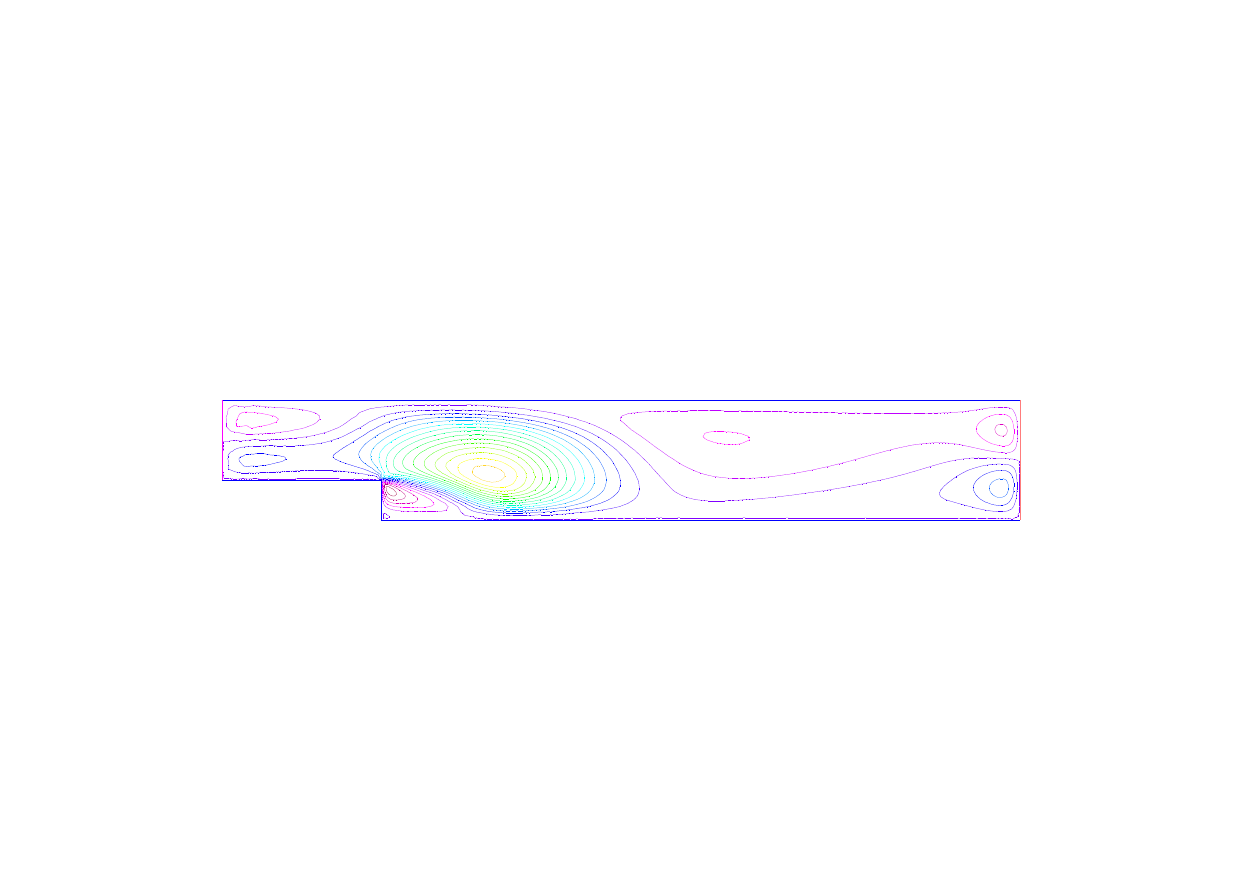}}
\subfigure[ $u_2$: $r=10$]
{\includegraphics[width=4cm]{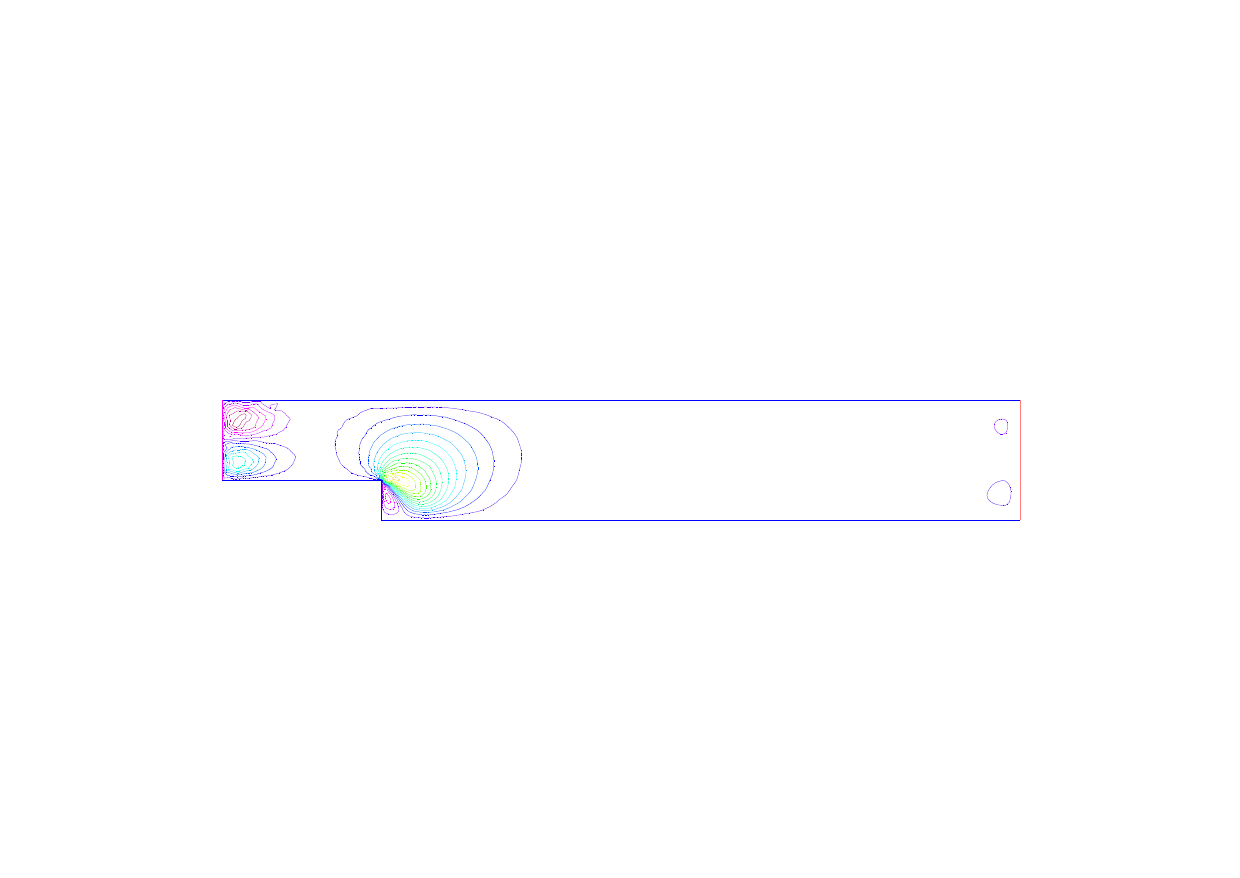}}
\subfigure[  $u_2$: $r=50$]
{\includegraphics[width=4cm]{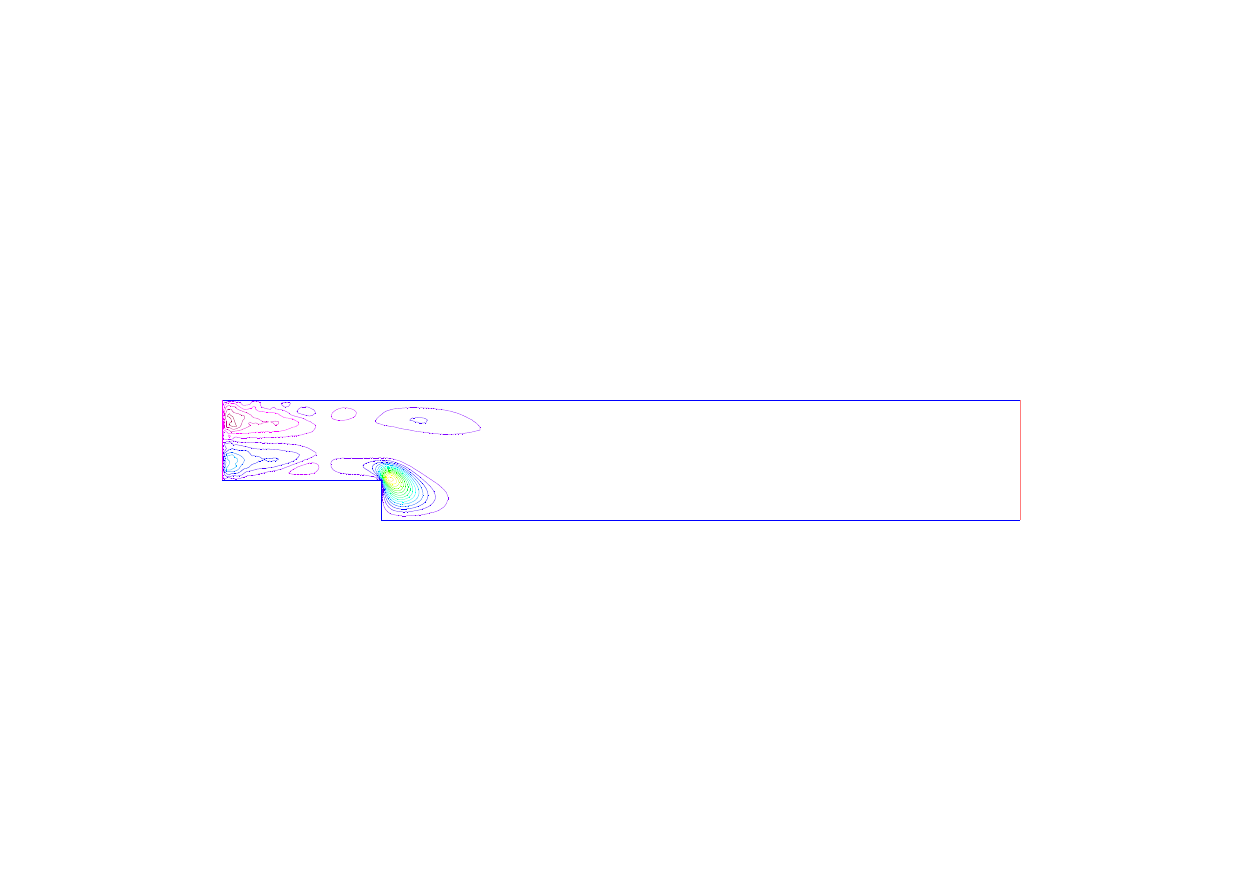}}\\

 \subfigure[pressure: $r=5$]
{\includegraphics[width=4cm]{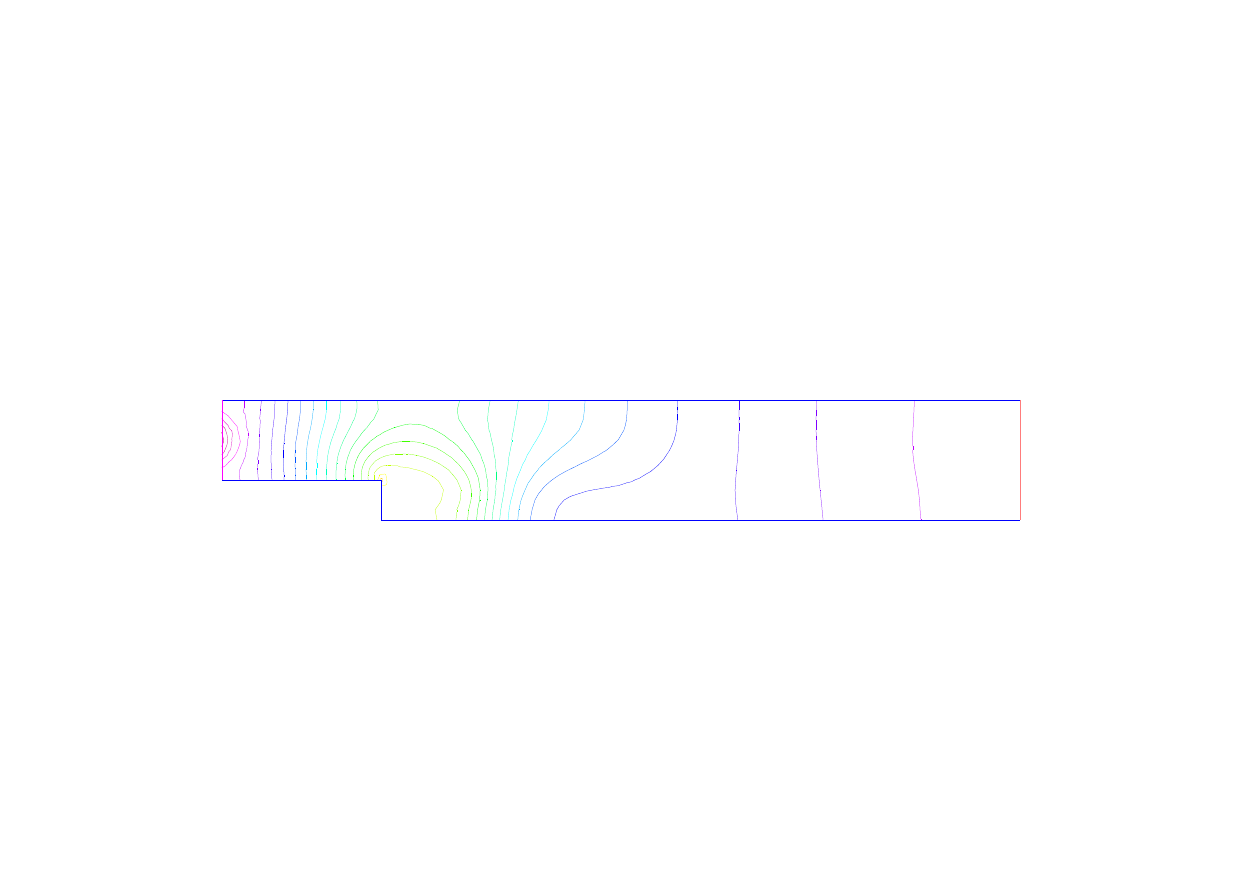}}
\subfigure[pressure: $r=10$]
{\includegraphics[width=4cm]{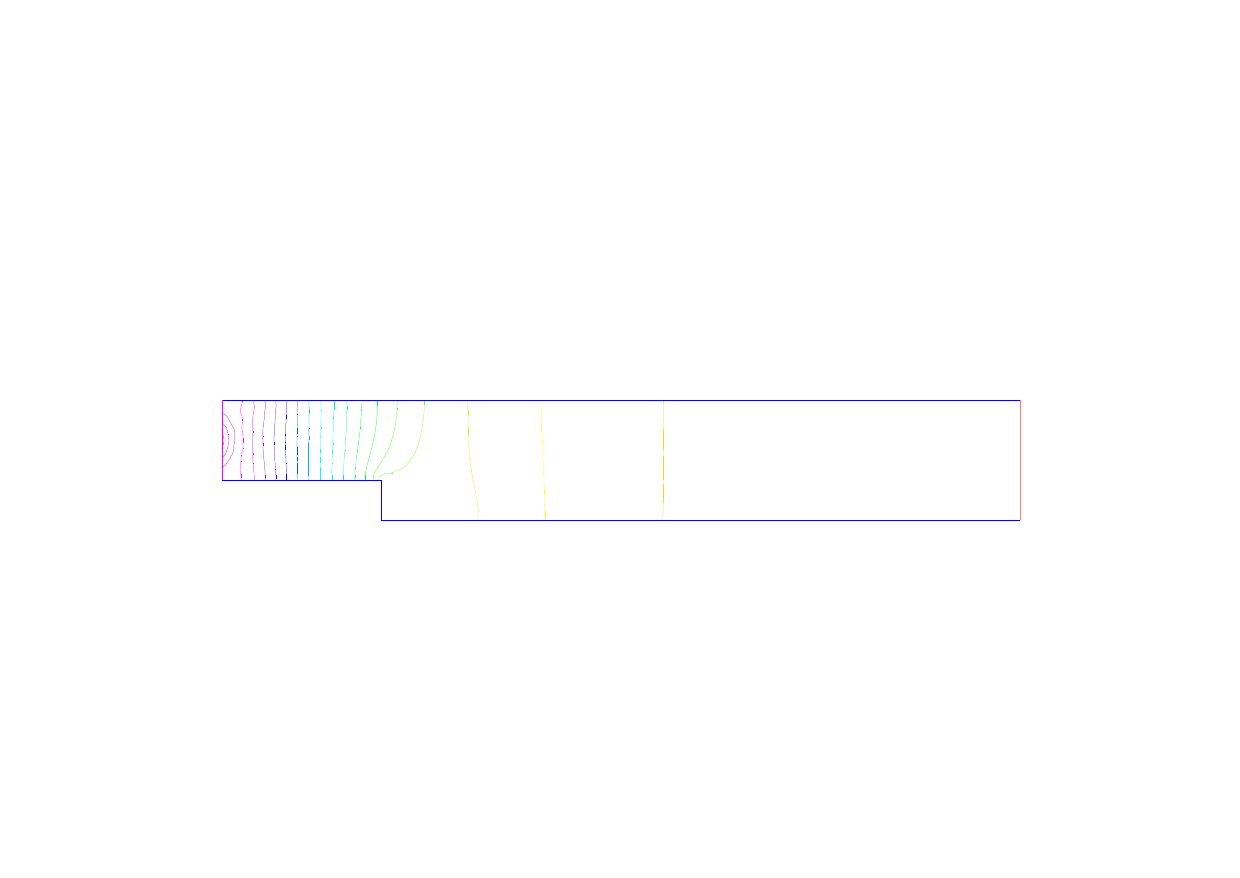}}
\subfigure[pressure: $r=50$]
{\includegraphics[width=4cm]{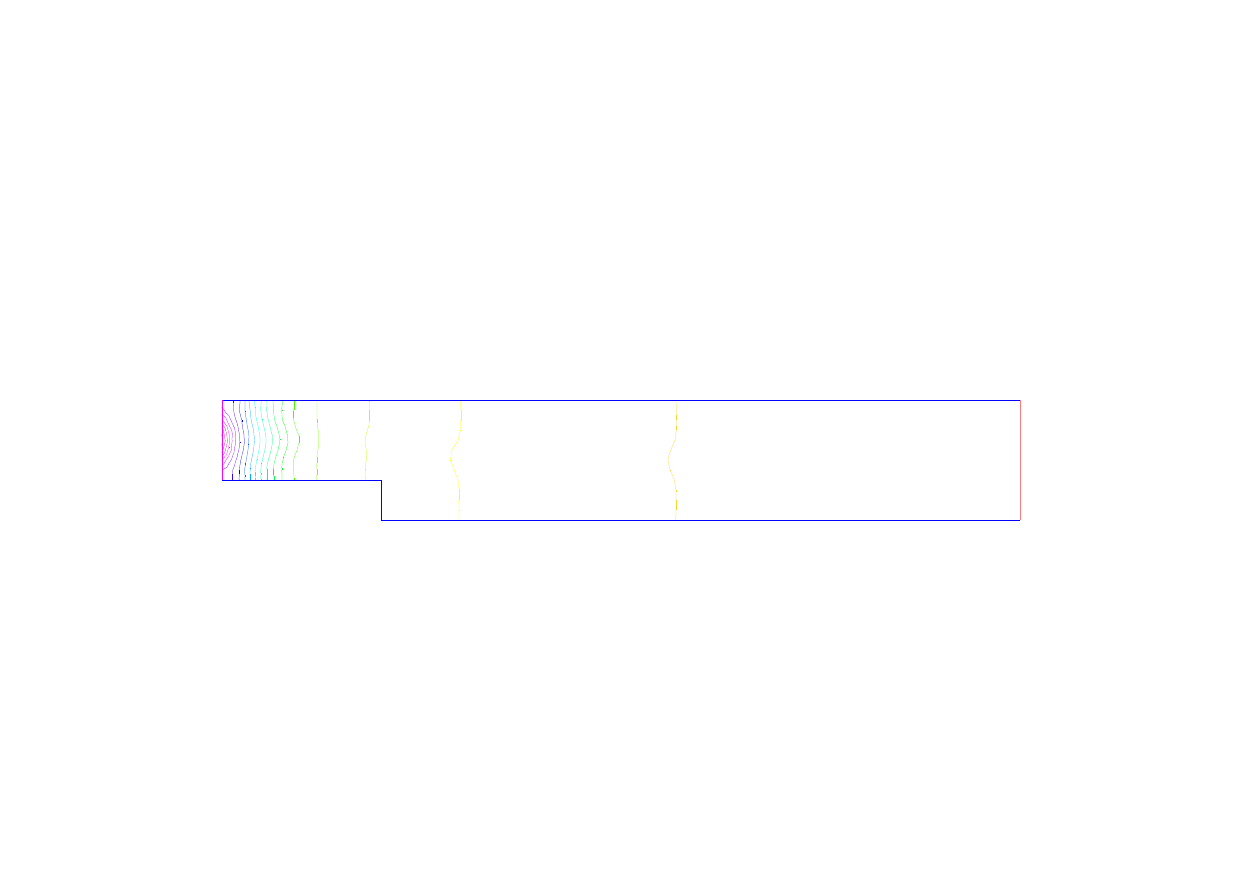}}
\caption{The velocity  $\bm u_h=(u_1,u_2)^T$   and pressure contours for Example \ref{EX7.4}:  $\alpha=1$ and diverse $r=5,10,50$}
\label{fig42:42}
\end{figure}

\begin{exam}[The backward-facing step flow problem]\label{EX7.4}
We consider a backward-facing step flow problem in  $\Omega=\Omega_1\setminus \Omega_2$, with  $\Omega_1=[-4, 16]\times[-1, 2]$ and  $\Omega_2=[-4, 0]\times[-1, 0]$; see Figure \ref{fig4:mesh3} for the   domain and its finite element mesh. We take $\nu=0.005$   and $\bm{f}=\bm{0}$.  The boundary conditions are as follows:
$$\bm{u}|_{y=-1}=\bm{u}|_{y=2}=\bm{u}|_{-4\leq x\leq 0,y=0}=\bm{u}|_{x=0, -1\leq y\leq 0}=\bm{0},  $$
$$\bm{u}|_{x=-4}=(y(2-y), 0 )^T, \quad  \left(-p +\nu\frac{\partial u_1}{\partial x}\right) |_{x=16}=0, \quad u_{2}|_{x=16}=0.$$

We compute the WG scheme \eqref{WG} with  $m=k=2$ in the following cases:
\begin{itemize}
\item [ I]. $\alpha=0$, i.e. the case of the Navier-Stokes  equations;

\item [ II].   $r=3.5$ and $\alpha=0.01, 0.1, 1$;

\item [ III]. $\alpha=1$ and $r=5, 10, 50$.
\end{itemize}
The obtained velocity and pressure approximations are shown   in Figures \ref{fig41:40}, \ref{fig41:41} and \ref{fig42:42}.  As a comparison,  the    numerical solutions  obtained with the Taylor-Hood element 
 are also shown for $\alpha=0$; see (a), (b)  and (c) in Figure \ref{fig41:40}. Similar to Example \ref{EX7.3}, we can see that  our method is effective and the damping effect is gradually enhanced as the parameters $\alpha$ and $r$   increase.

\end{exam}

\section{Conclusion}

We have developed  a class of WG methods of arbitrary order for the steady Brinkman-Forchheimer equations.  The methods  yield   globally divergence-free velocity and  are pressure robust. Optimal error estimates have been derived for  the  velocity and pressure approximations.  The proposed Oseen's iteration algorithm is unconditionally convergent. Numerical experiments have verified the theoretical analysis and demonstrated the robustness of the methods.







\bibliographystyle{siam}
\bibliography{mybib}{}

\begin{thebibliography}{10}

\bibitem{bookin1}
{\sc F.~Auricchio, L.~Beir{\~a}o~da Veiga, C.~Lovadina, and A.~Reali}, {\em The
  importance of the exact satisfaction of the incompressibility constraint in
  nonlinear elasticity: mixed {FEM}s versus {NURBS}-based approximations},
  Computer Methods in Applied Mechanics and Engineering, 199 (2010),
  pp.~314--323.

\bibitem{BorggaardJeff2012Tdot}
{\sc J.~Borggaard, T.~Iliescu, and J.~P. Roop}, {\em {T}wo-level discretization
  of the {N}avier-{S}tokes equations with r-{L}aplacian subgridscale
  viscosity}, Numerical methods for partial differential equations, 28 (3)
  (2012), pp.~1056--1078.

\bibitem{RT-1}
{\sc F.~Brezzi, D.~Boffi, L.~Demkowicz, R.~Durán, R.~Falk, and M.~Fortin},
  {\em Mixed Finite Elements, Compatibility Conditions, and Applications},
  Springer Berlin Heidelberg, 2008.

\bibitem{CaiXiaojing2008Wass}
{\sc X.~Cai and Q.~Jiu}, {\em {W}eak and strong solutions for the
  incompressible {N}avier-{S}tokes equations with damping}, Journal of
  Mathematical Analysis and Applications, 343 (2) (2008), pp.~799--809.

\bibitem{CaucaoSergio2022AtBs}
{\sc S.~Caucao, R.~Oyarzúa, S.~Villa-Fuentes, and I.~Yotov}, {\em {A}
  three-field {B}anach spaces-based mixed formulation for the unsteady
  {B}rinkman-{F}orchheimer equations}, Computer methods in applied mechanics
  and engineering, 394 (2022), p.~114895.

\bibitem{Caucao;2021}
{\sc S.~Caucao and I.~Yotov}, {\em {A} {B}anach space mixed formulation for the
  unsteady {B}rinkman-{F}orchheimer equations}, IMA Journal of Numerical
  Analysis, 41 (4) (2021), pp.~2708--2743.

\bibitem{CelebiA.O.2006Ocdo}
{\sc A.~Celebi, V.~Kalantarov, and D.~Ugurlu}, {\em {O}n continuous dependence
  on coefficients of the {B}rinkman-{F}orchheimer equations}, Applied
  Mathematics Letters, 19 (8) (2006), pp.~801--807.

\bibitem{CFX2016}
{\sc G.~Chen, M.~Feng, and X.~Xie}, {\em {R}obust globally divergence-free weak
  {G}alerkin methods for {S}tokes equations}, Journal of Computational
  Mathematics, 34 (5) (2016), pp.~549--572.

\bibitem{Chen-Xie2023}
{\sc G.~Chen and X.~Xie}, {\em {A}nalysis of a class of globally
  divergence-free {HDG} methods for stationary {N}avier-{S}tokes equations},
  Science China-Mathematics, 66 (2023).
\newblock \url{https://doi.org/10.1007/s11425-022-2077-7}.

\bibitem{CiarletP.G1978TFEM}
{\sc P.~G. Ciarlet}, {\em {T}he finite element method for {E}lliptic problems},
  vol.~v.Volume 4 of Studies in mathematics and its applications, Elsevier
  Science, San Diego, 1~ed., 1978.

\bibitem{CKS2007}
{\sc B.~Cockburn, G.~Kanschat, and D.~Sch\"{o}tzau}, {\em {A} note on
  discontinuous {G}alerkin divergence-free solutions of the {N}avier-{S}tokes
  equations}, Journal of Scientific Computing, 31 (2007), pp.~61--73.

\bibitem{Djoko;2014}
{\sc J.~K. Djoko and P.~A. Razafimandimby}, {\em {A}nalysis of the
  {B}rinkman-{F}orchheimer equations with slip boundary conditions}, Applicable
  Analisis, 93 (2014), pp.~1477--1494.

\bibitem{Girault.V;Raviart.P1986}
{\sc V.~Girault and P.-A. Raviart}, {\em Finite element methods for
  {N}avier-{S}tokes equations}, vol.~5 of Springer Series in Computational
  Mathematics, Springer-Verlag, Berlin, 1986.
\newblock Theory and algorithms.

\bibitem{HLX2019}
{\sc Y.~Han, H.~Li, and X.~Xie}, {\em {R}obust globally divergence-free weak
  {G}alerkin finite element methods for unsteady natural convection problems},
  Numerical Mathematics: Theory, Methods and Applications, 12 (2019),
  pp.~1266--1308.

\bibitem{HX2019}
{\sc Y.~Han and X.~Xie}, {\em {R}obust globally divergence-free weak {G}alerkin
  finite element methods for natural convection problems}, Communications in
  Computational Physics, 26 (2019), pp.~1039--1070.

\bibitem{HuXiaozhe2019AwGf}
{\sc X.~Hu, L.~Mu, and X.~Ye}, {\em {A} weak {G}alerkin finite element method
  for the {N}avier-{S}tokes equations}, Journal of computational and applied
  mathematics, 362 (2019), pp.~614--625.

\bibitem{JiangZaihong2012Abos}
{\sc Z.~Jiang}, {\em {A}symptotic behavior of strong solutions to the 3{D}
  {N}avier-{S}tokes equations with a nonlinear damping term}, Nonlinear
  analysis, 75 (13) (2012), pp.~5002--5009.

\bibitem{JLMNR2017}
{\sc V.~John, A.~Linke, C.~Merdon, M.~Neilan, and L.~Rebholz}, {\em {O}n the
  divergence constraint in mixed finite element methods for incompressible
  flows}, SIAM Review, 59 (3) (2017), pp.~492--544.

\bibitem{D1982Nonlinear}
{\sc D.~D. Joseph, D.~A. Nield, and G.~Papanicolaou}, {\em {N}onlinear equation
  governing flow in a saturated porous medium}, Water Resources Research, 18
  (4) (1982), pp.~1049--1052.

\bibitem{MR3511719}
{\sc C.~Lehrenfeld and J.~Sch\"{o}berl}, {\em High order exactly
  divergence-free hybrid discontinuous {G}alerkin methods for unsteady
  incompressible flows}, Computer Methods in Applied Mechanics and Engineering,
  307 (2016), pp.~339--361.

\bibitem{LiKwang-Ok2022Grft}
{\sc K.-O. Li and Y.-H. Kim}, {\em {G}lobal regularity for the 3{D}
  inhomogeneous incompressible {N}avier-{S}tokes equations with damping},
  Applications of mathematics (Prague), 68 (2) (2022), pp.~191--207.

\bibitem{2022Unconditional}
{\sc M.~Li, Z.~Li, and D.~Shi}, {\em Unconditional optimal error estimates for
  the transient {N}avier-{S}tokes equations with damping}, Advances in Applied
  Mathematics and Mechanics,  (2022), p.~14.

\bibitem{LiMinghao2019Tmfe}
{\sc M.~Li, D.~Shi, Z.~Li, and H.~Chen}, {\em {T}wo-level mixed finite element
  methods for the {N}avier-{S}tokes equations with damping}, Journal of
  mathematical analysis and applications, 470 (1) (2019), pp.~292--307.

\bibitem{LiYuanfei2014Cdft}
{\sc Y.~Li and C.~Lin}, {\em {C}ontinuous dependence for the nonhomogeneous
  {B}rinkman-{F}orchheimer equations in a semi-infinite pipe}, Applied
  mathematics and computation, 244 (2014), pp.~201--208.

\bibitem{LiZhenzhen2019Smfe}
{\sc Z.~Li, D.~Shi, and M.~Li}, {\em {S}tabilized mixed finite element methods
  for the {N}avier-{S}tokes equations with damping}, Mathematical methods in
  the applied sciences, 42 (2) (2019), pp.~605--619.

\bibitem{add4}
{\sc A.~Linke}, {\em {D}ivergence-free mixed finite elements for the
  incompressible a {N}avier-{S}tokes equation}, Ph.D. Dissertation, University
  of Erlangen,  (2007).

\bibitem{bookin3}
{\sc A.~Linke}, {\em Collision in a cross-shaped domain---a steady 2d
  {N}avier-{S}tokes example demonstrating the importance of mass conservation
  in {CFD}}, Computer Methods in Applied Mechanics and Engineering, 198 (2009),
  pp.~3278--3286.

\bibitem{MR3564690}
{\sc A.~Linke and C.~Merdon}, {\em Pressure-robustness and discrete {H}elmholtz
  projectors in mixed finite element methods for the incompressible
  {N}avier-{S}tokes equations}, Computer Methods in Applied Mechanics and
  Engineering, 311 (2016), pp.~304--326.

\bibitem{2019Mixed}
{\sc D.~Liu and K.~Li}, {\em {M}ixed finite element for two-dimensional
  incompressible convective {B}rinkman-{F}orchheimer equations}, Applied
  Mathematics and Mechanics(English Edition), 40 (6) (2019), pp.~889--910.

\bibitem{LiuHui2021Wotg}
{\sc H.~Liu, L.~Lin, and C.~Sun}, {\em {W}ell-posedness of the generalized
  {N}avier-{S}tokes equations with damping}, Applied mathematics letters, 121
  (2021), p.~107471.

\bibitem{QianLiu2021Saoa}
{\sc Q.~Liu and D.~Shi}, {\em {S}uperconvergent analysis of a nonconforming
  mixed finite element method for time-dependent damped {N}avier-{S}tokes
  equations}, Computational and applied mathematics, 40 (1) (2021).

\bibitem{Louaked;2017}
{\sc M.~Louaked, N.~Seloula, and S.~Trabelsi}, {\em {A}pproximation of the
  unsteady {B}rinkman-{F}orchheimer equations by the pressure stabilization
  method}, Numerical Methods for Partial Differential Equations, 33 (6) (2017),
  pp.~1949--1965.

\bibitem{MuLin2023ApwG}
{\sc L.~Mu}, {\em {A} pressure-robust weak {G}alerkin finite element method for
  {N}avier-{S}tokes equations}, Numerical methods for partial differential
  equations, 39 (3) (2023), pp.~2327--2354.

\bibitem{MuLin2018Addf}
{\sc L.~Mu, J.~Wang, X.~Ye, and S.~Zhang}, {\em {A} discrete divergence free
  weak {G}alerkin finite element method for the {S}tokes equations}, Applied
  numerical mathematics, 125 (2018), pp.~172--182.

\bibitem{add2}
{\sc M.~A. Olshanskii and A.~Reusken}, {\em Grad-div stabilization for {S}tokes
  equations}, Mathematics of Computation, 73 (2004), pp.~1699--1718.

\bibitem{PayneL.E.1999CaCD}
{\sc L.~E. Payne and B.~Straughan}, {\em {C}onvergence and continuous
  dependence for the {B}rinkman-{F}orchheimer equations}, Studies in applied
  mathematics (Cambridge), 102 (4) (1999), pp.~419--439.

\bibitem{PengHui2022WGmf}
{\sc H.~Peng and Q.~Zhai}, {\em {W}eak {G}alerkin method for the {S}tokes
  equations with damping}, Discrete and continuous dynamical systems. Series B,
  27 (4) (2022), p.~1853.

\bibitem{QiuHailong2019Msaf}
{\sc H.~Qiu and L.~Mei}, {\em {M}ulti-level stabilized algorithms for the
  stationary incompressible {N}avier-{S}tokes equations with damping}, Applied
  numerical mathematics, 143 (2019), pp.~188--202.

\bibitem{bookFEM}
{\sc Z.~Shi and M.~Wang}, {\em {F}inite element methods}, Series in Information
  and Computation Science, Science Press, Beijing, 58 (2013).

\bibitem{K1981Boundary}
{\sc K.~Vafai and C.~L. Tien}, {\em {B}oundary and inertia effects on flow and
  heat transfer in porous media}, International Journal of Heat and Mass
  Transfer, 24 (2) (1981), pp.~195--203.

\bibitem{WangJunping2013AwGf}
{\sc J.~Wang and X.~Ye}, {\em {A} weak {G}alerkin finite element method for
  second-order elliptic problems}, Journal of Computational and Applied
  Mathematics, 241 (1) (2013), pp.~103--115.

\bibitem{WangJunping2013AWGM}
{\sc J.~Wang and X.~Ye}, {\em {A} weak {G}alerkin mixed finite element method
  for second-order elliptic problems}, Mathematics of Computation,  (2014),
  pp.~2101--2126.

\bibitem{WangJunping2016AwGf}
{\sc J.~Wang and X.~Ye}, {\em {A} weak {G}alerkin finite element method for the
  {S}tokes equations}, Advances in computational mathematics, 42 (1) (2016),
  pp.~155--174.

\bibitem{WangRuishu2016AwGf}
{\sc R.~Wang, X.~Wang, Q.~Zhai, and R.~Zhang}, {\em {A} weak {G}alerkin finite
  element scheme for solving the stationary {S}tokes equations}, Journal of
  computational and applied mathematics, 302 (2016), pp.~171--185.

\bibitem{WassimEid2022Lapf}
{\sc E.~Wassim and Y.~Shang}, {\em {L}ocal and parallel finite element
  algorithms for the incompressible {N}avier-{S}tokes equations with damping},
  Discrete and continuous dynamical systems. Series B, 27 (11) (2022), p.~6823.

\bibitem{WassimEid2023Aptm}
{\sc E.~Wassim, B.~Zheng, and Y.~Shang}, {\em {A} parallel two-grid method
  based on finite element approximations for the 2{D}/3{D} {N}avier-{S}tokes
  equations with damping}, Engineering with computers,  (2023).
\newblock \url{https://doi.org/10.1007/s00366-023-01807-w}.

\bibitem{XZ2010}
{\sc X.~Xu and S.~Zhang}, {\em {A} new divergence-free interpolation operator
  with applications to the {D}arcy-{S}tokes-{B}rinkman equations}, SIAM Journal
  On Scientific Computing, 32 (2) (2010), pp.~855--874.

\bibitem{YangHuaijun2023Saot}
{\sc H.~Yang and X.~Jia}, {\em {S}uperconvergence analysis of the
  bilinear-constant scheme for two-dimensional incompressible convective
  {B}rinkman-{F}orchheimer equations}, Numerical methods for partial
  differential equations, 40 (1) (2024).

\bibitem{YangRong2021Aso3}
{\sc R.~Yang and X.-G. Yang}, {\em {A}symptotic stability of 3{D}
  {N}avier-{S}tokes equations with damping}, Applied mathematics letters, 116
  (2021), p.~107012.

\bibitem{ZhangLi2019Agdw}
{\sc L.~Zhang, M.~Feng, and J.~Zhang}, {\em {A} globally divergence-free weak
  {G}alerkin method for {B}rinkman equations}, Applied numerical mathematics,
  137 (2019), pp.~213--229.

\bibitem{ZZX2023}
{\sc M.~Zhang, T.~Zhang, and X.~Xie}, {\em {R}obust globally divergence-free
  weak {G}alerkin finite element method for incompressible
  {M}agnetohydrodynamics flow}, Communications in Nonlinear Science and
  Numerical Simulation, 131 (2024), p.~107810.

\bibitem{ZhangZujin2011Otuo}
{\sc Z.~Zhang, X.~Wu, and M.~Lu}, {\em {O}n the uniqueness of strong solution
  to the incompressible {N}avier-{S}tokes equations with damping}, Journal of
  Mathematical Analysis and Applications, 377 (1) (2011), pp.~414--419.

\bibitem{ZhengBo2021Tdsa}
{\sc B.~Zheng and Y.~Shang}, {\em {T}wo-level defect-correction stabilized
  algorithms for the simulation of 2{D}/3{D} steady {N}avier-{S}tokes equations
  with damping}, Applied numerical mathematics, 163 (2021), pp.~182--203.

\bibitem{ZhengBo2023Tsaf}
{\sc B.~Zheng and Y.~Shang}, {\em {T}wo-grid stabilized algorithms for the
  steady {N}avier-{S}tokes equations with damping}, Mathematical methods in the
  applied sciences, 46 (1) (2023), pp.~107--125.

\bibitem{ZCX2017}
{\sc X.~Zheng, G.~Chen, and X.~Xie}, {\em {A} divergence-free weak {G}alerkin
  method for quasi-{N}ewtonian {S}tokes flows}, Science China Mathematics, 60
  (8) (2017), pp.~1515--1528.

\bibitem{ZhongXin2019Anot}
{\sc X.~Zhong}, {\em {A} note on the uniqueness of strong solution to the
  incompressible {N}avier-{S}tokes equations with damping}, Electronic Journal
  of Qualitative Theory of Differential Equations, 2019 (15) (2019), pp.~1--4.

\bibitem{ZhouYong2012Rauf}
{\sc Y.~Zhou}, {\em {R}egularity and uniqueness for the 3{D} incompressible
  {N}avier-{S}tokes equations with damping}, Applied Mathematics Letters, 25
  (11) (2012), pp.~1822--1825.

\end{thebibliography}



\end{document}